\documentclass{amsart}

\usepackage{amsmath,amssymb,amsthm,texdraw,citesort}

\numberwithin{equation}{section}

\newsavebox{\tmpfiga}
\newsavebox{\tmpfigb}
\newsavebox{\tmpfigc}
\newsavebox{\tmpfigd}
\newsavebox{\tmpfige}
\newsavebox{\tmpfigf}
\newsavebox{\tmpfigg}
\newsavebox{\tmpfigh}
\newsavebox{\tmpfigi}
\newsavebox{\tmpfigj}
\newsavebox{\tmpfigk}
\newsavebox{\tmpfigl}
\newsavebox{\tmpfigm}
\newsavebox{\tmpfign}
\newsavebox{\tmpfigo}
\newsavebox{\tmpfigp}
\newsavebox{\tmpfigq}

\theoremstyle{plain}
\newtheorem{thm}{Theorem}[section]
\newtheorem{prop}[thm]{Proposition}

\theoremstyle{definition}
\newtheorem{df}[thm]{Definition}
\newtheorem{example}[thm]{Example}
\theoremstyle{remark}
\newtheorem{remark}[thm]{Remark}

\newcommand{\nc}{\newcommand}

\nc{\eit}{\tilde{e}_i}
\nc{\fit}{\tilde{f}_i}
\nc{\ali}{\alpha_i}
\nc{\csa}{\mathfrak{h}}
\nc{\op}{\oplus}
\nc{\ot}{\otimes}
\nc{\C}{\mathbf{C}}
\nc{\Z}{\mathbf{Z}}
\nc{\g}{\mathfrak{g}}
\nc{\slice}{\mathcal{s}^{(l)}}
\nc{\nslice}{\mathcal{S}^{(l)}}
\nc{\newcrystal}{\mathcal{C}^{(l)}}
\nc{\uq}{U_q} \nc{\uqp}{U_q'}
\nc{\ftil}{\tilde{f}}
\nc{\etil}{\tilde{e}}
\nc{\anone}{A_n^{(1)}}
\nc{\cnone}{C_n^{(1)}}
\nc{\dnone}{D_n^{(1)}}
\nc{\atwontwo}{A_{2n}^{(2)}}
\nc{\dnponetwo}{D_{n+1}^{(2)}}
\nc{\bnone}{B_n^{(1)}}
\nc{\atwonmonetwo}{A_{2n-1}^{(2)}}
\nc{\oldcrystal}{\mathcal{B}^{(l)}}
\nc{\defi}[1]{\emph{\textbf{#1}}}
\nc{\veps}{\varepsilon}
\nc{\vphi}{\varphi}
\nc{\La}{\Lambda}
\nc{\la}{\lambda}
\nc{\pathspace}{\mathcal{P}}
\nc{\pwspace}{\mathcal{F}}
\nc{\rpwspace}{\mathcal{Y}}
\nc{\gspath}{\mathbf{b}}
\nc{\path}{\mathbf{p}}
\nc{\spaceiso}{\Psi}
\nc{\spaceosi}{\Phi}
\nc{\gswall}{\mathbf{Y}}
\nc{\wall}{\mathbf{Y}}
\nc{\hwc}{\mathcal{B}}
\nc{\ba}{\mathbf{a}}
\nc{\bb}{\mathbf{b}}
\nc{\bc}{\mathbf{c}}
\nc{\bd}{\mathbf{d}}
\nc{\zo}{0\!|\!1}
\nc{\nn}{(\!n\text{-}\!1\!)\!|\!n}
\nc{\tf}{3\!|\!4}
\DeclareMathOperator{\wt}{wt}
\DeclareMathOperator{\cwt}{\overline{wt}}

\begin{document}

\title[Affine crystals of type $\dnone$ and Young walls]
{Affine crystals of type $\dnone$ and Young walls}

\author{Hyeonmi Lee}
\thanks{This work was supported in part by KOSEF Grant
        R01-2003-000-10012-0}
\address{Korea Institute for Advanced Study\\
         207-43 Cheongnyangni 2-dong, Dongdaemun-gu\\
         Seoul 130-722, Korea}
\email{hmlee@kias.re.kr}

\begin{abstract}
We give a new realization of arbitrary level perfect crystals and
arbitrary level irreducible highest weight crystals of type $\dnone$,
in the language of Young walls.
The notions of \emph{splitting of blocks}
and \emph{slices} play crucial roles
in the construction of crystals. The perfect crystals are realized
as the set of equivalence classes of slices.
And the irreducible highest weight crystals
are realized as the affine crystals consisting of
reduced proper
Young walls which, in turn, is a concatenation of slices.
\end{abstract}

\maketitle

%%%%%%%%%%%%%%%%%%%%%%%%%%%%%%%%%%%%%%%%%%%%%%%%%%%%%%%%%%%%%%%%%%%%%%%%%%
\section{Introduction}%
\label{sec:1}
The theory of Young walls is a combinatorial scheme for realizing
crystal bases for quantum affine algebras.
One may view this as a development of the Young tableaux theory,
which was used in realizing irreducible highest weight crystal for
the classical finite types~\cite{MR95c:17025}, in the affine direction.
At its root is the realization of crystal bases for irreducible
highest weight representations over the affine quantum group
of type $\anone$ which appeared in~\cite{MR91j:17021,MR93a:17015}.
Motivated by some insight from solvable lattice model theory,
their theory developed into the theory of \emph{perfect
crystals}~\cite{MR94j:17013,MR94a:17008,MR93g:17027}.
It gives a realization of crystal bases for irreducible highest
weight modules over any affine quantum groups in terms of \emph{paths}.
This path realization plays an important part in the theory
of Young walls.

The concept of Young walls first appeared in~\cite{MR2001j:17029}
for a very restricted case.
It was subsequently extended to cover level-$1$ highest weight
crystals of most classical affine types in~\cite{MR1881971}
and~\cite{Kang2000}.
The $\cnone$ case, which is missing from these works,
was more difficult to deal with,
because the level-$1$ perfect crystal
for this type is intrinsically of level-$2$.
The difficulty was resolved in \cite{CnOne} by introducing the
notions of \emph{slices} and \emph{splitting of blocks}.
It was also in this work that a realization of perfect crystals,
obtained in the process of Young wall realization,
was explicitly stated as a separate result.
This opened up the road for realization of arbitrary level
highest weight crystals and this was done in~\cite{kang:lee}
for most of the classical quantum affine algebras.
This time, type $\dnone$ withstood realization efforts.
The present paper addresses this gap.

In passing, we acknowledge the existence of a
realization for type $\anone$ arbitrary level highest
weight crystals given in~\cite{Savage}, which is very
similar to that appearing in~\cite{kang:lee}.

Let us explain the contents of this paper.
We first introduce the notions of splitting of blocks and slices,
and give a new realization of
level-$l$ perfect crystals as the equivalence classes of slices.
The two notions presented in this work
are a refined versions of those appearing
in previous works \cite{CnOne} and \cite{kang:lee}.
We then proceed to define the notion of Young walls,
proper Young walls, reduced proper Young walls, and
ground-state walls.
Finally, we prove that an arbitrary level irreducible
highest weight crystal may be realized as the affine crystal
consisting of reduced proper Young walls.

One immediate application of this work could be the extension of
the work by Savage\cite{Savage} on the correspondence between
geometric and combinatorial realizations of
crystal graphs to the $\dnone$ case.

%%%%%%%%%%%%%%%%%%%%%%%%%%%%%%%%%%%%%%%%%%%%%%%%%%%%%%%%%%%%%%%%%%%%%%%%%%
\section{Quantum group of type $\dnone$ and its perfect crystal}%
\label{sec:2}

In this section, we introduce notations and cite facts
that are crucial for our work.
We follow notations of the references cited
in the introduction.
First let us fix basic notations.
\begin{itemize}
\item $\g$ : affine Kac-Moody algebra of type $\dnone$ ($n\geq 4$).
\item $\uq(\g)$ : quantum group.
\item $I = \{0,1,\dots, n\}$ : index set.
\item $A = (a_{ij})_{i,j\in I}$ : generalized Cartan matrix for $\g$.
\item $P^{\vee} = \left(\op_{i\in I} \Z h_i\right) \op \Z d$ :
      dual weight lattice.
\item $\csa = \C \ot_{\Z} P^{\vee}$ : Cartan subalgebra.
\item $\ali$, $\delta$, $\La_i$ : simple root, null root, fundamental
      weight.
\item $P = \left( \op_{i\in I} \Z\La_i \right) \op \Z\delta$ :
      weight lattice.
\item $e_i, K_i^{\pm{1}}, f_i, q^d$ : generators of $\uq(\g)$.
\item $\uqp(\g)$ : subalgebra of $\uq(\g)$
      generated by $e_i, K_i^{\pm{1}}, f_i$ ($i\in I$).
\item $\bar{P} = \op_{i\in I} \Z\La_i$ : classical weight lattice.
\item $\wt, \cwt$ : (affine) weight, classical weight.
\item $\hwc(\la)$ : irreducible highest weight crystal of highest
      weight $\la$.
\item $\oldcrystal$ : level-$l$ perfect crystal of type $\dnone$.
\item $\fit, \eit$ : Kashiwara operators.
\item $\pathspace(\la)$ : the set of $\la$-paths
      (with crystal structure).
\end{itemize}

Below is a perfect crystal of type $\dnone$
introduced in~\cite{MR94a:17008,MR94j:17013,MR93g:17027,MR99h:17008}.
In the next section,
we give a new realization of this perfect crystal as a
certain set of blocks.

\begin{thm}\label{thm:22}
A level-$l$ perfect crystal of type $\dnone$
and the maps defining the crystal structure
is given as follows.
\begin{equation}
\oldcrystal =
 \left\{
 (x_1,\dots,x_n | \bar{x}_n,\dots, \bar{x}_1) \,\Big\vert\,
 \begin{aligned}
 & x_n=0 \textup{ or } \bar{x}_n=0, \ x_i, \bar{x}_i \in \Z_{\geq0}\\
 & \ \qquad \textstyle\sum_{i=1}^{n} (x_i + \bar{x}_i) = \textnormal{$l$}
 \end{aligned}
 \right\}.
\end{equation}

For $b = (x_1,\dots,x_n | \bar{x}_n,\dots, \bar{x}_1)$,
the action of the Kashiwara operator $\tilde{f}_0$, $\tilde{e}_0$
on $\oldcrystal$ are given by
\begin{equation*}
\begin{aligned}
\tilde{f}_0 b &=
  \begin{cases}
  (x_1,x_2+1,x_3,\dots,x_n|\bar{x}_n,\dots,\bar{x}_2,\bar{x}_1-1) &
  \textnormal{if $x_2 \geq \bar{x}_2$,}\\
  (x_1+1,x_2,\dots,x_n|\bar{x}_n,\dots,\bar{x}_3,\bar{x}_2-1,\bar{x}_1) &
  \textnormal{if $x_2 < \bar{x}_2$,}
  \end{cases}\\
\tilde{e}_0 b &=
  \begin{cases}
  (x_1,x_2-1,x_3,\dots,x_n|\bar{x}_n,\dots,\bar{x}_2,\bar{x}_1+1) &
  \textnormal{if $x_2 > \bar{x}_2$,}\\
  (x_1-1,x_2,\dots,x_n|\bar{x}_n,\dots,\bar{x}_3,\bar{x}_2+1,\bar{x}_1) &
  \textnormal{if $x_2 \leq \bar{x}_2$.}
  \end{cases}
\end{aligned}
\end{equation*}
For $i=1,\dots,n-2$,
\begin{equation*}
\begin{aligned}
\fit b &=
  \begin{cases}
  (x_1,\dots,x_{i-1},x_i-1,x_{i+1}+1,x_{i+2},\dots,x_n|
  \bar{x}_n,\dots,\bar{x}_1) &
  \textnormal{if $x_{i+1} \geq \bar{x}_{i+1}$,}\\
  (x_1,\dots,x_n|\bar{x}_n,\dots,\bar{x}_{i+2},\bar{x}_{i+1}-1,
  \bar{x}_{i}+1,\bar{x}_{i-1},\dots,\bar{x}_1) &
  \textnormal{if $x_{i+1} < \bar{x}_{i+1}$,}
  \end{cases}\\
\eit b &=
  \begin{cases}
  (x_1,\dots,x_{i-1},x_i+1,x_{i+1}-1,x_{i+2},\dots,x_n|
  \bar{x}_n,\dots,\bar{x}_1) &
  \textnormal{if $x_{i+1} > \bar{x}_{i+1}$,}\\
  (x_1,\dots,x_n|\bar{x}_n,\dots,\bar{x}_{i+2},\bar{x}_{i+1}+1,
  \bar{x}_{i}-1,\bar{x}_{i-1},\dots,\bar{x}_1) &
  \textnormal{if $x_{i+1} \leq \bar{x}_{i+1}$.}
  \end{cases}
\end{aligned}
\end{equation*}
And,
\begin{equation*}
\begin{aligned}
\tilde{f}_{n-1} b &=
  \begin{cases}
  (x_1,\dots,x_{n-2},x_{n-1}-1,x_n+1|\bar{x}_n,\dots,\bar{x}_1) &
  \textnormal{if $x_n\geq 0$, $\bar{x}_n=0$,}\\
  (x_1,\dots,x_n|\bar{x}_n-1,\bar{x}_{n-1}+1,\bar{x}_{n-2}\dots,\bar{x}_1) &
  \textnormal{if $x_n=0$, $\bar{x}_n\geq 1$,}
  \end{cases}\\
\tilde{e}_{n-1} b &=
  \begin{cases}
  (x_1,\dots,x_{n-2},x_{n-1}+1,x_n-1|\bar{x}_n,\dots,\bar{x}_1) &
  \textnormal{if $x_n\geq 1$, $\bar{x}_n=0$,}\\
  (x_1,\dots,x_n|\bar{x}_n+1,\bar{x}_{n-1}-1,\bar{x}_{n-2}\dots,\bar{x}_1) &
  \textnormal{if $x_n=0$, $\bar{x}_n\geq 0$,}
  \end{cases}
\end{aligned}
\end{equation*}
\begin{equation*}
\begin{aligned}
\tilde{f}_n b &=
  \begin{cases}
  (x_1,\dots,x_{n-1},x_n - 1|
  \bar{x}_n,\bar{x}_{n-1}+1,\bar{x}_{n-2},\dots,\bar{x}_1) &
  \textnormal{if $x_n\geq 1$, $\bar{x}_n=0$,}\\
  (x_1,\dots,x_{n-2},x_{n-1}-1,x_n|
  \bar{x}_n + 1,\bar{x}_{n-1},\dots,\bar{x}_1) &
  \textnormal{if $x_n=0$, $\bar{x}_n\geq 0$,}
  \end{cases}\\
\tilde{e}_n b  &=
  \begin{cases}
  (x_1,\dots,x_{n-1},x_n + 1|
  \bar{x}_n,\bar{x}_{n-1}-1,\bar{x}_{n-2},\dots,\bar{x}_1) &
  \textnormal{if $x_n\geq 0$, $\bar{x}_n=0$,}\\
  (x_1,\dots,x_{n-2},x_{n-1}+1,x_n|
  \bar{x}_n-1,\bar{x}_{n-1},\dots,\bar{x}_1) &
  \textnormal{if $x_n=0$, $\bar{x}_n\geq 1$.}
  \end{cases}\\
\end{aligned}
\end{equation*}
The remaining maps describing the crystal structure on $\oldcrystal$
are given below.
\begin{equation}
\begin{aligned}
&\vphi_0 (b) =
 \bar{x}_1 + (\bar{x}_2 - x_2)_+,\\
&\vphi_i (b) =
  x_i + (\bar{x}_{i+1} - x_{i+1})_+ \quad (i=1,\dots,n-2),\\
&\vphi_{n-1} (b) =
  x_{n-1} + {\bar x_n},\\
&\vphi_n (b) =
  x_{n-1} + x_n,\\
&\veps_0 (b) =
  x_1 + (x_2 - \bar{x}_2)_+,\\
&\veps_i (b) =
  \bar{x}_i + ({x}_{i+1} - \bar{x}_{i+1})_+ \quad (i=1,\dots,n-2),\\
&\veps_{n-1} (b) =
 \bar{x}_{n-1} + x_n,\\
&\veps_n (b) =
  \bar{x}_{n-1} + \bar{x}_n,\\
&\cwt(b) =
  \sum_{i=0}^n (\vphi_i(b) - \veps_i(b))\La_i.
\end{aligned}
\end{equation}
Here, $(x)_+ = \max(0,x)$.
\end{thm}

The action of the Kashiwara operators on Young walls,
which is our main objective, builds on the action of Kashiwara operator
on path elements.
Let us review the action of the Kashiwara operator $\eit$ and $\fit$
on a path element
\begin{equation}
\path = \cdots \ot p(j) \ot \cdots \ot p(1)\ot p(0) \in \pathspace(\la).
\end{equation}
\begin{enumerate}
\item Under each $p(j)$, write
      $\veps_i(p(j))$-many 1 followed by $\vphi_i(p(j))$-many 0.
\item From the (half-)infinite list of 0 and 1,
      successively cancel out each $(0,1)$ pair
      to obtain a sequence of 1 followed
      by some 0 (reading from left to right).
\item Let $\fit$ (and $\eit$) act on the $p(j)$
      corresponding to the left-most 0 (resp. right-most $1$) remaining.
      Set it to zero if no 0 (resp. $1$) remains.
\end{enumerate}

Following is the \emph{path realization} of irreducible highest
weight crystals.
Using this theorem,
we shall show that the set of Young walls satisfying some conditions
gives a new realization of the highest weight crystal graph.

\begin{thm}\label{thm:21}
The path space is isomorphic to
the irreducible highest weight crystal.
\begin{equation}
\mathcal{B}(\la) \cong \pathspace(\la).
\end{equation}
\end{thm}

%%%%%%%%%%%%%%%%%%%%%%%%%%%%%%%%%%%%%%%%%%%%%%%%%%%%%%%%%%%%%%%%%%%%%%%%%%
\section{Slices and splitting of blocks}%
\label{sec:3}

In this section,
we introduce the notions of \emph{splitting of blocks}
and \emph{slices}, and define a classical crystal structure
on the set $\newcrystal$ of equivalence classes of slices.
The classical crystal $\newcrystal$ is isomorphic to the level-$l$
perfect crystal $\oldcrystal$.
This fact will be shown in the next section.

The basic ingredient of our discussion will be
the following colored \defi{blocks} of two different shapes.
\savebox{\tmpfigb}{%
\begin{texdraw}
\fontsize{7}{7}\selectfont
\drawdim em
\textref h:C v:C
\setunitscale 1.9
\move(-1 0)\lvec(0 0)\lvec(0 1)\lvec(-1 1)\lvec(-1 0)
\move(0 0)\lvec(0.25 0.25)\lvec(0.25 1.25)
          \lvec(-0.75 1.25)\lvec(-1 1)
\move(0 1)\lvec(0.25 1.25)
\htext(-0.5 0.5){$i$}
\end{texdraw}
}%
\savebox{\tmpfigc}{%
\begin{texdraw}
\fontsize{7}{7}\selectfont
\drawdim em
\textref h:C v:C
\setunitscale 1.9
\move(-1 0)\lvec(0 0)\lvec(0 1)\lvec(-1 1)\lvec(-1 0)
\move(0 0)\lvec(0.5 0.5)\lvec(0.5 1.5)\lvec(-0.5 1.5)\lvec(-1 1)
\move(0 1.0)\lvec(0.5 1.5)
\htext(-0.5 0.5){$i$}
\end{texdraw}
}%
\begin{align*}
\raisebox{-0.5em}{\usebox{\tmpfigb}}
 &\ : \text{\ unit height, unit width, half-unit depth
             ($i=0,1,n-1,n$).}\\
\raisebox{-0.5em}{\usebox{\tmpfigc}}
 &\ : \text{\ unit height, unit width, unit depth ($i=2,\dots,n-2$).}
\end{align*}

\vspace{3mm}

For simplicity, we will use the following notation.\\[2mm]
\begin{center}
\begin{tabular}{rcl}
\raisebox{-0.4\height}{
\begin{texdraw}
\fontsize{7}{7}\selectfont
\drawdim em
\textref h:C v:C
\setunitscale 1.9
\move(-1 0)\lvec(0 0)\lvec(0 1)\lvec(-1 1)\lvec(-1 0)
\move(0 0)\lvec(0.5 0.5)\lvec(0.5 1.5)\lvec(-0.5 1.5)\lvec(-1 1)
\move(0 1.0)\lvec(0.5 1.5)
\htext(-0.5 0.5){$i$}
\end{texdraw}
}
& $\longleftrightarrow$ &
\raisebox{-0.4\height}{
\begin{texdraw}
\fontsize{7}{7}\selectfont
\drawdim em
\textref h:C v:C
\setunitscale 1.9
\move(-1 0)\lvec(0 0)\lvec(0 1)\lvec(-1 1)\lvec(-1 0)
\htext(-0.5 0.5){$i$}
\end{texdraw}
}
\end{tabular}
\quad
\begin{tabular}{rcl}
\raisebox{-0.4\height}{
\begin{texdraw}
\fontsize{7}{7}\selectfont
\drawdim em
\textref h:C v:C
\setunitscale 1.9
\move(-1 0)\lvec(0 0)\lvec(0 1)\lvec(-1 1)\lvec(-1 0)
\move(0 0)\lvec(0.25 0.25)\lvec(0.25 1.25)\lvec(-0.75 1.25)\lvec(-1 1)
\move(0 1)\lvec(0.25 1.25)
\lpatt(0.03 0.1)
\move(0 0)\lvec(-0.25 -0.25)\lvec(-1.25 -0.25)\lvec(-1 0)
\htext(-0.5 0.5){$j$}
\end{texdraw}
}
& $\longleftrightarrow$ &
\raisebox{-0.4\height}{
\begin{texdraw}
\fontsize{7}{7}\selectfont
\drawdim em
\textref h:C v:C
\setunitscale 1.9
\move(-1 0)\lvec(0 0)\lvec(0 1)\lvec(-1 1)\lvec(-1 0)
\move(0 1)\lvec(-1 0)
\htext(-0.25 0.25){$j$}
\end{texdraw}
}\\[4mm]
\raisebox{-0.4\height}{
\begin{texdraw}
\fontsize{7}{7}\selectfont
\drawdim em
\textref h:C v:C
\setunitscale 1.9
\move(-1 0)\lvec(0 0)\lvec(0 1)\lvec(-1 1)\lvec(-1 0)
\move(0 0)\lvec(0.25 0.25)\lvec(0.25 1.25)\lvec(-0.75 1.25)\lvec(-1 1)
\move(0 1)\lvec(0.25 1.25)
\lpatt(0.03 0.1)
\move(0.25 0.25)\lvec(0.5 0.5)\lvec(0.25 0.5)
\htext(-0.5 0.5){$i$}
\end{texdraw}
}
& $\longleftrightarrow$ &
\raisebox{-0.4\height}{
\begin{texdraw}
\fontsize{7}{7}\selectfont
\drawdim em
\textref h:C v:C
\setunitscale 1.9
\move(-1 0)\lvec(0 0)\lvec(0 1)\lvec(-1 1)\lvec(-1 0)
\move(0 1)\lvec(-1 0)
\htext(-0.75 0.75){$i$}
\end{texdraw}
}
\end{tabular}
\quad
\begin{tabular}{rcl}
\raisebox{-0.4\height}{
\begin{texdraw}
\fontsize{7}{7}\selectfont
\drawdim em
\textref h:C v:C
\setunitscale 1.9
\move(-1 0)\lvec(0 0)\lvec(0 1)\lvec(-1 1)\lvec(-1 0)
\move(0 0)\lvec(0.5 0.5)\lvec(0.5 1.5)\lvec(-0.5 1.5)\lvec(-1 1)
\move(0.25 0.25)\lvec(0.25 1.25)\lvec(-0.75 1.25)
\move(0 1.0)\lvec(0.5 1.5)
\htext(-0.5 0.5){$i$}
\end{texdraw}
}
& $\longleftrightarrow$ &
\raisebox{-0.4\height}{
\begin{texdraw}
\fontsize{7}{7}\selectfont
\drawdim em
\textref h:C v:C
\setunitscale 1.9
\move(-1 0)\lvec(0 0)\lvec(0 1)\lvec(-1 1)\lvec(-1 0)
\move(0 1)\lvec(-1 0)
\htext(-0.75 0.75){$i$}
\htext(-0.25 0.25){$j$}
\end{texdraw}
}
\end{tabular}
\end{center}
\vspace{1em}\noindent
That is, we use just the frontal view for most blocks.
And for half-unit depth blocks, we use triangles to cope with
their co-presence inside a single unit cube.

\begin{df}\label{df3.1}
We define a \defi{level-$1$ slice} of $\dnone$ type
to be the set of finitely many blocks of the above given shapes
stacked in one column of unit depth
following (either one of) the patterns given below.
\begin{center}
\begin{texdraw}
\fontsize{7}{7}\selectfont
\textref h:C v:C
\drawdim em
\setunitscale 1.9
\move(0 0)\rlvec(1 0)
\move(0 0)\rlvec(1 1)
\move(0 1)\rlvec(1 0)
\move(0 2)\rlvec(1 0)
\move(0 3.4)\rlvec(1 0)
\move(1 5.4)\rlvec(-0.5 -0.5)
\move(0 4.4)\rlvec(1 0)
\move(0 5.4)\rlvec(1 0)
\move(0 6.4)\rlvec(1 0)
\move(0 8.8)\rlvec(1 0)
\move(0 7.8)\rlvec(1 0)
\move(0 8.8)\rlvec(1 1)
\move(0 9.8)\rlvec(1 0)
\move(0 10.8)\rlvec(1 0)
\move(0 0)\rlvec(0 11.1)
\move(1 0)\rlvec(0 11.1)
\htext(0.75 0.25){$0$}
\htext(0.25 0.75){$1$}
\htext(0.5 1.5){$2$}
\vtext(0.5 2.8){$\cdots$}
\htext(0.5 3.9){$\!\!n\!\!-\!\!2\!\!$}
\htext(0.25 5.15){$n$}
\htext(0.55 4.65){$\!\!\!n\!\!-\!\!1\!\!\!$}
\htext(0.5 5.9){$\!\!n\!\!-\!\!2\!\!$}
\vtext(0.5 7.2){$\cdots$}
\htext(0.5 8.3){$2$}
\htext(0.75 9.05){$0$}
\htext(0.25 9.55){$1$}
\htext(0.5 10.3){$2$}
\textref h:L v:C
\htext(1.1 2.78){$\left.\rule{0pt}{5.5em}\right\}$
       \defi{supporting} blocks}
\textref h:R v:C
\htext(-0.1 7.2){\defi{covering} blocks
         $\left\{\rule{0pt}{5.5em}\right.$}
\htext(-0.1 0.5){covering blocks $\rightarrow$}
\end{texdraw}
\hspace{13em}
\begin{texdraw}
\fontsize{7}{7}\selectfont
\textref h:C v:C
\drawdim em
\setunitscale 1.9
\move(0 0)\rlvec(1 0)
\move(0 0)\rlvec(1 1)
\move(0 1)\rlvec(1 0)
\move(0 2)\rlvec(1 0)
\move(0 3.4)\rlvec(1 0)
\move(0 4.4)\rlvec(0.5 0.5)
\move(0 4.4)\rlvec(1 0)
\move(0 5.4)\rlvec(1 0)
\move(0 6.4)\rlvec(1 0)
\move(0 8.8)\rlvec(1 0)
\move(0 7.8)\rlvec(1 0)
\move(0 8.8)\rlvec(1 1)
\move(0 9.8)\rlvec(1 0)
\move(0 10.8)\rlvec(1 0)
\move(0 0)\rlvec(0 11.1)
\move(1 0)\rlvec(0 11.1)
\htext(0.75 0.25){$1$}
\htext(0.25 0.75){$0$}
\htext(0.5 1.5){$2$}
\vtext(0.5 2.8){$\cdots$}
\htext(0.5 3.9){$\!\!n\!\!-\!\!2\!\!$}
\htext(0.45 5.15){$\!\!\!n\!\!-\!\!1\!\!\!$}
\htext(0.75 4.65){$n$}
\htext(0.5 5.9){$\!\!n\!\!-\!\!2\!\!$}
\vtext(0.5 7.2){$\cdots$}
\htext(0.5 8.3){$2$}
\htext(0.75 9.05){$1$}
\htext(0.25 9.55){$0$}
\htext(0.5 10.3){$2$}
\textref h:L v:C
\htext(1.1 2.78){$\left.\rule{0pt}{5.5em}\right\}$
             \defi{supporting} blocks}
\textref h:R v:C
\htext(-0.1 7.2){\defi{covering} blocks
          $\left\{\rule{0pt}{5.5em}\right.$}
\htext(-0.1 0.5){covering blocks $\rightarrow$}
\end{texdraw}
\end{center}
In stacking the blocks,
no block should be placed on top of a column of half-unit depth.
\end{df}

We can see in the figure that this is a repeating pattern.
In these repeating patterns,
there are certain colors which appear twice in each cycle.
To distinguish the two positions containing the
same colored blocks,
we have given names to these positions and blocks.

We say that an $i$-block is a \defi{covering $i$-block}
(or a \defi{supporting $i$-block})
if it is closer to the $n-1$ and $n$-blocks
that sit below (resp. above) it
than to the $n-1$ and $n$-blocks above (resp. below) it.
A block sitting in a position that appears only once in each cycle
is regarded as being both a supporting block and a covering block.

An \defi{$i$-slot} is the top of a level-$1$ slice
where one may add an $i$-block.
The notions of \defi{covering $i$-slot}
and \defi{supporting $i$-slot} are self-explanatory.

We define a \defi{$\delta$-column} to be any set of blocks appearing
in a single cycle of the stacking pattern.
For a level-1 slice $c$, we define $c + \delta$
(and $c - \delta$) to be the level-$1$ slice obtained from
$c$ by adding (resp. removing) a $\delta$-column.
The following examples are for type $D_4^{(1)}$.
\begin{enumerate}
\item
$\delta$-columns following the first pattern of Definition \ref{df3.1} :\\
\savebox{\tmpfiga}
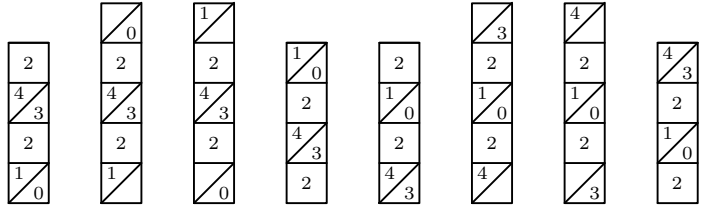
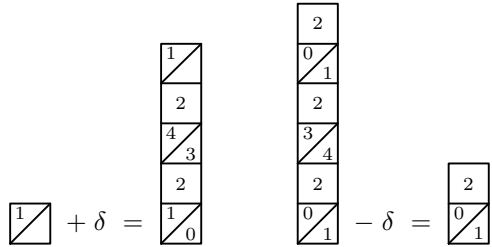
{
\begin{texdraw}
\fontsize{7}{7}\selectfont
\textref h:C v:C
\drawdim em
\setunitscale 1.9
\move(0 0)\lvec(1 0)\lvec(1 4)\lvec(0 4)\lvec(0 0)
\move(0 0)\lvec(1 1)
\move(0 1)\lvec(1 1)
\move(0 2)\lvec(1 2)
\move(0 2)\rlvec(1 1)
\move(0 3)\lvec(1 3)
\move(0 4)\lvec(1 4)
\htext(0.75 0.25){$0$}\htext(0.25 0.75){$1$}
\htext(0.5 1.5){$2$}\htext(0.25 2.75){$4$}
\htext(0.75 2.25){$3$}\htext(0.5 3.5){$2$}
\end{texdraw}
}%
\savebox{\tmpfigb}{
\begin{texdraw}
\fontsize{7}{7}\selectfont
\textref h:C v:C
\drawdim em
\setunitscale 1.9
\move(0 0)\rlvec(1 0)\rlvec(0 5)\rlvec(-1 0)\rlvec(0 -5)
\move(0 1)\rlvec(1 0)
\move(0 2)\rlvec(1 0)
\move(0 3)\rlvec(1 0)
\move(0 4)\rlvec(1 0)
\move(0 0)\rlvec(1 1)
\move(0 2)\rlvec(1 1)
\move(0 4)\rlvec(1 1)
\htext(0.75 4.25){$0$}\htext(0.25 0.75){$1$}
\htext(0.5 1.5){$2$}\htext(0.25 2.75){$4$}
\htext(0.75 2.25){$3$}\htext(0.5 3.5){$2$}
\end{texdraw}
}%
\savebox{\tmpfigc}{
\begin{texdraw}
\fontsize{7}{7}\selectfont
\textref h:C v:C
\drawdim em
\setunitscale 1.9
\move(0 0)\rlvec(1 0)\rlvec(0 5)\rlvec(-1 0)\rlvec(0 -5)
\move(0 1)\rlvec(1 0)
\move(0 2)\rlvec(1 0)
\move(0 3)\rlvec(1 0)
\move(0 4)\rlvec(1 0)
\move(0 0)\rlvec(1 1)
\move(0 2)\rlvec(1 1)
\move(0 4)\rlvec(1 1)
\htext(0.25 4.75){$1$}\htext(0.75 0.25){$0$}
\htext(0.5 1.5){$2$}\htext(0.25 2.75){$4$}
\htext(0.75 2.25){$3$}\htext(0.5 3.5){$2$}
\end{texdraw}
}%
\savebox{\tmpfigd}{
\begin{texdraw}
\fontsize{7}{7}\selectfont
\textref h:C v:C
\drawdim em
\setunitscale 1.9
\move(0 0)\rlvec(1 0)\rlvec(0 4)\rlvec(-1 0)\rlvec(0 -4)
\move(0 1)\rlvec(1 0)
\move(0 2)\rlvec(1 0)
\move(0 3)\rlvec(1 0)
\move(0 4)\rlvec(1 0)
\move(0 1)\rlvec(1 1)
\move(0 3)\rlvec(1 1)
\htext(0.25 3.75){$1$}\htext(0.75 3.25){$0$}
\htext(0.5 0.5){$2$}\htext(0.25 1.75){$4$}
\htext(0.75 1.25){$3$}\htext(0.5 2.5){$2$}
\end{texdraw}
}%
\savebox{\tmpfige}{
\begin{texdraw}
\fontsize{7}{7}\selectfont
\textref h:C v:C
\drawdim em
\setunitscale 1.9
\move(0 0)\rlvec(1 0)\rlvec(0 4)\rlvec(-1 0)\rlvec(0 -4)
\move(0 1)\rlvec(1 0)
\move(0 2)\rlvec(1 0)
\move(0 3)\rlvec(1 0)
\move(0 4)\rlvec(1 0)
\move(0 0)\rlvec(1 1)
\move(0 2)\rlvec(1 1)
\htext(0.25 2.75){$1$}\htext(0.75 2.25){$0$}
\htext(0.5 1.5){$2$}\htext(0.25 0.75){$4$}
\htext(0.75 0.25){$3$}\htext(0.5 3.5){$2$}
\end{texdraw}
}%
\savebox{\tmpfigf}{
\begin{texdraw}
\fontsize{7}{7}\selectfont
\textref h:C v:C
\drawdim em
\setunitscale 1.9
\move(0 0)\rlvec(1 0)\rlvec(0 5)\rlvec(-1 0)\rlvec(0 -5)
\move(0 1)\rlvec(1 0)
\move(0 2)\rlvec(1 0)
\move(0 3)\rlvec(1 0)
\move(0 4)\rlvec(1 0)
\move(0 0)\rlvec(1 1)
\move(0 2)\rlvec(1 1)
\move(0 4)\rlvec(1 1)
\htext(0.75 2.25){$0$}\htext(0.25 2.75){$1$}
\htext(0.5 1.5){$2$}\htext(0.25 0.75){$4$}
\htext(0.75 4.25){$3$}\htext(0.5 3.5){$2$}
\end{texdraw}
}%
\savebox{\tmpfigg}{
\begin{texdraw}
\fontsize{7}{7}\selectfont
\textref h:C v:C
\drawdim em
\setunitscale 1.9
\move(0 0)\rlvec(1 0)\rlvec(0 5)\rlvec(-1 0)\rlvec(0 -5)
\move(0 1)\rlvec(1 0)
\move(0 2)\rlvec(1 0)
\move(0 3)\rlvec(1 0)
\move(0 4)\rlvec(1 0)
\move(0 0)\rlvec(1 1)
\move(0 2)\rlvec(1 1)
\move(0 4)\rlvec(1 1)
\htext(0.25 2.75){$1$}\htext(0.75 2.25){$0$}
\htext(0.5 1.5){$2$}\htext(0.25 4.75){$4$}
\htext(0.75 0.25){$3$}\htext(0.5 3.5){$2$}
\end{texdraw}
}%
\savebox{\tmpfigh}{
\begin{texdraw}
\fontsize{7}{7}\selectfont
\textref h:C v:C
\drawdim em
\setunitscale 1.9
\move(0 0)\rlvec(1 0)\rlvec(0 4)\rlvec(-1 0)\rlvec(0 -4)
\move(0 1)\rlvec(1 0)
\move(0 2)\rlvec(1 0)
\move(0 3)\rlvec(1 0)
\move(0 4)\rlvec(1 0)
\move(0 1)\rlvec(1 1)
\move(0 3)\rlvec(1 1)
\htext(0.25 3.75){$4$}\htext(0.75 3.25){$3$}
\htext(0.5 0.5){$2$}\htext(0.25 1.75){$1$}
\htext(0.75 1.25){$0$}\htext(0.5 2.5){$2$}
\end{texdraw}
}%
\begin{center}
\raisebox{-7.1em}{\usebox{\tmpfiga}}
\quad
\raisebox{-7.1em}{\usebox{\tmpfigb}}
\quad
\raisebox{-7.1em}{\usebox{\tmpfigc}}
\quad
\raisebox{-7.1em}{\usebox{\tmpfigd}}
\quad
\raisebox{-7.1em}{\usebox{\tmpfige}}
\quad
\raisebox{-7.1em}{\usebox{\tmpfigf}}
\quad
\raisebox{-7.1em}{\usebox{\tmpfigg}}
\quad
\raisebox{-7.1em}{\usebox{\tmpfigh}}
\end{center}

\vspace{4mm}

\savebox{\tmpfiga}{
\begin{texdraw}
\fontsize{7}{7}\selectfont
\textref h:C v:C
\drawdim em
\setunitscale 1.9
\move(0 0)\rlvec(1 0)\rlvec(0 1)\rlvec(-1 0)\rlvec(0 -1)
\move(0 0)\lvec(1 1)
\htext(0.25 0.75){$1$}
\end{texdraw}
\raisebox{0.5em}{$\;+\;\delta\ =\ $}
\begin{texdraw}
\fontsize{7}{7}\selectfont
\textref h:C v:C
\drawdim em
\setunitscale 1.9
\move(0 0)\lvec(1 0)\lvec(1 5)\lvec(0 5)\lvec(0 0)
\move(0 0)\lvec(1 1)
\move(0 1)\lvec(1 1)
\move(0 2)\lvec(1 2)
\move(0 2)\rlvec(1 1)
\move(0 3)\lvec(1 3)
\move(0 4)\lvec(1 4)
\move(0 4)\lvec(1 5)
\htext(0.75 0.25){$0$}\htext(0.25 0.75){$1$}
\htext(0.5 1.5){$2$}\htext(0.75 2.25){$3$}
\htext(0.25 2.75){$4$}\htext(0.5 3.5){$2$}
\htext(0.25 4.75){$1$}
\end{texdraw}}
\savebox{\tmpfigb}{
\begin{texdraw}
\fontsize{7}{7}\selectfont
\textref h:C v:C
\drawdim em
\setunitscale 1.9
\move(0 0)\lvec(1 0)\lvec(1 6)\lvec(0 6)\lvec(0 0)
\move(0 0)\lvec(1 1)
\move(0 1)\lvec(1 1)
\move(0 2)\lvec(1 2)
\move(0 2)\rlvec(1 1)
\move(0 3)\lvec(1 3)
\move(0 4)\lvec(1 4)
\move(0 4)\lvec(1 5)
\move(0 5)\rlvec(1 0)
\htext(0.75 0.25){$1$}\htext(0.25 0.75){$0$}
\htext(0.5 1.5){$2$}\htext(0.75 2.25){$4$}
\htext(0.25 2.75){$3$}\htext(0.5 3.5){$2$}
\htext(0.25 4.75){$0$}
\htext(0.75 4.25){$1$}
\htext(0.5 5.5){$2$}
\end{texdraw}
\raisebox{0.5em}{$\;-\;\delta\ =\ $}
\begin{texdraw}
\fontsize{7}{7}\selectfont
\textref h:C v:C
\drawdim em
\setunitscale 1.9
\move(0 0)\rlvec(1 0)\rlvec(0 2)\rlvec(-1 0)\rlvec(0 -2)
\move(0 0)\lvec(1 1)
\move(0 1)\rlvec(1 0)
\htext(0.25 0.75){$0$}
\htext(0.75 0.25){$1$}
\htext(0.5 1.5){$2$}
\end{texdraw}}
\item
$c \pm \delta$ :
\begin{center}
\raisebox{-7.1em}{\usebox{\tmpfiga}}
\qquad\quad
\raisebox{-7.1em}{\usebox{\tmpfigb}}
\end{center}
\end{enumerate}

\vspace{3mm}

For two level-$1$ slices $c$ and $c'$,
we write $c \subset c'$ if $c'$ contains all blocks
which make up $c$.

\begin{df}
An ordered $l$-tuple $C = (c_1, \dots, c_l)$ of level-$1$ slices
is a \defi{level-$l$ pre-slice} for type $\dnone$
if it satisfies the following condition.
\begin{equation}\label{eq:inclusion}
c_1 \subset c_2 \subset \cdots \subset c_l \subset c_1 + \delta.
\end{equation}
Only one of the two stacking patterns
should be used in building a level-$l$ pre-slice.
Each $c_i$ is called the $i$th \defi{layer} of $C$.
\end{df}

We shall often just say \emph{pre-slice},
when dealing with level-$l$ pre-slices.
Mentally, we picture a level-$l$ pre-slice as a gathering of the columns
with the $i$th layer placed in front of the $(i+1)$th layer,
rather than as an ordered $l$-tuple.

To simplify drawings,
we shall use just the side view
when representing a pre-slice.
We explain this method used in drawing a pre-slice
with the following example for $D_4^{(1)}$-type.
\begin{center}
\raisebox{0.5em}{$(\: c_1\:=\: $}
\begin{texdraw}
\fontsize{7}{7}\selectfont
\textref h:C v:C
\drawdim em
\setunitscale 1.9
\move(0 0)\lvec(1 0)\lvec(1 2)\lvec(0 2)\lvec(0 0)
\move(0 0)\lvec(1 1)
\move(0 1)\lvec(1 1)
\htext(0.75 0.25){$0$} \htext(0.25 0.75){$1$}
\htext(0.5 1.5){$2$}
\end{texdraw}\;,
\raisebox{0.5em}{$c_2\:=\: $}
\begin{texdraw}
\fontsize{7}{7}\selectfont
\textref h:C v:C
\drawdim em
\setunitscale 1.9
\move(0 0)\lvec(1 0)\lvec(1 3)\lvec(0 3)\lvec(0 0)
\move(0 1)\lvec(1 1) \move(0 2)\lvec(1 2)
\move(0 0)\lvec(1 1)
\move(0 2)\rlvec(1 1)
\htext(0.75 0.25){$0$} \htext(0.25 0.75){$1$}
\htext(0.5 1.5){$2$} \htext(0.25 2.75){$4$}
\end{texdraw}\;,
\raisebox{0.5em}{$c_3\:=\: $}
\begin{texdraw}
\fontsize{7}{7}\selectfont
\textref h:C v:C
\drawdim em
\setunitscale 1.9
\move(0 0)\lvec(1 0)\lvec(1 3)\lvec(0 3)\lvec(0 0)
\move(0 1)\lvec(1 1) \move(0 2)\lvec(1 2)
\move(0 0)\lvec(1 1)
\move(0 2)\rlvec(1 1)
\htext(0.75 0.25){$0$} \htext(0.25 0.75){$1$}
\htext(0.5 1.5){$2$} \htext(0.75 2.25){$3$}
\htext(0.25 2.75){$4$}
\end{texdraw}\;,
\raisebox{0.5em}{$c_4\:=\: $}
\begin{texdraw}
\fontsize{7}{7}\selectfont
\textref h:C v:C
\drawdim em
\setunitscale 1.9
\move(0 0)\lvec(1 0)\lvec(1 5)\lvec(0 5)\lvec(0 0)
\move(0 1)\lvec(1 1) \move(0 2)\lvec(1 2)
\move(0 0)\lvec(1 1)
\move(0 3)\lvec(1 3) \move(0 4)\lvec(1 4)
\move(0 4)\lvec(1 5)
\move(0 2)\rlvec(1 1)
\htext(0.75 0.25){$0$} \htext(0.25 0.75){$1$}
\htext(0.5 1.5){$2$} \htext(0.75 2.25){$3$}
\htext(0.25 2.75){$4$} \htext(0.5 3.5){$2$}
\htext(0.75 4.25){$0$}
\end{texdraw}\;,
\raisebox{0.5em}{$c_5\:=\: $}
\begin{texdraw}
\fontsize{7}{7}\selectfont
\textref h:C v:C
\drawdim em
\setunitscale 1.9
\move(0 0)\lvec(1 0)\lvec(1 5)\lvec(0 5)\lvec(0 0)
\move(0 1)\lvec(1 1) \move(0 2)\lvec(1 2)
\move(0 0)\lvec(1 1)
\move(0 3)\lvec(1 3) \move(0 4)\lvec(1 4)
\move(0 4)\lvec(1 5)
\move(0 2)\rlvec(1 1)
\htext(0.75 0.25){$0$}\htext(0.25 0.75){$1$}
\htext(0.5 1.5){$2$}\htext(0.25 2.75){$4$}
\htext(0.75 2.25){$3$}\htext(0.5 3.5){$2$}
\htext(0.75 4.25){$0$}\htext(0.25 4.75){$1$}
\end{texdraw}
\raisebox{0.5em}{$\: )$}

\vspace{3mm}

\raisebox{0.6em}{$\ =\ \ $}
\raisebox{-1.3em}{\begin{texdraw}
\fontsize{7}{7}\selectfont
\textref h:C v:C
\drawdim em
\setunitscale 1.9
\move(0.4 -0.3)
  \bsegment
  \move(0 0)\lvec(-1 0)\lvec(-1 2)\lvec(0 2)\lvec(0 0)
  \move(0 1)\lvec(-1 1)
  \htext(-0.5 0.5){$1$}
  \htext(-0.5 1.5){$2$}
  \esegment
\move(0 0)
  \bsegment
  \move(-1 0)\lvec(-1 3)\lvec(-0 3)\lvec(-0 2)\lvec(-1 2)
  \htext(-0.5 2.5){$4$}
  \esegment
\move(-0.4 0.3)
  \bsegment
  \move(-1 0)\lvec(-1 2)\lvec(-0.85 2)
  \move(0 3)\rlvec(0 -0.15)
  \move(-0.2 3.15)\rlvec(-1 0)
  \esegment
\move(-0.8 0.6)
  \bsegment
  \move(-0 3)\lvec(-0 4)\lvec(-1 4)\lvec(-1 2)
  \move(-0 3)\lvec(-1 3)
  \move(-1 0)\rlvec(-0 2)
  \htext(-0.5 3.5){$2$}
  \esegment
\move(-1.0 0.75)
  \bsegment
  \move(-1 0)\lvec(-1 1)
  \move(-1 4)\lvec(-1 5)\lvec(-0 5)\lvec(-0 4)\lvec(-1 4)
  \htext(-0.5 4.5){$0$}
  \esegment
\move(-1.2 0.9)
  \bsegment
  \move(-1 0)\lvec(-1 5)\lvec(-0 5)
  \esegment
\move(-1.4 1.05)
  \bsegment
  \move(-1 0)\lvec(-1 1)
  \move(-1 4)\lvec(-1 5)\lvec(-0 5)
  \esegment
\move(-1.6 1.2)
  \bsegment
  \move(-1 0)\lvec(-1 5)\lvec(-0 5)
  \esegment

\move(-0.6 -0.3)\lvec(-2.6 1.2)
\move(-0.6 0.7)\lvec(-2.6 2.2)
\move(-0.6 1.7)\lvec(-2.6 3.2)
\move(0.4 1.7)\lvec(-0 2)

\move(-0.4 3.3)\lvec(-0.8 3.6)
\move(-1.4 3.3)\lvec(-2.6 4.2)
\move(-0.4 3.3)\rlvec(-1 0)
\move(-1.4 3.3)\rlvec(-0 -1)
\move(-0 3)\rlvec(-0.2 0.15)
\move(-1 3)\rlvec(-0.2 0.15)
\move(-0.2 3.15)\rlvec(-1 0)
\move(-1.2 3.15)\rlvec(-0 -1)

\move(-1.8 4.6)\lvec(-2.6 5.2)
\move(-0.8 4.6)\lvec(-1 4.75)
\move(-2 5.75)\rlvec(-0.6 0.45)
\move(-1 5.75)\rlvec(-0.6 0.45)

\move(-0.8 -0.15)\lvec(-0.8 0.85)
\move(-1.2 0.15)\lvec(-1.2 1.15)
\move(-1.6 0.45)\lvec(-1.6 1.45)
\move(-1.6 2.45)\lvec(-1.6 3.45)
\move(-2 2.75)\rlvec(-0 1)
\move(-2.4 3.05)\rlvec(-0 1)
\end{texdraw}}
\raisebox{0.6em}{$\quad =\ \ $}
\begin{texdraw}
\fontsize{7}{7}\selectfont
\textref h:C v:C
\drawdim em
\setunitscale 1.9
\move(0 0)
\bsegment
\move(0 0)\lvec(1 0)\lvec(1 5)\lvec(0 5)\lvec(0 0)
\move(0 1)\lvec(1 1) \move(0 2)\lvec(1 2)
\move(0.5 0)\lvec(0.5 1)
\move(0 3)\lvec(1 3) \move(0 4)\lvec(1 4)
\move(0.5 2)\rlvec(0 1)
\move(0.5 4)\lvec(0.5 5)
\htext(0.25 0.5){$0$}\htext(0.75 0.5){$1$}
\htext(0.5 1.5){$2$}\htext(0.25 2.5){$3$}
\htext(0.75 2.5){$4$}\htext(0.5 3.5){$2$}
\htext(0.25 4.5){$0$}\htext(0.75 4.5){$1$}
\esegment
\move(1 0)
\bsegment
\move(0 0)\lvec(1 0)\lvec(1 4)\lvec(0.5 4)\lvec(0.5 5)\lvec(0 5)
\move(0 1)\lvec(1 1)
\move(0 2)\lvec(1 2) \move(0.5 0)\lvec(0.5 1)
\move(0 3)\lvec(1 3)
\move(0.5 2)\rlvec(0 1)
\move(0 4)\lvec(1 4) \move(0.5 4)\lvec(0.5 5)
\htext(0.25 0.5){$0$} \htext(0.75 0.5){$1$}
\htext(0.5 1.5){$2$} \htext(0.25 2.5){$3$}
\htext(0.75 2.5){$4$} \htext(0.5 3.5){$2$}
\htext(0.25 4.5){$0$}
\esegment
\move(2 0)
\bsegment
\move(0 0)\lvec(1 0)\lvec(1 2)\lvec(0 2)\lvec(0 0)
\move(0 1)\lvec(1 1) \move(0 2)\lvec(1 2)
\move(0.5 0)\lvec(0.5 1)
\move(0.5 3)\lvec(1 3)\lvec(1 2)
\move(0 3)\rlvec(0.5 0)\rlvec(0 -1)
\htext(0.25 0.5){$0$} \htext(0.75 0.5){$1$}
\htext(0.5 1.5){$2$}\htext(0.25 2.5){$3$}
\htext(0.75 2.5){$4$}
\esegment
\move(3 0)
\bsegment
\move(0 0)\lvec(1 0)\lvec(1 2)\lvec(0 2)\lvec(0 0)
\move(0 1)\lvec(1 1) \move(0 2)\lvec(1 2)
\move(0.5 0)\lvec(0.5 1)
\move(0.5 2)\rlvec(0 1)\rlvec(0.5 0)\rlvec(0 -1)
\htext(0.25 0.5){$0$} \htext(0.75 0.5){$1$}
\htext(0.5 1.5){$2$}\htext(0.75 2.5){$4$}
\esegment
\move(4 0)
\bsegment
\move(0 0)\lvec(1 0)\lvec(1 2)\lvec(0 2)\lvec(0 0)
\move(0.5 0)\lvec(0.5 1) \move(0 1)\lvec(1 1)
\htext(0.25 0.5){$0$} \htext(0.75 0.5){$1$}
\htext(0.5 1.5){$2$}
\esegment
\end{texdraw}
\end{center}

\vspace{2mm}

Next we explain the notion of \defi{splitting an $i$-block}
in a level-$l$ pre-slice.
To do this, we first define a \defi{whole block} to be
any one of the unit depth blocks, or a gluing of two
half-unit depth blocks, as given by the following figures.\\
\begin{center}
\begin{texdraw}
\fontsize{7}{7}\selectfont
\textref h:C v:C
\drawdim em
\setunitscale 1.9
\move(0 0)
\bsegment
\move(0 0)\lvec(0 1)\lvec(1 1)\lvec(1 0)\lvec(0 0)\lvec(1 1)
\htext(0.25 0.75){$1$}
\htext(0.75 0.25){$0$}
\esegment
\rmove(2 0)
\bsegment
\move(0 0)\lvec(0 1)\lvec(1 1)\lvec(1 0)\lvec(0 0)\lvec(1 1)
\htext(0.25 0.75){$0$}
\htext(0.75 0.25){$1$}
\esegment
\rmove(2 0)
\bsegment
\move(0 0)\lvec(0 1)\lvec(1 1)\lvec(1 0)\lvec(0 0)
\htext(0.5 0.5){$2$}
\esegment
\rmove(3 0.5)
\htext{$\cdots\cdots$}
\rmove(2 -0.5)
\bsegment
\move(0 0)\lvec(0 1)\lvec(1 1)\lvec(1 0)\lvec(0 0)
\htext(0.5 0.5){$n\!\!-\!\!2$}
\esegment
\rmove(2 0)
\bsegment
\move(1 1)\lvec(1 0)\lvec(0 0)\lvec(0 1)\lvec(1 1)\lvec(0.5 0.5)
\htext(0.25 0.75){$n$}
\htext(0.55 0.25){$n\!\!-\!\!1$}
\esegment
\rmove(2 0)
\bsegment
\move(0 0)\lvec(0 1)\lvec(1 1)\lvec(1 0)\lvec(0 0)\lvec(0.5 0.5)
\htext(0.45 0.75){$n\!\!-\!\!1$}
\htext(0.75 0.25){$n$}
\esegment
\end{texdraw}
\end{center}
\vspace{2mm}
The first two whole blocks will be named as $\zo$-blocks
and the last two blocks shall be named $\nn$-blocks.

For each $i=\zo, 0,1,\dots,n, \nn$,
note that in any fixed pre-slice, due to the circular inclusion
condition~(\ref{eq:inclusion}) on level-$1$ slices,
there can be at most two heights
in which a covering (or supporting) $i$-block may appear
as the top block of a layer.
Similarly, there can be at most two heights
in which a supporting (or covering) $i$-slot may appear.

\begin{df}\label{3.3}
Fix an $i=\zo, 2, \dots, n-2,\nn$.
Suppose that in the pre-slice under consideration,
there are layers with the top
occupied by a covering (whole) $i$-block and layers whose top is
a supporting (whole) $i$-slot.
Choose the covering $i$-block lying in the fore-most layer
(i.e., the one with the smallest layer index)
among the ones sitting in the higher height,
and the supporting $i$-slot
lying in the rear-most layer (i.e., the one with the largest layer index)
among the ones situated in the lower height.

To \defi{split an $i$-block} in a level-$l$ pre-slice
means to break off the top half of the
chosen covering $i$-block
and to place it in the chosen supporting $i$-slot.

The action characterized by
interchanging \emph{covering} and
\emph{supporting} in the above description is also
defined to be \defi{splitting of an $i$-block} in a level-$l$
pre-slice.
\end{df}

\begin{remark}
It is possible to develop the rest of the theory of
Young walls for type $\dnone$ after defining
\emph{splitting of a block} as either breaking a supporting
block and placing it in a covering slot or vice versa.
We chose to do as the above in view of the symmetry
$\dnone$ Dynkin diagram exhibits.
\end{remark}

\begin{remark}
Note that the result obtained after of splitting is no longer
a pre-slice.
Hence, logically, it is possible to split only a single whole
block from a pre-slice.
But the above definition of splitting can easily be
applied to such split results also.
We shall allow for this extension of the above definition
to be used.
\end{remark}

\begin{remark}
Although not strictly true, one may think of the choices of block
and slot in the above definition for splitting as inevitable if one
wants to preserve the circular inclusion condition~(\ref{eq:inclusion}).
\end{remark}

\begin{remark}
The splitting of one block may render the splitting of
another block impossible.
So the choice of which block to split first can
affect the next choice for splitting.
\end{remark}

Dotted lines shall be used to denote broken blocks,
as can be seen in the next example.
\begin{example}
The following illustrates splitting of blocks for $D_4^{(1)}$-type
pre-slices.
\begin{enumerate}
\item
We've split a single pre-slice in two different ways below.
The upper right figure shows the splitting of a $\tf$-block.
The half $\tf$-block sitting in the third (resp. second)
layer is the broken off top (resp. bottom) part of the
$\tf$-block of the second layer shown in the left figure.
The lower right figure shows the splitting of a covering
$2$-block from the same pre-slice.
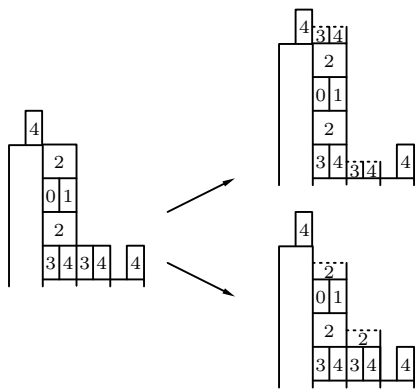
\begin{center}
\begin{texdraw}%
\drawdim in
\arrowheadsize l:0.065 w:0.03
\arrowheadtype t:F
\fontsize{7}{7}\selectfont
\textref h:C v:C
\drawdim em
\setunitscale 1.6
\move(0 0)
\bsegment
\move(0 -0.2)\rlvec(0 4.2)\rlvec(1 0)\rlvec(0 -4.2)
\move(0.5 4)\rlvec(0 1)\rlvec(0.5 0)\rlvec(0 -1)
\move(1 4)\rlvec(1 0)\rlvec(0 -4.2)
\move(1 3)\rlvec(1 0)
\move(1 2)\rlvec(1 0)
\move(1 1)\rlvec(2 0)\rlvec(0 -1.2)
\move(1 0)\rlvec(3 0)\rlvec(0 -0.2)
\move(3.5 0)\lvec(3.5 1)\lvec(4 1)\lvec(4 0)
\move(1.5 2)\rlvec(0 1)
\move(2.5 0)\rlvec(0 1)
\move(1.5 0)\rlvec(0 1)
\htext(0.75 4.5){$4$}
\htext(1.5 3.5){$2$}
\htext(1.25 2.5){$0$}
\htext(1.75 2.5){$1$}
\htext(1.5 1.5){$2$}
\htext(2.25 0.5){$3$}
\htext(2.75 0.5){$4$}
\htext(3.75 0.5){$4$}
\htext(1.25 0.5){$3$}
\htext(1.75 0.5){$4$}
\esegment
\move(8 3)
\bsegment
\move(0 -0.2)\rlvec(0 4.2)\rlvec(1 0)\rlvec(0 -4.2)
\move(0.5 4)\rlvec(0 1)\rlvec(0.5 0)\rlvec(0 -1)
\move(1 4)\rlvec(1 0)\rlvec(0 -4.2)
\move(1 3)\rlvec(1 0)
\move(1 2)\rlvec(1 0)
\move(1 1)\rlvec(1 0)
\move(3 -0.2)\rlvec(0 0.7)
\move(1 0)\rlvec(3 0)\rlvec(0 -0.2)
\move(3.5 0)\lvec(3.5 1)\lvec(4 1)\lvec(4 0)
\move(1.5 2)\rlvec(0 1)
\move(2.5 0)\rlvec(0 0.5)
\move(1.5 4.5)\rlvec(0 -0.5)
\move(2 4.5)\rlvec(0 -0.5)
\move(1.5 0)\rlvec(0 1)
\lpatt(0.05 0.15)
\move(2 0.5)\rlvec(1 0)
\move(1 4.5)\rlvec(1 0)
\htext(0.75 4.5){$4$}
\htext(1.5 3.5){$2$}
\htext(1.25 2.5){$0$}
\htext(1.75 2.5){$1$}
\htext(1.5 1.5){$2$}
\htext(2.25 0.25){$3$}
\htext(2.75 0.25){$4$}
\htext(3.75 0.5){$4$}
\htext(1.25 4.25){$3$}
\htext(1.75 4.25){$4$}
\htext(1.25 0.5){$3$}
\htext(1.75 0.5){$4$}
\esegment
\move(8 -3)
\bsegment
\move(0 -0.2)\rlvec(0 4.2)\rlvec(1 0)\rlvec(0 -4.2)
\move(0.5 4)\rlvec(0 1)\rlvec(0.5 0)\rlvec(0 -1)
\move(2 3.5)\rlvec(0 -3.7)
\move(1 3)\rlvec(1 0)
\move(1 2)\rlvec(1 0)
\move(1 1)\rlvec(1 0)
\move(3 -0.2)\rlvec(0 1.7)
\move(1 0)\rlvec(3 0)\rlvec(0 -0.2)
\move(1.5 2)\rlvec(0 1)
\move(2.5 0)\rlvec(0 1)
\move(1.5 0)\rlvec(0 1)
\move(2 1)\rlvec(1 0)\rlvec(0 -1)
\move(3.5 0)\lvec(3.5 1)\lvec(4 1)\lvec(4 0)
\lpatt(0.05 0.15)
\move(1 3.5)\rlvec(1 0)
\move(2 1.5)\rlvec(1 0)
\htext(0.75 4.5){$4$}
\htext(1.5 3.25){$2$}
\htext(1.25 2.5){$0$}
\htext(1.75 2.5){$1$}
\htext(1.5 1.5){$2$}
\htext(2.5 1.25){$2$}
\htext(2.25 0.5){$3$}
\htext(2.75 0.5){$4$}
\htext(3.75 0.5){$4$}
\htext(1.25 0.5){$3$}
\htext(1.75 0.5){$4$}
\esegment

\move(4.7 2)\ravec(2 1)
\move(4.7 0.5)\ravec(2 -1)
\end{texdraw}%
\end{center}
The circular inclusion condition~(\ref{eq:inclusion}) is no longer
met in a strict sense, but one may still consider it true
in a weaker (volume) sense.

\vspace{2mm}

\item This example shows four different split results
      of a single pre-slice. None of the results allow
      further splitting.\\[2mm]
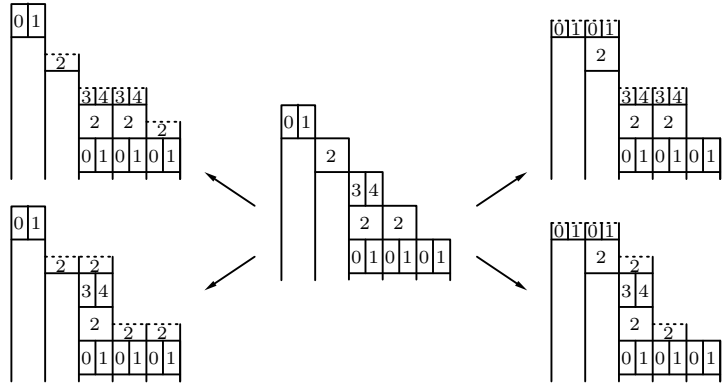
\begin{center}
\begin{texdraw}%
\drawdim in
\arrowheadsize l:0.065 w:0.03
\arrowheadtype t:F
\fontsize{7}{7}\selectfont
\textref h:C v:C
\drawdim em
\setunitscale 1.6
\move(0 0)
\bsegment
\move(0 -0.2)\rlvec(0 5.2)\rlvec(1 0)\rlvec(0 -5.2)
\move(0 4)\rlvec(2 0)\rlvec(0 -4.2)
\move(1 3)\rlvec(2 0)\rlvec(0 -3.2)
\move(2 2)\rlvec(2 0)\rlvec(0 -2.2)
\move(2 1)\rlvec(3 0)\rlvec(0 -1.2)
\move(2 0)\rlvec(3 0)
\move(0.5 4)\rlvec(0 1)
\move(2.5 2)\rlvec(0 1)
\move(2.5 0)\rlvec(0 1)
\move(3.5 0)\rlvec(0 1)
\move(4.5 0)\rlvec(0 1)
\htext(0.25 4.5){$0$}
\htext(0.75 4.5){$1$}
\htext(1.5 3.5){$2$}
\htext(2.25 2.5){$3$}
\htext(2.75 2.5){$4$}
\htext(2.5 1.5){$2$}
\htext(3.5 1.5){$2$}
\htext(2.25 0.5){$0$}
\htext(2.75 0.5){$1$}
\htext(3.25 0.5){$0$}
\htext(3.75 0.5){$1$}
\htext(4.25 0.5){$0$}
\htext(4.75 0.5){$1$}
\esegment
\move(-8 3)
\bsegment
\move(0 -0.2)\rlvec(0 5.2)\rlvec(1 0)\rlvec(0 -5.2)
\move(0 4)\rlvec(1 0)
\move(2 3.5)\rlvec(0 -3.7)
\move(5 1.5)\rlvec(0 -1.7)
\move(1 3)\rlvec(1 0)
\move(3 2.5)\rlvec(0 -2.7)
\move(4 2.5)\rlvec(0 -2.7)
\move(2 2)\rlvec(2 0)
\move(2 1)\rlvec(3 0)\rlvec(0 -1.2)
\move(2 0)\rlvec(3 0)
\move(0.5 4)\rlvec(0 1)
\move(2.5 2)\rlvec(0 0.5)
\move(3.5 2)\rlvec(0 0.5)
\move(2.5 0)\rlvec(0 1)
\move(3.5 0)\rlvec(0 1)
\move(4.5 0)\rlvec(0 1)
\lpatt(0.05 0.15)
\move(1 3.5)\rlvec(1 0)
\move(4 1.5)\rlvec(1 0)
\move(2 2.5)\rlvec(2 0)
\htext(0.25 4.5){$0$}
\htext(0.75 4.5){$1$}
\htext(1.5 3.25){$2$}
\htext(2.25 2.25){$3$}
\htext(2.75 2.25){$4$}
\htext(3.25 2.25){$3$}
\htext(3.75 2.25){$4$}
\htext(2.5 1.5){$2$}
\htext(3.5 1.5){$2$}
\htext(2.25 0.5){$0$}
\htext(2.75 0.5){$1$}
\htext(3.25 0.5){$0$}
\htext(3.75 0.5){$1$}
\htext(4.25 0.5){$0$}
\htext(4.75 0.5){$1$}
\htext(4.5 1.25){$2$}
\esegment
\move(-8 -3)
\bsegment
\move(0 -0.2)\rlvec(0 5.2)\rlvec(1 0)\rlvec(0 -5.2)
\move(0 4)\rlvec(1 0)
\move(2 3.5)\rlvec(0 -3.7)
\move(5 1.5)\rlvec(0 -1.7)
\move(1 3)\rlvec(2 0)
\move(3 3.5)\rlvec(0 -3.7)
\move(4 1.5)\rlvec(0 -1.7)
\move(2 2)\rlvec(1 0)
\move(2 1)\rlvec(3 0)\rlvec(0 -1.2)
\move(2 0)\rlvec(3 0)
\move(0.5 4)\rlvec(0 1)
\move(2.5 2)\rlvec(0 1)
\move(2.5 0)\rlvec(0 1)
\move(3.5 0)\rlvec(0 1)
\move(4.5 0)\rlvec(0 1)
\lpatt(0.05 0.15)
\move(1 3.5)\rlvec(2 0)
\move(3 1.5)\rlvec(2 0)
\htext(0.25 4.5){$0$}
\htext(0.75 4.5){$1$}
\htext(1.5 3.25){$2$}
\htext(2.5 3.25){$2$}
\htext(2.25 2.5){$3$}
\htext(2.75 2.5){$4$}
\htext(2.5 1.5){$2$}
\htext(3.5 1.25){$2$}
\htext(2.25 0.5){$0$}
\htext(2.75 0.5){$1$}
\htext(3.25 0.5){$0$}
\htext(3.75 0.5){$1$}
\htext(4.25 0.5){$0$}
\htext(4.75 0.5){$1$}
\htext(4.5 1.25){$2$}
\esegment
\move(8 3)
\bsegment
\move(0 -0.2)\rlvec(0 4.7)
\move(1 4.5)\rlvec(0 -4.7)
\move(0 4)\rlvec(2 0)
\move(2 4.5)\rlvec(0 -4.7)
\move(5 1)\rlvec(0 -1.2)
\move(1 3)\rlvec(1 0)
\move(3 2.5)\rlvec(0 -2.7)
\move(4 2.5)\rlvec(0 -2.7)
\move(2 2)\rlvec(2 0)
\move(2 1)\rlvec(3 0)\rlvec(0 -1.2)
\move(2 0)\rlvec(3 0)
\move(0.5 4)\rlvec(0 0.5)
\move(1.5 4)\rlvec(0 0.5)
\move(2.5 2)\rlvec(0 0.5)
\move(3.5 2)\rlvec(0 0.5)
\move(2.5 0)\rlvec(0 1)
\move(3.5 0)\rlvec(0 1)
\move(4.5 0)\rlvec(0 1)
\lpatt(0.05 0.15)
\move(0 4.5)\rlvec(2 0)
\move(2 2.5)\rlvec(2 0)
\htext(0.25 4.25){$0$}
\htext(0.75 4.25){$1$}
\htext(1.25 4.25){$0$}
\htext(1.75 4.25){$1$}
\htext(1.5 3.5){$2$}
\htext(2.25 2.25){$3$}
\htext(2.75 2.25){$4$}
\htext(3.25 2.25){$3$}
\htext(3.75 2.25){$4$}
\htext(2.5 1.5){$2$}
\htext(3.5 1.5){$2$}
\htext(2.25 0.5){$0$}
\htext(2.75 0.5){$1$}
\htext(3.25 0.5){$0$}
\htext(3.75 0.5){$1$}
\htext(4.25 0.5){$0$}
\htext(4.75 0.5){$1$}
\esegment
\move(8 -3)
\bsegment
\move(0 -0.2)\rlvec(0 4.7)
\move(1 4.5)\rlvec(0 -4.7)
\move(0 4)\rlvec(2 0)
\move(2 4.5)\rlvec(0 -4.7)
\move(5 1)\rlvec(0 -1.2)
\move(1 3)\rlvec(2 0)
\move(3 3.5)\rlvec(0 -3.7)
\move(4 1.5)\rlvec(0 -1.7)
\move(2 2)\rlvec(1 0)
\move(2 1)\rlvec(3 0)\rlvec(0 -1.2)
\move(2 0)\rlvec(3 0)
\move(0.5 4)\rlvec(0 0.5)
\move(1.5 4)\rlvec(0 0.5)
\move(2.5 2)\rlvec(0 1)
\move(2.5 0)\rlvec(0 1)
\move(3.5 0)\rlvec(0 1)
\move(4.5 0)\rlvec(0 1)
\lpatt(0.05 0.15)
\move(0 4.5)\rlvec(2 0)
\move(2 3.5)\rlvec(1 0)
\move(3 1.5)\rlvec(1 0)
\htext(0.25 4.25){$0$}
\htext(0.75 4.25){$1$}
\htext(1.25 4.25){$0$}
\htext(1.75 4.25){$1$}
\htext(1.5 3.5){$2$}
\htext(2.5 3.25){$2$}
\htext(2.25 2.5){$3$}
\htext(2.75 2.5){$4$}
\htext(2.5 1.5){$2$}
\htext(3.5 1.25){$2$}
\htext(2.25 0.5){$0$}
\htext(2.75 0.5){$1$}
\htext(3.25 0.5){$0$}
\htext(3.75 0.5){$1$}
\htext(4.25 0.5){$0$}
\htext(4.75 0.5){$1$}
\esegment

\move(5.8 2)\ravec(1.5 1)
\move(5.8 0.5)\ravec(1.5 -1)
\move(-0.8 2)\ravec(-1.5 1)
\move(-0.8 0.5)\ravec(-1.5 -1)
\end{texdraw}%
\end{center}

\vspace{2mm}

\item The following figure shows an incorrect $\zo$-block splitting
      of the above pre-slice.
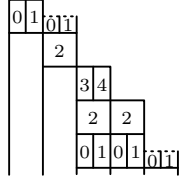
\begin{center}
\begin{texdraw}%
\drawdim in
\arrowheadsize l:0.065 w:0.03
\arrowheadtype t:F
\fontsize{7}{7}\selectfont
\textref h:C v:C
\drawdim em
\setunitscale 1.6
\move(0 -0.2)\rlvec(0 5.2)\rlvec(1 0)\rlvec(0 -5.2)
\move(0 4)\rlvec(2 0)\rlvec(0 -4.2)
\move(2 4)\rlvec(0 0.5)
\move(1 3)\rlvec(2 0)\rlvec(0 -3.2)
\move(2 2)\rlvec(2 0)\rlvec(0 -2.2)
\move(2 1)\rlvec(2 0)
\move(5 0.5)\rlvec(0 -0.7)
\move(2 0)\rlvec(3 0)
\move(0.5 4)\rlvec(0 1)
\move(2.5 2)\rlvec(0 1)
\move(2.5 0)\rlvec(0 1)
\move(3.5 0)\rlvec(0 1)
\move(4.5 0)\rlvec(0 0.5)
\move(1.5 4)\rlvec(0 0.5)
\lpatt(0.05 0.15)
\move(1 4.5)\rlvec(1 0)
\move(4 0.5)\rlvec(1 0)
\htext(0.25 4.5){$0$}
\htext(0.75 4.5){$1$}
\htext(1.5 3.5){$2$}
\htext(2.25 2.5){$3$}
\htext(2.75 2.5){$4$}
\htext(2.5 1.5){$2$}
\htext(3.5 1.5){$2$}
\htext(2.25 0.5){$0$}
\htext(2.75 0.5){$1$}
\htext(3.25 0.5){$0$}
\htext(3.75 0.5){$1$}
\htext(4.25 0.25){$0$}
\htext(4.75 0.25){$1$}
\htext(1.25 4.25){$0$}
\htext(1.75 4.25){$1$}
\end{texdraw}%
\end{center}
\end{enumerate}
\end{example}

\vspace{2mm}

As we can see in the above example, given a level-$l$ pre-slice,
it may be possible to split it in many different ways.
Any one of the (possibly many) results of splitting a given
pre-slice $C$,
which does not allow further splitting, is said to be
a \defi{split form} of $C$.
If $C$ does not allow any splitting at all, $C$ itself is the
split form of $C$.

We saw in the above example that
if we can split a covering $i$-block for
$i= 2,\dots, \nn$
(or a supporting $i$-block for $i=\zo, \dots, n-2$),
then it is always possible to split
a supporting $(i-1)$-block
(resp. covering $(i+1)$-block) from the same pre-slice, instead.

\begin{df}
Chose one of the two given patterns for stacking blocks
given in Definition~\ref{df3.1}.
The set consisting of all split forms of all level-$l$ pre-slices
of type $\dnone$ following the unique pattern will be denoted by $\nslice$.
The elements of $\nslice$ are called
\defi{level-$l$ slices} of type $\dnone$.
\end{df}

\begin{remark}
The theory on slices to be developed will not depend on
which of the two patterns one has chosen.
\end{remark}

\begin{remark}
We are not imposing the circular inclusion condition
given by~\eqref{eq:inclusion} to slices.
In fact, as the previous example of a splitting shows,
the circular inclusion condition now holds true only in a weaker
sense.
\end{remark}

\begin{remark}
The notions of \emph{$i$-slot}, \emph{$\delta$-column}, and
\emph{layer}, defined for pre-slices,
naturally carries over to those of slices.
Some care must be exercised, however.
For example, $\delta$-columns should allow for halves of blocks
to add up to a $\delta$,
and we should now consider halves of $i$-slots.
\end{remark}

Now we would like to caution the readers that the top of each layer is
an $i$-slot if the layer viewed as a level-$1$ slice contains
an $i$-slot at the top.
Those that are visually oriented might fall into the error
of viewing only the rear-most
$i$-slot among those of the same height as an $i$-slot.

We now define the action of Kashiwara operators on the set
$\nslice$ of slices.
Let $C$ be a level-$l$ slice and fix an index $i\in I$.

\noindent
[$\fit$ for $i = 2, \cdots, n-2$]\hfill\\
All slots considered below are $i$-slots.
We shall deal with the following four cases separately.
\begin{itemize}
\item[($\textup{F}_{\textup{two}}$)] $C$ contains covering full-slots
     (as opposed to halves of slots) and also
     contains supporting full-slots.
\item[($\textup{F}_{\textup{one}}$)] $C$ does contain full-slots but they are
     of one type, i.e., covering or supporting.
\item[(H)] $C$ contains no full-slots, but does contain
     half-slots.
     In such cases, there will be an equal number of
     covering half-slots and supporting half-slots.
\item[(N)] $C$ contains no full-slot
     or even a half-slot.
\end{itemize}
Break an $i$-block into two halves and, for each of the four case, position them
as instructed below.
\begin{itemize}
\item[($\textup{F}_{\textup{two}}$)]
     Place one half-block in the rear-most covering full-slot among
     those in the lower height.
     Do the same on supporting full-slots.
\item[($\textup{F}_{\textup{one}}$)]
     Place a half-block in the rear-most full-slot among those in the
     lower height.
     This will create an (upper) half-slot.
     Considering just the (upper) half-slots of the same
     (covering or supporting) type, place the other half-block in the
     rear-most half-slot among those in the lower height.
\item[(H)]
     Place a half-block in the rear-most covering half-slot among those
     in the lower height.
     Do the same in a supporting half-slot.
\item[(N)] We define $\fit C=0$.
\end{itemize}
[$\eit$ for $i = 2, \cdots, n-2$]\hfill\\
All blocks considered below are $i$-blocks that sit at the top of some
layer.
Replace the word \emph{slot} with \emph{block}
in the above division of classes to define four new cases,
which we denote by the same label.
For each case, remove two half-blocks as described below.
\begin{itemize}
\item[($\textup{F}_{\textup{two}}$)]
  Remove the upper half-block from the
  fore-most covering full-block among those in the higher height.
  Do the same from a supporting full-block.
\item[($\textup{F}_{\textup{one}}$)]
  Remove the upper half-block from the
  fore-most full-block among those in the higher height.
  This will create an (lower) half-block.
  Considering just the (lower) half-blocks of the same
  (covering or supporting) type, remove the fore-most half-block among
  those in the higher height.
\item[(H)]
  Remove the fore-most covering half-block among those
  in the higher height.
  Do the same with a supporting half-block.
  \item[(N)] We define $\eit C=0$.
\end{itemize}

\noindent
[$\fit$ for $i = 0, 1, n-1, n$]\hfill\\
For $i = 0, 1$, a block of unspecified color refers to a
$\zo$-block and for $i = n-1, n$, a block is an $\nn$-block.
First, un-split every possible block from $C$.
\begin{itemize}
\item[(F)]
  If the result contains $i$-slots,
  place an $i$-block in the rear-most $i$-slot
  among those in the lower height.
  Then split every block possible.
\item[(N)] If the result contains no $i$-slot,
  we define $\fit C=0$.
\end{itemize}
\noindent
[$\eit$ for $i = 0, 1, n-1, n$]\hfill\\
For $i = 0, 1$, a block of unspecified color refers to a
$\zo$-block, and for $i = n-1, n$, a block is an $\nn$-block.
Un-split every possible block.
\begin{itemize}
\item[(F)] If the result contains $i$-blocks, remove the
  fore-most $i$-block among those in the higher height.
  Then split every block possible.
\item[(N)] If the result contains no $i$-block,
  we define $\eit C=0$.
\end{itemize}

\begin{remark}
In a slice, unlike the case for other whole blocks,
the result of un-splitting all $\zo$ and $\nn$-blocks,
mentioned in above actions of Kashiwara operators for
$i = 0, 1, n-1, n$, is unique.
\end{remark}

Let $C= (c_1, \dots, c_l)$ be a level-$l$ slice.
We define the slices $C \pm \delta$ by
\begin{equation}\label{eq:add a delta}
\begin{aligned}
C+\delta & = (c_2, \dots, c_l, c_1 + \delta), \\
C-\delta & = (c_l - \delta, c_1, \cdots, c_{l-1}).
\end{aligned}
\end{equation}

\noindent
We say that two slices $C$ and $C'$ are \defi{related}, denoted by
$C\sim C'$, if one of the two slices may be obtained from the
other by adding finitely many $\delta$'s. Let
\begin{equation}\label{eq:rel}
\newcrystal = \nslice/\sim
\end{equation}
be the set of equivalence classes of level-$l$ slices.
We will use the same symbol $C$ for the equivalence class
containing the level-$l$ slice $C$.
By abuse of terminology,
the equivalence class containing a slice $C$
will be often referred to as the {\it slice} $C$.

Note that the map $C(\in \nslice) \mapsto C + \delta$ commutes with
the action of Kashiwara operators. Hence the Kashiwara operators
are well-defined on $\newcrystal$. We define
\begin{equation}
\begin{aligned}
\vphi_i ( C ) &= \max\{ n \mid \fit^n C \in \newcrystal \},\\
\veps_i ( C ) &= \max\{ n \mid \eit^n C \in \newcrystal \},\\
\cwt( C ) &= \sum_i \big(\vphi_i(C) - \veps_i(C)\big) \La_i.
\end{aligned}
\end{equation}
Then it is tedious but straightforward to prove the following
proposition.

\begin{prop}
The Kashiwara operators, together with the maps $\vphi_i$,
$\veps_i$ $(i\in I)$, $\cwt$, define a $\uq'(\g)$-crystal
structure on the set $\newcrystal$.
\end{prop}

%%%%%%%%%%%%%%%%%%%%%%%%%%%%%%%%%%%%%%%%%%%%%%%%%%%%%%%%%%%%%%%%%%%%%%%%%%
\section{New realization of perfect crystals}%
\label{sec:4}

In this section, we will show that the $\uq'(\g)$-crystal
$\newcrystal$ gives a new realization of the level-$l$ perfect
crystal $\oldcrystal$ described in Section~\ref{sec:2}. We first
define a \emph{canonical map} $\psi: \oldcrystal \rightarrow
\newcrystal$ as follows.

\vspace{2mm}

\noindent
Recall that every element $b \in \oldcrystal$ is of the form
$b=(x_1, \dots, x_n|\bar{x}_n, \cdots, \bar{x}_1)$
with $x_n=0 \textup{ or } \bar{x}_n=0$, $x_i, \bar{x}_i \in \Z_{\geq0}$,
and $\textstyle\sum_{i=1}^{n} (x_i + \bar{x}_i) = \textnormal{$l$}$.
Set
{\allowdisplaybreaks
\begin{equation}\label{eq:bnone}
\begin{aligned}
x'_1 &= \text{max} \{ 0, x_1-\bar x_1 \}, \\
& \qquad \ \vdots \\
x'_{n-1} &= \text{max} \{0, x_{n-1}-\bar x_{n-1} \}, \\
x'_n &= x_n, \\
\bar x'_n &= \bar x_n, \\
\bar x'_{n-1} &= \text{max} \{0, \bar x_{n-1}- x_{n-1} \}, \\
& \qquad \ \vdots \\
\bar x'_1 &= \text{max} \{ 0, \bar x_1-x_1 \}, \\
x^*_{\zo} &= \text{min} \{x_1, \bar x_1 \}, \\
x^*_2 &= \text{min} \{x_2, \bar x_2 \}, \\
& \qquad \ \vdots \\
x^*_{n-2} &= \text{min} \{x_{n-2}, \bar x_{n-2} \},\\
x^*_{\nn} &= \text{min} \{x_{n-1}, \bar x_{n-1} \}.
\end{aligned}
\end{equation}}% end of allowdisplaybreaks
For ease of writing, we shall temporarily use the following notation
to denote a gathering of blocks or halves of blocks stacked
in a single column of unit depth, i.e., something
which may appear as a layer of a slice.
Each notation is defined by its top part,
described by the figure to its right :
\begin{center}
\allowdisplaybreaks
\begin{align*}
s_1 \ &= \quad
\raisebox{-0.6em}{\begin{texdraw}
\fontsize{7}{7}\selectfont
\textref h:C v:C
\drawdim em
\setunitscale 1.9
\move(0 -0.2)\lvec(0 1)\lvec(1 1)\lvec(1 -0.2) \move(0 0)\lvec(1 0)
\move(0 0)\lvec(1 1) \move(0 0)\lvec(1 0) \htext(0.75 0.25){$0$}
\end{texdraw}%
}
 \quad \textup{ or } \quad
\raisebox{-0.6em}{\begin{texdraw}
\fontsize{7}{7}\selectfont
\textref h:C v:C
\drawdim em
\setunitscale 1.9
\move(0 -0.2)\lvec(0 1)\lvec(1 1)\lvec(1 -0.2)
\move(0 0)\lvec(1 0)
\move(0 0)\lvec(1 1)
\move(0 0)\lvec(1 0)
\htext(0.25 0.75){$0$}
\end{texdraw}}\\[0.3em]
s_2 \ &= \quad
\raisebox{-0.6em}{\begin{texdraw}
\fontsize{7}{7}\selectfont
\textref h:C v:C
\drawdim em
\setunitscale 1.9
\move(0 -0.2)\lvec(0 1)\lvec(1 1)\lvec(1 -0.2)
\move(0 0)\lvec(1 0)
\move(0 0)\lvec(1 1)
\move(0 0)\lvec(1 0)
\htext(0.75 0.25){$0$}
\htext(0.25 0.75){$1$}
\end{texdraw}}%
 \quad \textup{ or } \quad
\raisebox{-0.6em}{\begin{texdraw}
\fontsize{7}{7}\selectfont
\textref h:C v:C
\drawdim em
\setunitscale 1.9
\move(0 -0.2)\lvec(0 1)\lvec(1 1)\lvec(1 -0.2)
\move(0 0)\lvec(1 0)
\move(0 0)\lvec(1 1)
\move(0 0)\lvec(1 0)
\htext(0.25 0.75){$0$}
\htext(0.75 0.25){$1$}
\end{texdraw}}\\[0.3em]
s_i \ &= \quad
\raisebox{-0.8em}{\begin{texdraw}
\fontsize{7}{7}\selectfont
\textref h:C v:C
\drawdim em
\setunitscale 1.9
\move(0 -0.2)\lvec(0 2)\lvec(1 2)\lvec(1 -0.2)
\move(0 0)\rlvec(1 0)
\move(0 1)\rlvec(1 0)
\htext(0.5 0.5){$i\!\!-\!\!2$}
\htext(0.5 1.5){$i\!\!-\!\!1$}
\end{texdraw}}%
 \qquad \qquad \qquad \qquad \text{ for $i=3,\cdots,n-1$, }\\[0.3em]
s_n \ &= \quad
\raisebox{-0.8em}{\begin{texdraw}
\fontsize{7}{7}\selectfont
\textref h:C v:C
\drawdim em
\setunitscale 1.9
\move(0 -0.2)\lvec(0 2)\lvec(1 2)\lvec(1 -0.2)
\move(0 0)\rlvec(1 0)
\move(0 1)\rlvec(1 0)
\move(1 2)\rlvec(-0.5 -0.5)
\htext(0.5 0.5){$n\!\!-\!\!2$}
\htext(0.55 1.25){$n\!\!-\!\!1$}
\end{texdraw}}%
 \quad \textup{ or } \quad
\raisebox{-0.8em}{\begin{texdraw}
\fontsize{7}{7}\selectfont
\textref h:C v:C
\drawdim em
\setunitscale 1.9
\move(0 -0.2)\lvec(0 2)\lvec(1 2)\lvec(1 -0.2)
\move(0 0)\rlvec(1 0)
\move(0 1)\rlvec(1 0)
\move(0 1)\rlvec(0.5 0.5)
\htext(0.5 0.5){$n\!\!-\!\!2$}
\htext(0.45 1.75){$n\!\!-\!\!1$}
\end{texdraw}}%
\\[0.3em]
\bar{s}_n \ &= \quad
\raisebox{-0.6em}{\begin{texdraw}
\fontsize{7}{7}\selectfont
\textref h:C v:C
\drawdim em
\setunitscale 1.9
\move(0 -0.2)\lvec(0 1)\lvec(1 1)\lvec(1 -0.2)
\move(0 0)\rlvec(1 1)
\move(0 0)\rlvec(1 0)
\htext(0.25 0.75){$n$}
\end{texdraw}}%
 \quad \textup{ or } \quad
\raisebox{-0.6em}{\begin{texdraw}
\fontsize{7}{7}\selectfont
\textref h:C v:C
\drawdim em
\setunitscale 1.9
\move(0 -0.2)\lvec(0 1)\lvec(1 1)\lvec(1 -0.2)
\move(0 0)\rlvec(1 1)
\move(0 0)\rlvec(1 0)
\htext(0.75 0.25){$n$}
\end{texdraw}}\\[0.3em]
\bar{s}_{n\!-\!1} \ &= \quad
\raisebox{-0.6em}{\begin{texdraw}
\fontsize{7}{7}\selectfont
\textref h:C v:C
\drawdim em
\setunitscale 1.9
\move(0 -0.2)\lvec(0 1)\lvec(1 1)\lvec(1 -0.2)
\move(1 1)\rlvec(-0.5 -0.5)
\move(0 0)\rlvec(1 0)
\htext(0.25 0.75){$n$}
\htext(0.55 0.25){$n\!\!-\!\!1$}
\end{texdraw}}%
 \quad \textup{ or } \quad
\raisebox{-0.6em}{\begin{texdraw}
\fontsize{7}{7}\selectfont
\textref h:C v:C
\drawdim em
\setunitscale 1.9
\move(0 -0.2)\lvec(0 1)\lvec(1 1)\lvec(1 -0.2)
\move(0 0)\rlvec(0.5 0.5)
\move(0 0)\rlvec(1 0)
\htext(0.45 0.75){$n\!\!-\!\!1$}
\htext(0.75 0.25){$n$}
\end{texdraw}}\\[0.3em]
\bar{s}_i \ &= \quad
\raisebox{-0.8em}{\begin{texdraw}
\fontsize{7}{7}\selectfont
\textref h:C v:C
\drawdim em
\setunitscale 1.9
\move(0 -0.2)\lvec(0 2)\lvec(1 2)\lvec(1 -0.2)
\move(0 0)\rlvec(1 0)
\move(0 1)\rlvec(1 0)
\htext(0.5 0.5){$i\!\!+\!\!1$}
\htext(0.5 1.5){$i$}
\end{texdraw}}%
 \qquad \qquad \qquad \qquad \text{ for $i=2,\cdots,n-2$, }\\[0.3em]
\bar{s}_1 \ &= \quad
\raisebox{-0.6em}{\begin{texdraw}
\fontsize{7}{7}\selectfont
\textref h:C v:C
\drawdim em
\setunitscale 1.9
\move(0 -0.2)\lvec(0 1)\lvec(1 1)\lvec(1 -0.2)
\move(0 0)\lvec(1 0)
\move(0 0)\lvec(1 1)
\move(0 0)\lvec(1 0)
\htext(0.25 0.75){$1$}
\end{texdraw}}%
 \quad \textup{ or } \quad
\raisebox{-0.6em}{\begin{texdraw}
\fontsize{7}{7}\selectfont
\textref h:C v:C
\drawdim em
\setunitscale 1.9
\move(0 -0.2)\lvec(0 1)\lvec(1 1)\lvec(1 -0.2)
\move(0 0)\lvec(1 0)
\move(0 0)\lvec(1 1)
\move(0 0)\lvec(1 0)
\htext(0.75 0.25){$1$}
\end{texdraw}}\\[0.3mm]
t_{\zo}=\bar{t}_{\zo} \ &= \quad
  \raisebox{-0.8em}{\begin{texdraw}
  \fontsize{7}{7}\selectfont
  \textref h:C v:C
  \drawdim em
  \setunitscale 1.9
  \move(0 -0.2)\lvec(0 1.5)
  \move(1 1.5)\lvec(1 -0.2)
  \move(0 1)\rlvec(1 0)
  \move(0 0)\lvec(1 0)
  \move(0 0)\lvec(1 0)
  \move(0 1)\rlvec(1 0.5)
  \htext(0.8 1.17){$0$}
  \htext(0.2 1.33){$1$}
  \htext(0.5 0.5){$2$}
  \lpatt(0.05 0.15)
  \move(0 1.5)\rlvec(1 0)
\end{texdraw}}
 \quad \textup{ or } \quad
\raisebox{-0.8em}{\begin{texdraw}
\fontsize{7}{7}\selectfont
\textref h:C v:C
\drawdim em
\setunitscale 1.9
\move(0 -0.2)\lvec(0 1.5)
\move(1 1.5)\lvec(1 -0.2)
\move(0 1)\rlvec(1 0)
\move(0 0)\lvec(1 0)
\move(0 0)\lvec(1 0)
\move(0 1)\rlvec(1 0.5)
\htext(0.8 1.17){$1$}
\htext(0.2 1.33){$0$}
\htext(0.5 0.5){$2$}
\lpatt(0.05 0.15)
\move(0 1.5)\rlvec(1 0)
\end{texdraw}}\\[0.3em]
t_2 \ &= \quad
\raisebox{-0.8em}{\begin{texdraw}
\fontsize{7}{7}\selectfont
\textref h:C v:C
\drawdim em
\setunitscale 1.9
\move(0 -0.2)\lvec(0 1.5)
\move(1 1.5)\lvec(1 -0.2)
\move(0 1)\rlvec(1 0)
\move(0 0)\lvec(1 0)
\move(0 0)\lvec(1 0)
\move(0 0)\rlvec(1 1)
\htext(0.5 1.25){$2$}
\htext(0.75 0.25){$0$}
\htext(0.25 0.75){$1$}
\lpatt(0.05 0.15)
\move(0 1.5)\rlvec(1 0)
\end{texdraw}}
 \quad \textup{ or } \quad
\raisebox{-0.8em}{\begin{texdraw}
\fontsize{7}{7}\selectfont
\textref h:C v:C
\drawdim em
\setunitscale 1.9
\move(0 -0.2)\lvec(0 1.5)
\move(1 1.5)\lvec(1 -0.2)
\move(0 1)\rlvec(1 0)
\move(0 0)\lvec(1 0)
\move(0 0)\lvec(1 0)
\move(0 0)\rlvec(1 1)
\htext(0.5 1.25){$2$}
\htext(0.75 0.25){$1$}
\htext(0.25 0.75){$0$}
\lpatt(0.05 0.15)
\move(0 1.5)\rlvec(1 0)
\end{texdraw}}\\[0.3em]
t_i \ &= \quad
\raisebox{-0.8em}{\begin{texdraw}
\fontsize{7}{7}\selectfont
\textref h:C v:C
\drawdim em
\setunitscale 1.9
\move(0 -0.2)\lvec(0 1.5)
\move(1 1.5)\lvec(1 -0.2)
\move(0 1)\rlvec(1 0)
\move(0 0)\lvec(1 0)
\move(0 0)\lvec(1 0)
\htext(0.5 1.25){$i$}
\htext(0.5 0.5){$i\!\!-\!\!1$}
\lpatt(0.05 0.15)
\move(0 1.5)\rlvec(1 0)
\end{texdraw}}
 \qquad \qquad \qquad \qquad \text{ for $i=3,\cdots,n-2$, }\\[0.3em]
t_{\nn}=\bar{t}_{\nn} \ &= \quad
  \raisebox{-0.8em}{\begin{texdraw}
  \fontsize{7}{7}\selectfont
  \textref h:C v:C
  \drawdim em
  \setunitscale 1.9
  \move(0 -0.2)\lvec(0 1.5)
  \move(1 1.5)\lvec(1 -0.2)
  \move(0 1)\rlvec(1 0)
  \move(0 0)\lvec(1 0)
  \move(0 0)\lvec(1 0)
  \move(1 1.5)\rlvec(-0.5 -0.25)
  \htext(0.55 1.15){$n\!\!-\!\!1$}
  \htext(0.15 1.4){$n$}
  \htext(0.5 0.5){$n\!\!-\!\!2$}
  \lpatt(0.05 0.15)
  \move(0 1.5)\rlvec(1 0)
\end{texdraw}}
 \quad \textup{ or } \quad
\raisebox{-0.8em}{\begin{texdraw}
\fontsize{7}{7}\selectfont
\textref h:C v:C
\drawdim em
\setunitscale 1.9
\move(0 -0.2)\lvec(0 1.5)
\move(1 1.5)\lvec(1 -0.2)
\move(0 1)\rlvec(1 0)
\move(0 0)\lvec(1 0)
\move(0 0)\lvec(1 0)
\move(0 1)\rlvec(0.5 0.25)
\htext(0.85 1.1){$n$}
\htext(0.45 1.35){$n\!\!-\!\!1$}
\htext(0.5 0.5){$n\!\!-\!\!2$}
\lpatt(0.05 0.15)
\move(0 1.5)\rlvec(1 0)
\end{texdraw}}\\[0.3em]
\bar{t}_{n\!-\!2} \ &= \quad
  \raisebox{-0.8em}{\begin{texdraw}
  \fontsize{7}{7}\selectfont
  \textref h:C v:C
  \drawdim em
  \setunitscale 1.9
  \move(0 -0.2)\lvec(0 1.5)
  \move(1 1.5)\lvec(1 -0.2)
  \move(0 1)\rlvec(1 0)
  \move(0 0)\lvec(1 0)
  \move(0 0)\lvec(1 0)
  \move(1 1)\rlvec(-0.5 -0.5)
  \htext(0.55 0.25){$n\!\!-\!\!1$}
  \htext(0.25 0.75){$n$}
  \htext(0.5 1.25){$n\!\!-\!\!2$}
  \lpatt(0.05 0.15)
  \move(0 1.5)\rlvec(1 0)
\end{texdraw}}
 \quad \textup{ or } \quad
\raisebox{-0.8em}{\begin{texdraw}
\fontsize{7}{7}\selectfont
\textref h:C v:C
\drawdim em
\setunitscale 1.9
\move(0 -0.2)\lvec(0 1.5)
\move(1 1.5)\lvec(1 -0.2)
\move(0 1)\rlvec(1 0)
\move(0 0)\lvec(1 0)
\move(0 0)\lvec(1 0)
\move(0 0)\rlvec(0.5 0.5)
\htext(0.75 0.25){$n$}
\htext(0.45 0.75){$n\!\!-\!\!1$}
\htext(0.5 1.25){$n\!\!-\!\!2$}
\lpatt(0.05 0.15)
\move(0 1.5)\rlvec(1 0)
\end{texdraw}}\\[0.3em]
\bar{t}_i \ &= \quad
\raisebox{-0.8em}{\begin{texdraw}
\fontsize{7}{7}\selectfont
\textref h:C v:C
\drawdim em
\setunitscale 1.9
\move(0 -0.2)\lvec(0 1.5)
\move(1 1.5)\lvec(1 -0.2)
\move(0 1)\rlvec(1 0)
\move(0 0)\lvec(1 0)
\move(0 0)\lvec(1 0)
\htext(0.5 0.5){$i\!\!+\!\!1$}
\htext(0.5 1.25){$i$}
\lpatt(0.05 0.15)
\move(0 1.5)\rlvec(1 0)
\end{texdraw}}
 \qquad \qquad \qquad \qquad \text{ for $i=2,\cdots,n-3$. }
\end{align*}
\end{center}
Here, the top of $t_i$ and $\bar{t}_i$ are covered by a broken half of
a whole $i$-block.

Image of the above $b$ under the map $\psi$ is defined
to be the equivalence class of a level-$l$ slice obtained
by pasting together
\begin{equation}
\begin{cases}
\textnormal{$x'_i$-many $s_i$} & \quad
  i=1,\dots,n,\\
\textnormal{$\bar x'_i$-many $\bar{s}_i$} & \quad
  i=1,\dots,n,\\
\textnormal{$x^*_i$-many $t_i$} & \quad
  i=\zo,2,\dots,n-2,\nn,\\
\textnormal{$x^*_i$-many $\bar{t}_i$} & \quad
  i=\zo,2,\dots,n-2,\nn.
\end{cases}
\end{equation}
In particular, we are gathering together $2x^*_{\zo}$-many
$t_{\zo}=\bar{t}_{\zo}$'s and $2x^*_{\nn}$-many $t_{\nn}=\bar{t}_{\nn}$'s.
Note that this equivalence class does not depend on the
way we have pasted the $s_i$'s, $\bar{s}_i$'s, $t_i$'s, and $\bar{t}_i$'s
together, as long as the pasted result forms a slice.
It is easy to verify that the map $\psi$ is well-defined.

\vspace{1mm}

Now, we define a new map $\phi: \newcrystal \rightarrow \oldcrystal$
to see that the map $\psi$ is a bijection.
Any element of $\newcrystal$ may be obtained
by pasting together some number of $s_i$, $\bar{s}_i$, $t_i$,
and $\bar{t}_i$.
We consider the general element $C$ of $\newcrystal$
obtained by a pasting together of $y_i$-many $s_i$'s,
$\bar y_i$-many $\bar s_i$'s, and $y^*_i$-many $t_i$'s and $\bar t_i$'s,
for each $i$.
The no-more-splitting property of slices tell us that
$y_i \bar{y}_i=0$,
and we know
$2(y^*_{\zo} + \sum_{i=2}^{n-2} y^*_i + y^*_{\nn})
+ \sum_{i=1}^{n} (y_i + \bar{y}_i) = l$.
For such an element $C$, we define
$\phi(C)=(x_1, \dots, x_n | \bar x_n, \dots, \bar x_1)$, where
{\allowdisplaybreaks
\begin{equation}\label{eq:bnone}
\allowdisplaybreaks
\begin{aligned}
x_1 &= y_1 + y^*_{\zo},\\
x_2 &= y_2 + y^*_2,\\
& \ \ \vdots\\
x_{n-2} &= y_{n-2} + y^*_{n-2},\\
x_{n-1} &= y_{n-1} + y^*_{\nn},\\
x_n &= y_n,\\
\bar{x}_n &= \bar{y}_n, \\
\bar{x}_{n-1} &= \bar{y}_{n-1} + y^*_{\nn},\\
\bar{x}_{n-2} &= \bar{y}_{n-2} + y^*_{n-2},\\
& \ \ \vdots\\
\bar{x}_2 &= \bar{y}_2 + y^*_2\\
\bar{x}_1 &= \bar{y}_1 + y^*_{\zo}.
\end{aligned}
\end{equation}
It is now easy to verify that the map $\phi$ is well-defined
and that it is actually the inverse of $\psi$.

It is almost obvious that the maps $\cwt$, $\vphi_i$, and $\veps_i$
commute with $\psi$.
It still remains to show that $\psi$
commutes with the Kashiwara operators.
However, this may be done through a lengthy
but a straightforward case-by-case verification.
Therefore, we obtain a new realization of level-$l$ perfect
crystals as the set of equivalence classes of level-$l$ slices.

\begin{thm}\label{thm:37}
The map $\psi: \oldcrystal \rightarrow \newcrystal$ defined above
is an isomorphism of $\uqp(\mathfrak{g})$-crystals.
\end{thm}

We close this section with an example.

\begin{example}
The following is a drawing of a portion of the level-$2$ perfect crystal
for type $D_4^{(1)}$ given in terms of elements of $\newcrystal$.
We have omitted all arrows that connect elements presented here
with those not presented here.\\
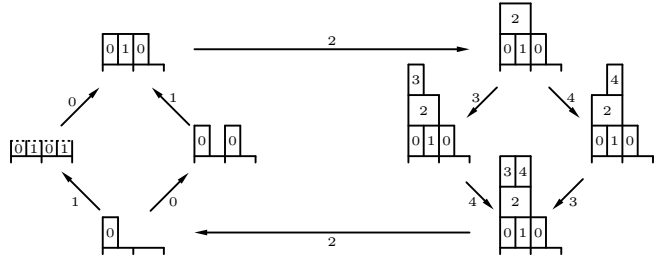
\begin{center}
\begin{texdraw}
\drawdim in
\arrowheadsize l:0.065 w:0.03
\arrowheadtype t:F
\fontsize{5}{5}\selectfont
\textref h:C v:C
\drawdim em
\setunitscale 1.7
\move(0 0)
 \bsegment
 \move(0 0)
  \bsegment
  \move(0 -0.2)\rlvec(0 0.7)
  \move(1 -0.2)\rlvec(0 0.7)
  \move(2 -0.2)\rlvec(0 0.7)
  \move(0.5 0)\rlvec(0 0.5)
  \move(1.5 0)\rlvec(0 0.5)
  \move(0 0)\rlvec(2 0)
  \lpatt(0.05 0.15)
  \move(0 0.5)\rlvec(2 0)
  \htext(0.25 0.25){$0$}
  \htext(0.75 0.25){$1$}
  \htext(1.25 0.25){$0$}
  \htext(1.75 0.25){$1$}
  \esegment
 \move(3 3)
  \bsegment
  \move(0 -0.2)\rlvec(0 1.2)
  \move(1 -0.2)\rlvec(0 1.2)
  \move(2 -0.2)\rlvec(0 0.2)
  \move(1.5 0)\rlvec(0 1)
  \move(0.5 0)\rlvec(0 1)
  \move(0 0)\rlvec(2 0)
  \move(0 1)\rlvec(1.5 0)
  \htext(0.25 0.5){$0$}
  \htext(0.75 0.5){$1$}
  \htext(1.25 0.5){$0$}
  \esegment
 \move(3 -3)
  \bsegment
  \move(0 -0.2)\rlvec(0 1.2)
  \move(1 -0.2)\rlvec(0 0.2)
  \move(2 -0.2)\rlvec(0 0.2)
  \move(0.5 0)\rlvec(0 1)
  \move(0 0)\rlvec(2 0)
  \move(0 1)\rlvec(0.5 0)
  \htext(0.25 0.5){$0$}
  \esegment
 \move(6 0)
  \bsegment
  \move(0 -0.2)\rlvec(0 1.2)
  \move(1 -0.2)\rlvec(0 1.2)
  \move(1.5 0)\rlvec(0 1)
  \move(2 -0.2)\rlvec(0 0.2)
  \move(0.5 0)\rlvec(0 1)
  \move(0 0)\rlvec(2 0)
  \move(0 1)\rlvec(0.5 0)
  \move(1 1)\rlvec(0.5 0)
  \htext(0.25 0.5){$0$}
  \htext(1.25 0.5){$0$}
  \esegment
 \move(1.65 1)\ravec(1.25 1.25)
 \htext(2 1.8){$0$}
 \move(2.9 -1.75)\ravec(-1.25 1.25)
 \htext(2.1 -1.5){$1$}
 \move(5.82 1.03)\ravec(-1.22 1.22)
 \htext(5.3 1.9){$1$}
 \move(4.57 -1.8)\ravec(1.25 1.25)
 \htext(5.3 -1.5){$0$}
 \esegment
\move(13 0)
 \bsegment
 \move(3 3)
  \bsegment
  \move(0 -0.2)\rlvec(0 2.2)
  \move(1 -0.2)\rlvec(0 2.2)
  \move(2 -0.2)\rlvec(0 0.2)
  \move(0 0)\rlvec(2 0)
  \move(0 2)\rlvec(1 0)
  \move(0 1)\rlvec(1.5 0)\rlvec(0 -1)
  \move(0.5 0)\rlvec(0 1)
  \htext(0.5 1.5){$2$}
  \htext(0.25 0.5){$0$}
  \htext(0.75 0.5){$1$}
  \htext(1.25 0.5){$0$}
  \esegment
 \move(0 0)
  \bsegment
  \move(0 -0.2)\rlvec(0 3.2)
  \move(1 -0.2)\rlvec(0 2.2)
  \move(2 -0.2)\rlvec(0 0.2)
  \move(0 0)\rlvec(2 0)
  \move(0.5 0)\rlvec(0 1)
  \move(1.5 0)\rlvec(0 1)
  \move(0 1)\rlvec(1.5 0)
  \move(0 2)\rlvec(1 0)
  \move(0 3)\rlvec(0.5 0)\rlvec(0 -1)
  \htext(0.25 2.5){$3$}
  \htext(0.25 0.5){$0$}
  \htext(1.25 0.5){$0$}
  \htext(0.75 0.5){$1$}
  \htext(0.5 1.5){$2$}
  \esegment
 \move(6 0)
  \bsegment
  \move(0 -0.2)\rlvec(0 2.2)
  \move(1 -0.2)\rlvec(0 3.2)
  \move(2 -0.2)\rlvec(0 0.2)
  \move(0 0)\rlvec(2 0)
  \move(0.5 0)\rlvec(0 1)
  \move(1.5 0)\rlvec(0 1)
  \move(0 1)\rlvec(1.5 0)
  \move(0 2)\rlvec(1 0)
  \move(1 3)\rlvec(-0.5 0)\rlvec(0 -1)
  \htext(0.75 2.5){$4$}
  \htext(0.25 0.5){$0$}
  \htext(1.25 0.5){$0$}
  \htext(0.75 0.5){$1$}
  \htext(0.5 1.5){$2$}
  \esegment
 \move(3 -3)
  \bsegment
  \move(0 -0.2)\rlvec(0 2.2)
  \move(1 -0.2)\rlvec(0 3.2)
  \move(2 -0.2)\rlvec(0 0.2)
  \move(0 0)\rlvec(2 0)
  \move(0.5 0)\rlvec(0 1)
  \move(1.5 0)\rlvec(0 1)
  \move(0 1)\rlvec(1.5 0)
  \move(0 2)\rlvec(1 0)
  \move(1 3)\rlvec(-0.5 0)\rlvec(0 -1)
  \move(0 2)\rlvec(0 1)\rlvec(0.5 0)
  \htext(0.75 2.5){$4$}
  \htext(0.25 2.5){$3$}
  \htext(0.25 0.5){$0$}
  \htext(1.25 0.5){$0$}
  \htext(0.75 0.5){$1$}
  \htext(0.5 1.5){$2$}
  \esegment
 \move(2.9 2.25)\ravec(-1 -1)
 \htext(2.2 1.9){$3$}
 \move(1.9 -0.85)\ravec(0.9 -0.9)
 \htext(2.1 -1.5){$4$}
 \move(4.6 2.25)\ravec(1 -1)
 \htext(5.3 1.9){$4$}
 \move(5.7 -0.75)\ravec(-1 -1)
 \htext(5.4 -1.5){$3$}
 \esegment

\move(6 3.5)\ravec(9 0)
\htext(10.5 3.8){$2$}
\move(15 -2.5)\ravec(-9 0)
\htext(10.5 -2.8){$2$}
\end{texdraw}
\end{center}
\end{example}

\vspace{2mm}

%%%%%%%%%%%%%%%%%%%%%%%%%%%%%%%%%%%%%%%%%%%%%%%%%%%%%%%%%%%%%%%%%%%%%%%%%%
\section{Combinatorics of Young walls}%
\label{sec:5}

In this section, we introduce \emph{level-$l$ Young walls}
for type $\dnone$.
Roughly speaking, these are constructed by
lining up, side by side, level-$l$ slices defined in Section~\ref{sec:3}.
In the next section, we will show that
the set of level-$l$ Young walls satisfying certain conditions
forms the highest weight crystal.

The pattern for building Young walls is given below.\\
\begin{equation}\label{ywallpattern}
\raisebox{-8em}{\begin{texdraw}
\fontsize{7}{7}\selectfont
\textref h:C v:C
\drawdim em
\setunitscale 1.9
\nc{\dtri}{
\bsegment
\esegment
}
\move(0 0)\dtri \move(-1 0)\dtri \move(-2 0)\dtri \move(-3 0)\dtri
\move(0 0)\rlvec(-4.3 0)
\move(0 1)\rlvec(-4.3 0)
\move(0 2)\rlvec(-4.3 0)
\move(0 3.5)\rlvec(-4.3 0)
\move(0 4.5)\rlvec(-4.3 0)
\move(0 5.5)\rlvec(-4.3 0)
\move(0 6.5)\rlvec(-4.3 0)
\move(0 8)\rlvec(-4.3 0)
\move(0 9)\rlvec(-4.3 0)
\move(0 10)\rlvec(-4.3 0)
\move(0 11)\rlvec(-4.3 0)
\move(0 0)\rlvec(0 11.3)
\move(-1 0)\rlvec(0 11.3)
\move(-2 0)\rlvec(0 11.3)
\move(-3 0)\rlvec(0 11.3)
\move(-4 0)\rlvec(0 11.3)
\move(-1 0)\rlvec(1 1)
\move(-2 0)\rlvec(1 1)
\move(-3 0)\rlvec(1 1)
\move(-4 0)\rlvec(1 1)
\move(-1 9)\rlvec(1 1)
\move(-2 9)\rlvec(1 1)
\move(-3 9)\rlvec(1 1)
\move(-4 9)\rlvec(1 1)
\move(0 5.5)\rlvec(-0.5 -0.5)
\move(-2 4.5)\rlvec(0.5 0.5)
\move(-2 5.5)\rlvec(-0.5 -0.5)
\move(-4 4.5)\rlvec(0.5 0.5)
\htext(-0.3 0.25){$0$} \htext(-0.75 0.75){$1$}
\htext(-0.5 1.5){$2$}
\vtext(-0.5 2.75){$\cdots$}
\htext(-0.5 4){$n\!\!-\!\!2$}
\htext(-0.5 6){$n\!\!-\!\!2$}
\htext(-0.5 8.5){$2$}
\htext(-0.3 9.25){$0$} \htext(-0.75 9.75){$1$}
\htext(-0.5 10.5){$2$}
\htext(-2.3 0.25){$0$} \htext(-2.75 0.75){$1$}
\htext(-2.5 1.5){$2$}
\vtext(-2.5 2.75){$\cdots$}
\htext(-2.5 4){$n\!\!-\!\!2$}
\htext(-2.5 6){$n\!\!-\!\!2$}
\htext(-2.5 8.5){$2$}
\htext(-2.3 9.25){$0$} \htext(-2.75 9.75){$1$}
\htext(-2.5 10.5){$2$}
\htext(-1.3 0.25){$1$} \htext(-1.75 0.75){$0$}
\htext(-1.5 1.5){$2$}
\vtext(-1.5 2.75){$\cdots$}
\htext(-1.5 4){$n\!\!-\!\!2$}
\htext(-1.5 6){$n\!\!-\!\!2$}
\htext(-1.5 8.5){$2$}
\htext(-1.3 9.25){$1$} \htext(-1.75 9.75){$0$}
\htext(-1.5 10.5){$2$}
\htext(-3.3 0.25){$1$} \htext(-3.75 0.75){$0$}
\htext(-3.5 1.5){$2$}
\vtext(-3.5 2.75){$\cdots$}
\htext(-3.5 4){$n\!\!-\!\!2$}
\htext(-3.5 6){$n\!\!-\!\!2$}
\htext(-3.5 8.5){$2$}
\htext(-3.3 9.25){$1$} \htext(-3.75 9.75){$0$}
\htext(-3.5 10.5){$2$}

\htext(-0.45 4.75){$n\!\!-\!\!1$}
\htext(-0.75 5.25){$n$}
\htext(-1.25 4.75){$n$}
\htext(-1.55 5.25){$n\!\!-\!\!1$}
\htext(-2.45 4.75){$n\!\!-\!\!1$}
\htext(-2.75 5.25){$n$}
\htext(-3.25 4.75){$n$}
\htext(-3.55 5.25){$n\!\!-\!\!1$}

\vtext(-0.5 7.25){$\cdots$} \vtext(-1.5 7.25){$\cdots$}
\vtext(-2.5 7.25){$\cdots$} \vtext(-3.5 7.25){$\cdots$}
\end{texdraw}}
\end{equation}

\vspace{2mm}

\begin{df}
A \defi{level-$l$ Young wall of type $\dnone$} is a
concatenation of level-$l$ slices, extending (possibly infinitely)
to the left,
that satisfies the following conditions.
\begin{enumerate}
\item It is stacked in the pattern given above.
\item When we look at the same layer of two consecutive level-$l$ slices,
    there is no free space to the right of any block belonging to the
    left slice.
\end{enumerate}
\end{df}

\begin{remark}\label{expl}
The second condition in the above definition for
level-$l$ Young walls needs more explanation.
We do this here.
\begin{enumerate}
\item When we are dealing only with unbroken blocks,
      what we mean by having a free space to the right is trivial.
      The following lists all nontrivial cases of having
      a free space to the right.
\savebox{\tmpfigb}{\begin{texdraw}
\fontsize{7}{7}\selectfont
\textref h:C v:C
\drawdim em
\setunitscale 1.9
\move(0 -0.2)\lvec(0 2)\lvec(1 2)\lvec(1 -0.2)
\move(0 1)\lvec(1 1)
\move(0 1)\rlvec(1 1)
\move(1 -0.2)\rlvec(0 1.7)
\move(2 -0.2)\rlvec(0 1.7)
\move(1 1)\rlvec(1 0)
\lpatt(0.05 0.15)
\move(1 1.5)\rlvec(1 0)
\htext(0.25 1.75){$n$}
\htext(0.75 1.25){$n\!\!-\!\!1$}
\end{texdraw}%
}%
\savebox{\tmpfigc}{\begin{texdraw}
\fontsize{7}{7}\selectfont
\textref h:C v:C
\drawdim em
\setunitscale 1.9
\move(0 -0.2)\lvec(0 2)\lvec(1 2)\lvec(1 -0.2)
\move(0 1)\lvec(1 1)
\move(0 1)\rlvec(1 1)
\move(1 -0.2)\rlvec(0 1.7)
\move(2 -0.2)\rlvec(0 1.7)
\move(1 1)\rlvec(1 0)
\lpatt(0.05 0.15)
\move(1 1.5)\rlvec(1 0)
\htext(0.75 1.25){$n$}
\htext(0.25 1.75){$n\!\!-\!\!1$}
\end{texdraw}%
}%
\begin{itemize}
\item The left is a whole block and the right is a broken half of a
      whole block.\\[0.3mm]
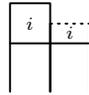
\begin{center}
\begin{texdraw}
\fontsize{7}{7}\selectfont
\textref h:C v:C
\drawdim em
\setunitscale 1.9
\move(0 -0.2)\lvec(0 2)\lvec(1 2)\lvec(1 -0.2)
\move(0 1)\lvec(1 1)
\move(1 -0.2)\rlvec(0 1.7)
\move(2 -0.2)\rlvec(0 1.7)
\move(1 1)\rlvec(1 0)
\lpatt(0.05 0.15)
\move(1 1.5)\rlvec(1 0)
\htext(0.5 1.5){$i$}
\htext(1.5 1.25){$i$}
\end{texdraw}\\[2mm]%
\end{center}
with $i = \zo, 2, \dots, n-2, \nn$.
\item Recall that unlike other whole blocks, for $j=\zo$ and $\nn$,
      there is a unique way to \emph{un-split} all broken halves of
      whole $j$-blocks contained in a slice.
      Hence it makes sense to talk of whether a broken half of a
      whole $j$-block is the upper half or the lower half.

      The next nontrivial case is when the left is a single block
      of half-unit depth and the right is the \emph{upper} broken half of
      a whole $j$-block.\\[0.3mm]
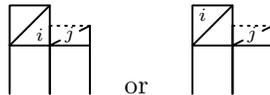
\begin{center}
\begin{texdraw}
\fontsize{7}{7}\selectfont
\textref h:C v:C
\drawdim em
\setunitscale 1.9
\move(0 -0.2)\lvec(0 2)\lvec(1 2)\lvec(1 -0.2)
\move(0 1)\lvec(1 1)
\move(0 1)\rlvec(1 1)
\move(1 -0.2)\rlvec(0 1.7)
\move(2 -0.2)\rlvec(0 1.7)
\move(1 1)\rlvec(1 0)
\move(1 1)\rlvec(0.25 0.125)
\move(2 1.5)\rlvec(-0.25 -0.125)
\lpatt(0.05 0.15)
\move(1 1.5)\rlvec(1 0)
\htext(0.75 1.25){$i$}
\htext(1.5 1.25){$j$}
\end{texdraw}%
\quad \text{ or } \quad
\begin{texdraw}
\fontsize{7}{7}\selectfont
\textref h:C v:C
\drawdim em
\setunitscale 1.9
\move(0 -0.2)\lvec(0 2)\lvec(1 2)\lvec(1 -0.2)
\move(0 1)\lvec(1 1)
\move(0 1)\rlvec(1 1)
\move(1 -0.2)\rlvec(0 1.7)
\move(2 -0.2)\rlvec(0 1.7)
\move(1 1)\rlvec(1 0)
\move(1 1)\rlvec(0.25 0.125)
\move(2 1.5)\rlvec(-0.25 -0.125)
\lpatt(0.05 0.15)
\move(1 1.5)\rlvec(1 0)
\htext(0.25 1.75){$i$}
\htext(1.5 1.25){$j$}
\end{texdraw}\\[2mm]%
\end{center}
      Here, $i = 0, 1, n-1, n$ and $j= \zo,\nn$.
      Note that since these must be blocks stacked following the pattern
      for building Young walls, the color of the broken half on the right
      will depend on the color $i$.
\item The right is a single block of half-unit depth and the left
      is the \emph{lower} broken half of a whole $j$-block.\\[0.5mm]
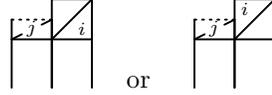
\begin{center}
\begin{texdraw}
\fontsize{7}{7}\selectfont
\textref h:C v:C
\drawdim em
\setunitscale 1.9
\move(1 -0.2)\lvec(1 2)\lvec(2 2)\lvec(2 -0.2)
\move(1 1)\rlvec(1 0)
\move(1 1)\rlvec(1 1)
\move(0 -0.2)\rlvec(0 1.7)
\move(1 -0.2)\rlvec(0 1.7)
\move(0 1)\rlvec(1 0)
\move(0 1)\rlvec(0.25 0.125)
\move(1 1.5)\rlvec(-0.25 -0.125)
\htext(0.5 1.25){$j$}
\lpatt(0.05 0.15)
\move(0 1.5)\rlvec(1 0)
\htext(1.75 1.25){$i$}
\end{texdraw}%
\quad \text{ or } \quad
\begin{texdraw}
\fontsize{7}{7}\selectfont
\textref h:C v:C
\drawdim em
\setunitscale 1.9
\move(1 -0.2)\lvec(1 2)\lvec(2 2)\lvec(2 -0.2)
\move(1 1)\rlvec(1 0)
\move(1 1)\rlvec(1 1)
\move(0 -0.2)\rlvec(0 1.7)
\move(1 -0.2)\rlvec(0 1.7)
\move(0 1)\rlvec(1 0)
\move(0 1)\rlvec(0.25 0.125)
\move(1 1.5)\rlvec(-0.25 -0.125)
\htext(0.5 1.25){$j$}
\lpatt(0.05 0.15)
\move(0 1.5)\rlvec(1 0)
\htext(1.25 1.75){$i$}
\end{texdraw}%
\end{center}
\vspace{1mm}
Here, $i = 0, 1, n-1, n$ and $j= \zo,\nn$.
\end{itemize}

\item
We state explicitly that the following can appear as
some layer of two consecutive level-$l$ slices
that make up a Young wall.\\
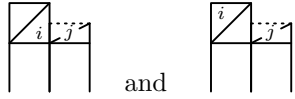
\begin{center}
\begin{texdraw}
\fontsize{7}{7}\selectfont
\textref h:C v:C
\drawdim em
\setunitscale 1.9
\move(0 -0.2)\lvec(0 2)\lvec(1 2)\lvec(1 -0.2)
\move(0 1)\lvec(1 1)
\move(0 1)\rlvec(1 1)
\move(1 -0.2)\rlvec(0 1.7)
\move(2 -0.2)\rlvec(0 1.7)
\move(1 1)\rlvec(1 0)
\move(1 1)\rlvec(0.25 0.125)
\move(2 1.5)\rlvec(-0.25 -0.125)
\htext(1.5 1.25){$j$}
\lpatt(0.05 0.15)
\move(1 1.5)\rlvec(1 0)
\htext(0.75 1.25){$i$}
\end{texdraw}%
\quad \text{ and } \quad
\begin{texdraw}
\fontsize{7}{7}\selectfont
\textref h:C v:C
\drawdim em
\setunitscale 1.9
\move(0 -0.2)\lvec(0 2)\lvec(1 2)\lvec(1 -0.2)
\move(0 1)\lvec(1 1)
\move(0 1)\rlvec(1 1)
\move(1 -0.2)\rlvec(0 1.7)
\move(2 -0.2)\rlvec(0 1.7)
\move(1 1)\rlvec(1 0)
\move(1 1)\rlvec(0.25 0.125)
\move(2 1.5)\rlvec(-0.25 -0.125)
\htext(1.5 1.25){$j$}
\lpatt(0.05 0.15)
\move(1 1.5)\rlvec(1 0)
\htext(0.25 1.75){$i$}
\end{texdraw}%
\end{center}
\vspace{1mm}
with $i = 0, 1, n-1, n$ and $j= \zo,\nn$ and where
the right is a lower half.\\
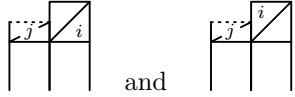
\begin{center}
\begin{texdraw}
\fontsize{7}{7}\selectfont
\textref h:C v:C
\drawdim em
\setunitscale 1.9
\move(1 -0.2)\lvec(1 2)\lvec(2 2)\lvec(2 -0.2)
\move(1 1)\rlvec(1 0)
\move(1 1)\rlvec(1 1)
\move(0 -0.2)\rlvec(0 1.7)
\move(1 -0.2)\rlvec(0 1.7)
\move(0 1)\rlvec(1 0)
\move(0 1)\rlvec(0.25 0.125)
\move(1 1.5)\rlvec(-0.25 -0.125)
\htext(0.5 1.25){$j$}
\lpatt(0.05 0.15)
\move(0 1.5)\rlvec(1 0)
\htext(1.75 1.25){$i$}
\end{texdraw}%
\quad \text{ and } \quad
\begin{texdraw}
\fontsize{7}{7}\selectfont
\textref h:C v:C
\drawdim em
\setunitscale 1.9
\move(1 -0.2)\lvec(1 2)\lvec(2 2)\lvec(2 -0.2)
\move(1 1)\rlvec(1 0)
\move(1 1)\rlvec(1 1)
\move(0 -0.2)\rlvec(0 1.7)
\move(1 -0.2)\rlvec(0 1.7)
\move(0 1)\rlvec(1 0)
\move(0 1)\rlvec(0.25 0.125)
\move(1 1.5)\rlvec(-0.25 -0.125)
\htext(0.5 1.25){$j$}
\lpatt(0.05 0.15)
\move(0 1.5)\rlvec(1 0)
\htext(1.25 1.75){$i$}
\end{texdraw}%
\end{center}
\vspace{1mm}
with $i = 0, 1, n-1, n$ and $j= \zo,\nn$ and where
the left is an upper half.
\end{enumerate}
\end{remark}

\begin{df} \hfill
\begin{enumerate}
\item A level-$l$ Young wall $\mathbf{Y}$ is said to be
      \defi{proper} if for each layer of $\mathbf{Y}$,
      none of the columns which is of integer height
      and whose top is of unit depth have the same height.
\item A column in a level-$l$ proper Young wall is said to contain a
      \defi{removable $\delta$} if one may remove a $\delta$ from that
      column and still obtain a proper Young wall.
\item A level-$l$ proper Young wall is said to be \defi{reduced}
      if none of its columns contain a removable $\delta$.
\end{enumerate}
\end{df}

We will denote by $\pwspace$ and $\rpwspace$ the set of all
level-$l$ proper Young walls and reduced proper Young walls,
respectively.

Let $\mathbf{Y}$ be a level-$l$ proper Young wall
and let $C$ be a column of $\mathbf{Y}$.
Recall that a column $C$ is a level-$l$ slice and
that for each $i\in I$, $\vphi_i(C)$ (and $\veps_i(C)$)
is the largest integer $k \ge 0$
such that $\fit^k (C) \neq 0$ (resp. $\eit^k (C) \neq 0$)
in $\newcrystal$.

We now define the action of Kashiwara operators $\fit$, $\eit$
$(i\in I)$ on $\mathbf{Y}$ as follows.

\begin{enumerate}
\item For each column $C$ of $\mathbf{Y}$, write
      $\veps_i(C)$-many $1$'s followed by
      $\vphi_i(C)$-many $0$'s under $C$. This sequence is
      called the \defi{$i$-signature of $C$}.
\item From this sequence of 1's and 0's,
      cancel out each (0,1)-pair to obtain a sequence of 1's
      followed by 0's (reading from left to right). This sequence is
      called the \defi{$i$-signature of $\mathbf{Y}$}.
\item We define $\fit \mathbf{Y}$ to be the proper Young wall
      obtained from $\mathbf{Y}$ by replacing the column $C$
      corresponding the leftmost 0 in the $i$-signature of $\mathbf{Y}$
      with the column $\fit C$.
\item We define $\eit \mathbf{Y}$ to be the proper Young wall
      obtained from $\mathbf{Y}$ by replacing the column $C$
      corresponding the rightmost $1$ in the $i$-signature of
      $\mathbf{Y}$ with the column $\eit C$.
\item If there is no $0$ (or $1$) in the $i$-signature of
      $\mathbf{Y}$, we define $\fit \mathbf{Y} = 0$ (resp. $\eit
      \mathbf{Y} = 0$).
\end{enumerate}
\begin{proof}[Sketch of proof for well-definedness]
Let us now show that the above definition of Kashiwara operators
on proper Young walls is well-defined.

First, we fix some notations.
Denote by $C$ the column corresponding to the leftmost $0$
in the $i$-signature of a proper Young wall $\mathbf{Y}$.
The column sitting to the right of column $C$ will be denoted by $C'$.

Suppose that $\fit\mathbf{Y}$ is not a proper Young wall.
In fact, it could be that $\fit\mathbf{Y}$ is not even a Young wall.
In such a case, the following statement would be true.
\begin{itemize}
\item
There is a free space to the right of some block or half-block
in some layer of the two columns of $\fit\mathbf{Y}$,
that correspond to columns $C$ and $C'$ of $\mathbf{Y}$.
\end{itemize}
If the result is a Young wall, but just not proper, then
the following statement would be true.
\begin{itemize}
\item
The columns of $\fit \bf{Y}$ corresponding to columns $C$ and $C'$
of $\mathbf{Y}$
contain a layer in which the tops are of unit depth and of
the same integer height.
\end{itemize}

When $i = 2, \cdots, n-2$, the following lists all possible
forms for $\fit Y$ which satisfies one of the above two statements.
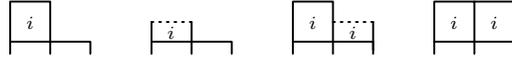
\savebox{\tmpfiga}{\begin{texdraw}
\fontsize{7}{7}\selectfont
\textref h:C v:C
\drawdim em
\setunitscale 1.9
\move(0 -0.3)\lvec(0 1)\lvec(1 1)\lvec(1 -0.3)
\move(0 0)\rlvec(2 0)\rlvec(0 -0.3)
\htext(0.5 0.5){$i$}
\end{texdraw}%
}%
\savebox{\tmpfigb}{\begin{texdraw}
\fontsize{7}{7}\selectfont
\textref h:C v:C
\drawdim em
\setunitscale 1.9
\move(0 -0.3)\lvec(0 0.5)
\move(1 -0.3)\lvec(1 0.5)
\move(0 0)\lvec(2 0)\rlvec(0 -0.3)
\lpatt(0.05 0.15)
\move(0 0.5)\rlvec(1 0)
\htext(0.5 0.25){$i$}
\end{texdraw}%
}%
\savebox{\tmpfigc}{\begin{texdraw}
\fontsize{7}{7}\selectfont
\textref h:C v:C
\drawdim em
\setunitscale 1.9
\move(0 -0.3)\lvec(0 1)\lvec(1 1)\lvec(1 -0.3)
\move(0 0)\rlvec(2 0)\rlvec(0 -0.3)
\move(2 -0.2)\lvec(2 0.5)
\lpatt(0.05 0.15)
\move(1 0.5)\rlvec(1 0)
\htext(0.5 0.5){$i$}
\htext(1.5 0.25){$i$}
\end{texdraw}%
}%
\savebox{\tmpfigd}{\begin{texdraw}
\fontsize{7}{7}\selectfont
\textref h:C v:C
\drawdim em
\setunitscale 1.9
\move(0 -0.3)\lvec(0 1)\lvec(1 1)\lvec(1 -0.3)
\move(0 0)\rlvec(2 0)\rlvec(0 -0.3)
\move(2 -0.2)\lvec(2 1)\rlvec(-1 0)
\htext(0.5 0.5){$i$}
\htext(1.5 0.5){$i$}
\end{texdraw}%
}%
\begin{center}
\usebox{\tmpfiga}
\qquad
\usebox{\tmpfigb}
\qquad
\usebox{\tmpfigc}
\qquad
\usebox{\tmpfigd}
\end{center}
When $i = 0,1,n-1,n$, the following lists all possible
forms for $\fit Y$ which satisfies one of the above two statements.
\savebox{\tmpfiga}{\begin{texdraw}
\fontsize{7}{7}\selectfont
\textref h:C v:C
\drawdim em
\setunitscale 1.9
\move(0 -0.3)\lvec(0 1)\lvec(1 1)\lvec(1 -0.3)
\move(0 0)\rlvec(2 0)\rlvec(0 -0.3)
\move(0 0)\rlvec(1 1)
\htext(0.75 0.25){$i$}
\end{texdraw}%
}%
\savebox{\tmpfigb}{\begin{texdraw}
\fontsize{7}{7}\selectfont
\textref h:C v:C
\drawdim em
\setunitscale 1.9
\move(0 -0.3)\lvec(0 1)\lvec(1 1)\lvec(1 -0.3)
\move(0 0)\rlvec(2 0)\rlvec(0 -0.3)
\move(0 0)\rlvec(1 1)
\htext(0.25 0.75){$i$}
\end{texdraw}%
}%
\savebox{\tmpfigc}{\begin{texdraw}
\fontsize{7}{7}\selectfont
\textref h:C v:C
\drawdim em
\setunitscale 1.9
\move(0 -0.3)\lvec(0 1)\lvec(1 1)\lvec(1 -0.3)
\move(0 0)\rlvec(2 0)\rlvec(0 -0.3)
\move(2 -0.2)\lvec(2 1)\rlvec(-1 0)
\move(0 0)\rlvec(1 1)
\move(1 0)\rlvec(1 1)
\htext(0.75 0.25){$i$}
\htext(1.25 0.75){$i$}
\end{texdraw}%
}%
\savebox{\tmpfigd}{\begin{texdraw}
\fontsize{7}{7}\selectfont
\textref h:C v:C
\drawdim em
\setunitscale 1.9
\move(0 -0.3)\lvec(0 1)\lvec(1 1)\lvec(1 -0.3)
\move(0 0)\rlvec(2 0)\rlvec(0 -0.3)
\move(2 -0.2)\lvec(2 1)\rlvec(-1 0)
\move(0 0)\rlvec(1 1)
\move(1 0)\rlvec(1 1)
\htext(0.25 0.75){$i$}
\htext(1.75 0.25){$i$}
\end{texdraw}%
}%
\savebox{\tmpfige}{\begin{texdraw}
\fontsize{7}{7}\selectfont
\textref h:C v:C
\drawdim em
\setunitscale 1.9
\move(0 -0.3)\lvec(0 1)\lvec(1 1)\lvec(1 -0.3)
\move(0 0)\rlvec(2 0)\rlvec(0 -0.3)
\move(0 0)\rlvec(1 1)
\move(2 0)\lvec(2 0.5)
\move(1 0)\rlvec(0.25 0.125)
\move(2 0.5)\rlvec(-0.25 -0.125)
\htext(0.75 0.25){$i$}
\htext(1.5 0.25){$j$}
\lpatt(0.05 0.15)
\move(1 0.5)\rlvec(1 0)
\end{texdraw}%
}%
\savebox{\tmpfigf}{\begin{texdraw}
\fontsize{7}{7}\selectfont
\textref h:C v:C
\drawdim em
\setunitscale 1.9
\move(0 -0.3)\lvec(0 1)\lvec(1 1)\lvec(1 -0.3)
\move(0 0)\rlvec(2 0)\rlvec(0 -0.3)
\move(0 0)\rlvec(1 1)
\move(2 0)\lvec(2 0.5)
\htext(0.25 0.75){$i$}
\move(1 0)\rlvec(0.25 0.125)
\move(2 0.5)\rlvec(-0.25 -0.125)
\htext(1.5 0.25){$j$}
\lpatt(0.05 0.15)
\move(1 0.5)\rlvec(1 0)
\end{texdraw}%
}%
\savebox{\tmpfigg}{\begin{texdraw}
\fontsize{7}{7}\selectfont
\textref h:C v:C
\drawdim em
\setunitscale 1.9
\move(0 -0.3)\rlvec(0 0.8)
\move(1 -0.3)\rlvec(0 0.8)
\move(0 0)\rlvec(2 0)\rlvec(0 -0.3)
\move(2 0)\lvec(2 0)
\move(0 0)\rlvec(0.25 0.125)
\move(1 0.5)\rlvec(-0.25 -0.125)
\htext(0.5 0.25){$j$}
\lpatt(0.05 0.15)
\move(0 0.5)\rlvec(1 0)
\end{texdraw}%
}%
\savebox{\tmpfigh}{\begin{texdraw}
\fontsize{7}{7}\selectfont
\textref h:C v:C
\drawdim em
\setunitscale 1.9
\move(0 -0.3)\rlvec(0 0.8)
\move(1 -0.3)\rlvec(0 0.8)
\move(0 0)\rlvec(2 0)\rlvec(0 -0.3)
\move(2 0)\rlvec(0 1)\rlvec(-1 0)\rlvec(0 -0.5)
\move(2 1)\rlvec(-1 -1)
\move(0 0)\rlvec(0.25 0.125)
\move(1 0.5)\rlvec(-0.25 -0.125)
\htext(0.5 0.25){$j$}
\htext(1.75 0.25){$i$}
\lpatt(0.05 0.15)
\move(0 0.5)\rlvec(1 0)
\end{texdraw}%
}%
\savebox{\tmpfigi}{\begin{texdraw}
\fontsize{7}{7}\selectfont
\textref h:C v:C
\drawdim em
\setunitscale 1.9
\move(0 -0.3)\rlvec(0 0.8)
\move(1 -0.3)\rlvec(0 0.8)
\move(0 0)\rlvec(2 0)\rlvec(0 -0.3)
\move(2 0)\rlvec(0 1)\rlvec(-1 0)\rlvec(0 -0.5)
\move(2 1)\rlvec(-1 -1)
\move(0 0)\rlvec(0.25 0.125)
\move(1 0.5)\rlvec(-0.25 -0.125)
\htext(0.5 0.25){$j$}
\htext(1.25 0.75){$i$}
\lpatt(0.05 0.15)
\move(0 0.5)\rlvec(1 0)
\end{texdraw}%
}%
\savebox{\tmpfigm}{\begin{texdraw}
\fontsize{7}{7}\selectfont
\textref h:C v:C
\drawdim em
\setunitscale 1.9
\move(0 -0.3)\lvec(0 1)\lvec(1 1)\lvec(1 -0.3)
\move(0 0)\rlvec(2 0)\rlvec(0 -0.3)
\move(0 0)\rlvec(0.25 0.25)
\move(1 1)\rlvec(-0.25 -0.25)
\htext(0.5 0.5){$j$}
\end{texdraw}
}%
\savebox{\tmpfigj}{\begin{texdraw}
\fontsize{7}{7}\selectfont
\textref h:C v:C
\drawdim em
\setunitscale 1.9
\move(0 -0.3)\lvec(0 1)\lvec(1 1)\lvec(1 -0.3)
\move(0 0)\rlvec(2 0)\rlvec(0 -0.3)
\move(2 0)\rlvec(0 1)\rlvec(-1 0)
\move(2 1)\rlvec(-1 -1)
\htext(1.75 0.25){$i$}
\move(0 0)\rlvec(0.25 0.25)
\move(1 1)\rlvec(-0.25 -0.25)
\htext(0.5 0.5){$j$}
\end{texdraw}
}%
\savebox{\tmpfigk}{\begin{texdraw}
\fontsize{7}{7}\selectfont
\textref h:C v:C
\drawdim em
\setunitscale 1.9
\move(0 -0.3)\lvec(0 1)\lvec(1 1)\lvec(1 -0.3)
\move(0 0)\rlvec(2 0)\rlvec(0 -0.3)
\move(2 0)\rlvec(0 1)\rlvec(-1 0)
\move(2 1)\rlvec(-1 -1)
\move(0 0)\rlvec(0.25 0.25)
\move(1 1)\rlvec(-0.25 -0.25)
\htext(0.5 0.5){$j$}
\htext(1.25 0.75){$i$}
\end{texdraw}
}%
\savebox{\tmpfigl}{\begin{texdraw}
\fontsize{7}{7}\selectfont
\textref h:C v:C
\drawdim em
\setunitscale 1.9
\move(0 -0.3)\lvec(0 1)\lvec(1 1)\lvec(1 -0.3)
\move(0 0)\rlvec(2 0)\rlvec(0 -0.3)
\move(2 0)\lvec(2 0.5)
\move(0 0)\rlvec(0.25 0.25)
\move(1 1)\rlvec(-0.25 -0.25)
\htext(0.5 0.5){$j$}
\move(1 0)\rlvec(0.25 0.125)
\move(2 0.5)\rlvec(-0.25 -0.125)
\htext(1.5 0.25){$j$}
\lpatt(0.05 0.15)
\move(1 0.5)\rlvec(1 0)
\end{texdraw}
}%
\savebox{\tmpfign}{\begin{texdraw}
\fontsize{7}{7}\selectfont
\textref h:C v:C
\drawdim em
\setunitscale 1.9
\move(0 -0.3)\lvec(0 1)\lvec(1 1)\lvec(1 -0.3)
\move(0 0)\rlvec(2 0)\rlvec(0 -0.3)
\move(2 -0.2)\lvec(2 1)\rlvec(-1 0)
\move(0 0)\rlvec(0.25 0.25)
\move(1 1)\rlvec(-0.25 -0.25)
\htext(0.5 0.5){$j$}
\move(1 0)\rlvec(0.25 0.25)
\move(2 1)\rlvec(-0.25 -0.25)
\htext(1.5 0.5){$j$}
\end{texdraw}%
}%
\begin{center}
\raisebox{-2.2em}{\usebox{\tmpfige}}
\qquad
\raisebox{-2.2em}{\usebox{\tmpfigf}}
\qquad
\raisebox{-2.2em}{\usebox{\tmpfigh}}
\qquad
\raisebox{-2.2em}{\usebox{\tmpfigi}}\\[3mm]
\raisebox{-2.2em}{\usebox{\tmpfiga}}
\qquad
\raisebox{-2.2em}{\usebox{\tmpfigb}}
\qquad
\raisebox{-2.2em}{\usebox{\tmpfigc}}
\qquad
\raisebox{-2.2em}{\usebox{\tmpfigd}}\\[3mm]
\raisebox{-2.2em}{\usebox{\tmpfigg}}
\qquad
\raisebox{-2.2em}{\usebox{\tmpfigm}}
\qquad
\raisebox{-2.2em}{\usebox{\tmpfigj}}
\qquad
\raisebox{-2.2em}{\usebox{\tmpfigk}}
\qquad
\raisebox{-2.2em}{\usebox{\tmpfigl}}\\[3mm]
\raisebox{-2.2em}{\usebox{\tmpfign}}
\end{center}
Here, $j=\zo, \nn$. As mentioned in Remark~\ref{expl},
the right (or left) of the first (resp. last) two diagrams
in the first row is the upper (resp. lower) broken half
of a whole $j$-block.

\noindent
In each of these cases, it is possible to obtain one of the following
two conclusions.
\begin{enumerate}
\item The columns $C$ and $C'$ of the proper Young wall $\mathbf{Y}$
      contains a layer in which the tops are of unit depth and of
      the same integer height.
\item $\vphi_i(C) \leq \veps_i(C')$.
\end{enumerate}
The first of these conclusion is a violation of the properness
of $\mathbf{Y}$.
As for the second,
since $C$ is column corresponding the leftmost $0$
in the $i$-signature of $\mathbf{Y}$
we must have $\vphi_i(C) > \veps_i(C')$.
Both of these conclusions brings us to contradictions,
and hence the resulting $\fit \mathbf{Y}$ must have been a proper Young
wall.
\end{proof}

We define the maps $\cwt: \pwspace \rightarrow \bar P$,
$\vphi_i, \veps_i : \pwspace \rightarrow \Z$ by setting
\begin{equation}
\begin{aligned}
\cwt(\mathbf{Y}) & = \sum_{i\in I}
(\vphi_i(\mathbf{Y})-\veps_i(\mathbf{Y})) \La_i, \\
\vphi_i(\mathbf{Y}) & = \text{the number of 0's in the
$i$-signature of $\mathbf{Y}$},\\
\veps_i(\mathbf{Y}) & = \text{the number of 1's in the
$i$-signature of $\mathbf{Y}$}.
\end{aligned}
\end{equation}
Then it is straightforward to verify that the following theorem
holds.

\begin{thm}
The set $\pwspace$ of all level-$l$ proper Young walls, together
with the maps $\cwt$, $\fit$, $\eit$, $\vphi_i$, $\veps_i$ $(i\in
I)$, forms a $\uq'(\g)$-crystal.
\end{thm}

%%%%%%%%%%%%%%%%%%%%%%%%%%%%%%%%%%%%%%%%%%%%%%%%%%%%%%%%%%%%%%%%%%%%%%%%%%
\section{Young wall realization of $\hwc(\la)$}%
\label{sec:6}

In this section, we give a new realization of arbitrary level
irreducible highest weight crystals in terms of reduced proper
Young walls.

Let $\la = a_0 \La_0 + a_1 \La_1 + \cdots + a_n \La_n$ be
a dominant integral weight of level-$l$
so that
\begin{equation}
l=a_0 +a_1 +2(a_2 + \cdots + a_{n-2}) + a_{n-1} +a_n.
\end{equation}
We would like to define the \defi{ground-state wall $\mathbf{Y}_{\la}$
of weight $\la$}.
Depending on whether $a_0 \geq a_1$ or $a_0 \leq a_1$ and
$a_{n-1} \geq a_n$ or $a_{n-1} \geq a_n$, there are four
different forms of the ground-state wall.

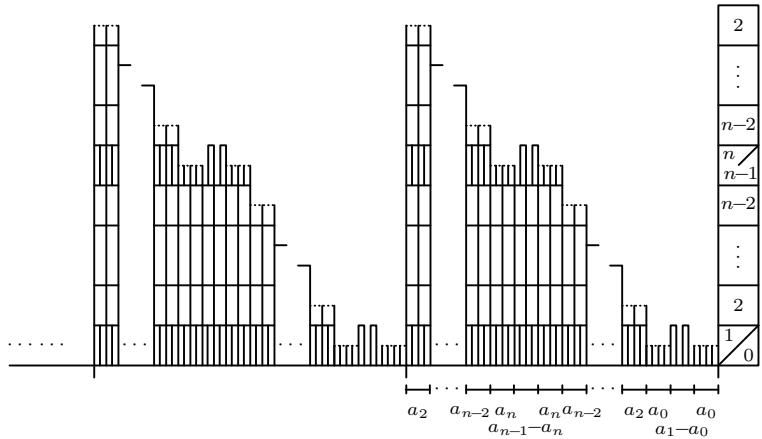
\vspace{3mm}
\begin{center}
\begin{texdraw}
\fontsize{7}{7}\selectfont
\textref h:C v:C
\drawdim em
\setunitscale 1.9
\move(0 0)
\bsegment
\move(0 0)\rlvec(1 0)\rlvec(0 9)\rlvec(-1 0)\rlvec(0 -9)
\move(0 1)\rlvec(1 0)
\move(0 2)\rlvec(1 0)
\move(0 3.5)\rlvec(1 0)
\move(0 4.5)\rlvec(1 0)
\move(0 5.5)\rlvec(1 0)
\move(0 6.5)\rlvec(1 0)
\move(0 8)\rlvec(1 0)
\move(0 0)\rlvec(1 1)
\move(1 5.5)\rlvec(-0.5 -0.5)
\htext(0.75 0.25){$0$}
\htext(0.25 0.75){$1$}
\htext(0.5 1.5){$2$}
\vtext(0.5 2.75){$\cdots$}
\htext(0.5 4){$n\!\!-\!\!2$}
\htext(0.55 4.75){$n\!\!-\!\!1$}
\htext(0.25 5.25){$n$}
\htext(0.5 6){$n\!\!-\!\!2$}
\vtext(0.5 7.25){$\cdots$}
\htext(0.5 8.5){$2$}
\move(0 0)\rlvec(0 -0.3)
\esegment

\move(0 0)
\bsegment
\move(-0.6 0)
  \bsegment
  \move(0.3 0)\rlvec(0 0.5)
  \move(0.15 0)\rlvec(0 0.5)
  \move(0.45 0)\rlvec(0 0.5)
  \move(0 0)\rlvec(0 0.5)
  \lpatt(0.02 0.1)
  \move(0 0.5)\rlvec(0.6 0)
  \esegment
\move(-1.2 0)
  \bsegment
  \move(0.45 0)\rlvec(0 1)\rlvec(-0.15 0)\rlvec(0 -1)
  \move(0.15 0)\rlvec(0 1)\rlvec(-0.15 0)\rlvec(0 -1)
  \esegment
\move(-1.8 0)
  \bsegment
  \move(0.3 0)\rlvec(0 0.5)
  \move(0.15 0)\rlvec(0 0.5)
  \move(0.45 0)\rlvec(0 0.5)
  \move(0 0)\rlvec(0 0.5)
  \lpatt(0.02 0.1)
  \move(0 0.5)\rlvec(0.6 0)
  \esegment
\move(-2.4 0)
  \bsegment
  \move(0.6 0)\rlvec(0 1.5)
  \move(0.3 0)\rlvec(0 1.5)
  \move(0.15 0)\rlvec(0 1)
  \move(0.45 0)\rlvec(0 1)
  \move(0.6 1)\rlvec(-0.6 0)
  \lpatt(0.02 0.1)
  \move(0 1.5)\rlvec(0.6 0)
\esegment
\move(-2.4 0)
\bsegment
\move(0 0)\rlvec(0 2.5)\rlvec(-0.3 0)
\esegment
\move(-3.9 0)
  \bsegment
  \move(0.6 3)\rlvec(0.3 0)
  \move(0.6 0)\rlvec(0 4)
  \move(0.3 0)\rlvec(0 4)
  \move(0 1)\rlvec(0.6 0)
  \move(0 2)\rlvec(0.6 0)
  \move(0 3.5)\rlvec(0.6 0)
  \move(0.15 0)\rlvec(0 1)
  \move(0.45 0)\rlvec(0 1)
  \lpatt(0.02 0.1)
  \move(0 4)\rlvec(0.6 0)
  \esegment
\move(-4.5 0)
  \bsegment
  \move(0.6 0)\rlvec(0 5)
  \move(0.3 0)\rlvec(0 5)
  \move(0 4.5)\rlvec(0.6 0)
  \move(0 1)\rlvec(0.6 0)
  \move(0 2)\rlvec(0.6 0)
  \move(0 3.5)\rlvec(0.6 0)
  \move(0.15 0)\rlvec(0 1)
  \move(0.45 0)\rlvec(0 1)
  \move(0.15 4.5)\rlvec(0 0.5)
  \move(0.45 4.5)\rlvec(0 0.5)
  \lpatt(0.02 0.1)
  \move(0 5)\rlvec(0.6 0)
  \esegment
\move(-5.1 0)
\bsegment
\move(0.6 0)\rlvec(0 4.5)\rlvec(-0.6 0)
\move(0.3 0)\rlvec(0 4.5)
\move(0 1)\rlvec(0.6 0)
\move(0 2)\rlvec(0.6 0)
\move(0 3.5)\rlvec(0.6 0)
\move(0.15 0)\rlvec(0 1)
\move(0.45 0)\rlvec(0 1)
\move(0.6 4.5)\rlvec(0 1)\rlvec(-0.15 0)\rlvec(0 -1)
\move(0.3 4.5)\rlvec(0 1)\rlvec(-0.15 0)\rlvec(0 -1)
\esegment
\move(-5.7 0)
  \bsegment
  \move(0.6 0)\rlvec(0 5)
  \move(0.3 0)\rlvec(0 5)
  \move(0 4.5)\rlvec(0.6 0)
  \move(0 1)\rlvec(0.6 0)
  \move(0 2)\rlvec(0.6 0)
  \move(0 3.5)\rlvec(0.6 0)
  \move(0.15 0)\rlvec(0 1)
  \move(0.45 0)\rlvec(0 1)
  \move(0.15 4.5)\rlvec(0 0.5)
  \move(0.45 4.5)\rlvec(0 0.5)
  \lpatt(0.02 0.1)
  \move(0 5)\rlvec(0.6 0)
\esegment
\move(-6.3 0)
  \bsegment
  \move(0.6 0)\rlvec(0 6)
  \move(0.3 0)\rlvec(0 6)
  \move(0 4.5)\rlvec(0.6 0)
  \move(0 1)\rlvec(0.6 0)
  \move(0 2)\rlvec(0.6 0)
  \move(0 3.5)\rlvec(0.6 0)
  \move(0.15 0)\rlvec(0 1)
  \move(0.45 0)\rlvec(0 1)
  \move(0 5.5)\rlvec(0.6 0)
  \move(0 0)\rlvec(0 7)\rlvec(-0.3 0)
  \move(0.15 4.5)\rlvec(0 1)
  \move(0.45 4.5)\rlvec(0 1)
  \lpatt(0.02 0.1)
  \move(0 6)\rlvec(0.6 0)
\esegment
\move(-7.8 0)
  \bsegment
  \move(0.6 0)\rlvec(0 8.5)
  \move(0.3 0)\rlvec(0 8.5)
  \move(0 4.5)\rlvec(0.6 0)
  \move(0 1)\rlvec(0.6 0)
  \move(0 6.5)\rlvec(0.6 0)
  \move(0 2)\rlvec(0.6 0)
  \move(0 8)\rlvec(0.6 0)
  \move(0.6 7.5)\rlvec(0.3 0)
  \move(0 3.5)\rlvec(0.6 0)
  \move(0.15 0)\rlvec(0 1)
  \move(0.45 0)\rlvec(0 1)
  \move(0 5.5)\rlvec(0.6 0)
  \move(0.15 4.5)\rlvec(0 1)
  \move(0.45 4.5)\rlvec(0 1)
  \lpatt(0.02 0.1)
  \move(0 8.5)\rlvec(0.6 0)
\esegment
\move(0 0)\rlvec(-7.8 0)\rlvec(0 8.5)
\move(-7.8 0)\rlvec(0 -0.3)
\htext(-6.75 0.5){$\cdots$}
\htext(-2.85 0.5){$\cdots$}

\move(0 -0.5)\rlvec(0 -0.2)
\move(0 -0.6)\rlvec(-0.6 0)    \htext(-0.3 -1.2){$a_0$}
\move(-0.6 -0.5)\rlvec(0 -0.2)
\move(-0.6 -0.6)\rlvec(-0.6 0) \htext(-0.9 -1.6){$a_1\!\!-\!\!a_0$}
\move(-1.2 -0.5)\rlvec(0 -0.2)
\move(-1.2 -0.6)\rlvec(-0.6 0) \htext(-1.5 -1.2){$a_0$}
\move(-1.8 -0.5)\rlvec(0 -0.2)
\move(-1.8 -0.6)\rlvec(-0.6 0) \htext(-2.1 -1.2){$a_2$}
\move(-2.4 -0.5)\rlvec(0 -0.2)
\htext(-2.85 -0.6){$\cdots$}
\move(-3.3 -0.5)\rlvec(0 -0.2)
\move(-3.3 -0.6)\rlvec(-0.6 0) \htext(-3.4 -1.2){$a_{n\!-\!2}$}
\move(-3.9 -0.5)\rlvec(0 -0.2)
\move(-3.9 -0.6)\rlvec(-0.6 0) \htext(-4.2 -1.2){$a_n$}
\move(-4.5 -0.5)\rlvec(0 -0.2)
\move(-4.5 -0.6)\rlvec(-0.6 0) \htext(-4.8 -1.6){$a_{n\!-\!1}\!\!-\!\!a_n$}
\move(-5.1 -0.5)\rlvec(0 -0.2)
\move(-5.1 -0.6)\rlvec(-0.6 0) \htext(-5.3 -1.2){$a_n$}
\move(-5.7 -0.5)\rlvec(0 -0.2)
\move(-5.7 -0.6)\rlvec(-0.6 0) \htext(-6.2 -1.2){$a_{n\!-\!2}$}
\move(-6.3 -0.5)\rlvec(0 -0.2)
\htext(-6.75 -0.6){$\cdots$}
\move(-7.2 -0.5)\rlvec(0 -0.2)
\move(-7.2 -0.6)\rlvec(-0.6 0) \htext(-7.5 -1.2){$a_2$}
\move(-7.8 -0.5)\rlvec(0 -0.2)
\esegment

\move(-7.8 0)
\bsegment
\move(-0.6 0)
  \bsegment
  \move(0.3 0)\rlvec(0 0.5)
  \move(0.15 0)\rlvec(0 0.5)
  \move(0.45 0)\rlvec(0 0.5)
  \move(0 0)\rlvec(0 0.5)
  \lpatt(0.02 0.1)
  \move(0 0.5)\rlvec(0.6 0)
  \esegment
\move(-1.2 0)
  \bsegment
  \move(0.45 0)\rlvec(0 1)\rlvec(-0.15 0)\rlvec(0 -1)
  \move(0.15 0)\rlvec(0 1)\rlvec(-0.15 0)\rlvec(0 -1)
  \esegment
\move(-1.8 0)
  \bsegment
  \move(0.3 0)\rlvec(0 0.5)
  \move(0.15 0)\rlvec(0 0.5)
  \move(0.45 0)\rlvec(0 0.5)
  \move(0 0)\rlvec(0 0.5)
  \lpatt(0.02 0.1)
  \move(0 0.5)\rlvec(0.6 0)
  \esegment
\move(-2.4 0)
  \bsegment
  \move(0.6 0)\rlvec(0 1.5)
  \move(0.3 0)\rlvec(0 1.5)
  \move(0.15 0)\rlvec(0 1)
  \move(0.45 0)\rlvec(0 1)
  \move(0.6 1)\rlvec(-0.6 0)
  \lpatt(0.02 0.1)
  \move(0 1.5)\rlvec(0.6 0)
  \esegment
\move(-2.4 0)
  \bsegment
  \move(0 0)\rlvec(0 2.5)\rlvec(-0.3 0)
  \esegment
\move(-3.9 0)
  \bsegment
  \move(0.6 3)\rlvec(0.3 0)
  \move(0.6 0)\rlvec(0 4)
  \move(0.3 0)\rlvec(0 4)
  \move(0 1)\rlvec(0.6 0)
  \move(0 2)\rlvec(0.6 0)
  \move(0 3.5)\rlvec(0.6 0)
  \move(0.15 0)\rlvec(0 1)
  \move(0.45 0)\rlvec(0 1)
  \lpatt(0.02 0.1)
  \move(0 4)\rlvec(0.6 0)
  \esegment
\move(-4.5 0)
  \bsegment
  \move(0.6 0)\rlvec(0 5)
  \move(0.3 0)\rlvec(0 5)
  \move(0 4.5)\rlvec(0.6 0)
  \move(0 1)\rlvec(0.6 0)
  \move(0 2)\rlvec(0.6 0)
  \move(0 3.5)\rlvec(0.6 0)
  \move(0.15 0)\rlvec(0 1)
  \move(0.45 0)\rlvec(0 1)
  \move(0.15 4.5)\rlvec(0 0.5)
  \move(0.45 4.5)\rlvec(0 0.5)
  \lpatt(0.02 0.1)
  \move(0 5)\rlvec(0.6 0)
  \esegment
\move(-5.1 0)
  \bsegment
  \move(0.6 0)\rlvec(0 4.5)\rlvec(-0.6 0)
  \move(0.3 0)\rlvec(0 4.5)
  \move(0 1)\rlvec(0.6 0)
  \move(0 2)\rlvec(0.6 0)
  \move(0 3.5)\rlvec(0.6 0)
  \move(0.15 0)\rlvec(0 1)
  \move(0.45 0)\rlvec(0 1)
  \move(0.6 4.5)\rlvec(0 1)\rlvec(-0.15 0)\rlvec(0 -1)
  \move(0.3 4.5)\rlvec(0 1)\rlvec(-0.15 0)\rlvec(0 -1)
  \esegment
\move(-5.7 0)
  \bsegment
  \move(0.6 0)\rlvec(0 5)
  \move(0.3 0)\rlvec(0 5)
  \move(0 4.5)\rlvec(0.6 0)
  \move(0 1)\rlvec(0.6 0)
  \move(0 2)\rlvec(0.6 0)
  \move(0 3.5)\rlvec(0.6 0)
  \move(0.15 0)\rlvec(0 1)
  \move(0.45 0)\rlvec(0 1)
  \move(0.15 4.5)\rlvec(0 0.5)
  \move(0.45 4.5)\rlvec(0 0.5)
  \lpatt(0.02 0.1)
  \move(0 5)\rlvec(0.6 0)
  \esegment
\move(-6.3 0)
  \bsegment
  \move(0.6 0)\rlvec(0 6)
  \move(0.3 0)\rlvec(0 6)
  \move(0 4.5)\rlvec(0.6 0)
  \move(0 1)\rlvec(0.6 0)
  \move(0 2)\rlvec(0.6 0)
  \move(0 3.5)\rlvec(0.6 0)
  \move(0.15 0)\rlvec(0 1)
  \move(0.45 0)\rlvec(0 1)
  \move(0 5.5)\rlvec(0.6 0)
  \move(0 0)\rlvec(0 7)\rlvec(-0.3 0)
  \move(0.15 4.5)\rlvec(0 1)
  \move(0.45 4.5)\rlvec(0 1)
  \lpatt(0.02 0.1)
  \move(0 6)\rlvec(0.6 0)
  \esegment
\move(-7.8 0)
  \bsegment
  \move(0.6 0)\rlvec(0 8.5)
  \move(0.3 0)\rlvec(0 8.5)
  \move(0 4.5)\rlvec(0.6 0)
  \move(0 1)\rlvec(0.6 0)
  \move(0 6.5)\rlvec(0.6 0)
  \move(0 2)\rlvec(0.6 0)
  \move(0 8)\rlvec(0.6 0)
  \move(0.6 7.5)\rlvec(0.3 0)
  \move(0 3.5)\rlvec(0.6 0)
  \move(0.15 0)\rlvec(0 1)
  \move(0.45 0)\rlvec(0 1)
  \move(0 5.5)\rlvec(0.6 0)
  \move(0.15 4.5)\rlvec(0 1)
  \move(0.45 4.5)\rlvec(0 1)
  \lpatt(0.02 0.1)
  \move(0 8.5)\rlvec(0.6 0)
  \esegment
\move(0 0)\rlvec(-7.8 0)\rlvec(0 8.5)
\move(-7.8 0)\rlvec(0 -0.3)
\htext(-6.75 0.5){$\cdots$}
\htext(-2.85 0.5){$\cdots$}
\esegment

\move(-15.6 0)\rlvec(-2.1 0)
\htext(-16.6 0.5){$\cdots$}
\htext(-17.4 0.5){$\cdots$}
\end{texdraw}%
\end{center}

\vspace{4mm}

\noindent
We have drawn just one of these here.
It corresponds to the case when
$a_0 \leq a_1$ and $a_{n-1}\geq a_n$.
We have drawn the left-side-view of each of the slices
that make up the Young wall.
At the right end, we have drawn the pattern for stacking
the blocks, so as to show the color of the blocks placed at
each height.
We believe that the readers can easily fill in
the ground-state walls for the other cases.
Note that these are reduce proper Young walls of level-$l$.
We caution the readers that the pattern on the right
is just for the even $i$th columns.
As given by figure~(\ref{ywallpattern}),
the odd columns will be stacked in a pattern with $0,1$
exchanged and $n-1,n$ exchanged.

A level-$l$ proper Young wall obtained by adding finitely many
blocks to the ground-state wall $\gswall_{\la}$ is said to have
been \defi{built on $\gswall_{\la}$}.
We denote by $\pwspace(\la)$ (and $\rpwspace(\la)$)
the set of all proper (resp. reduced proper) Young walls
built on $\gswall_{\la}$.

It is clear that the set $\pwspace(\la)$
forms a $\uqp(\g)$-subcrystal of $\pwspace$.
For each $\mathbf{Y} \in \pwspace(\la)$, we define its
\defi{affine weight} to be
\begin{equation}
\wt(\mathbf{Y}) = \la - \sum_{i=0}^{n} k_i \ali,
\end{equation}
where $k_i$ is the number of $i$-blocks that have been added to
$\gswall_{\la}$. Then we obtain\,:
\begin{prop}
The set $\pwspace(\la)$ of all level-$l$ proper Young walls built
on $\gswall_{\la}$, together with the maps $\eit$, $\fit$,
$\veps_i$, $\vphi_i$ $(i\in I)$ given in Section~\ref{sec:5}
and $\wt$, forms a $\uq(\g)$-crystal.
\end{prop}

Recall that the set $\pathspace(\la)$ of all $\la$-paths gives a
realization of the irreducible highest weight crystal $\hwc(\la)$
(see Section~\ref{sec:2}).
We will show that $\rpwspace(\la)$ is a new realization of
$\hwc(\la)$ by giving a $\uq(\g)$-crystal isomorphism
$\Phi: \rpwspace(\la)
\overset{\sim}\longrightarrow \pathspace(\la)$.

Let $\mathbf{Y} = (\mathbf{Y}(k))_{k=0}^{\infty}$ be a reduced
proper Young wall built on $\mathbf{Y}_{\la}$.
Here, $\mathbf{Y}(k)$ denotes the $k$th column of $\mathbf{Y}$.
Using the $\uq'(\g)$-crystal isomorphism $\psi: \oldcrystal \overset{\sim}
\longrightarrow \newcrystal$ given in Theorem~\ref{thm:37}, we get a map
$\Phi: \rpwspace(\la) \longrightarrow \pathspace(\la)$ which is
defined by
\begin{equation} \label{eq:main}
\Phi(\mathbf{Y}) = (\psi^{-1}(\mathbf{Y}(k)))_{k=0}^{\infty}.
\end{equation}
Note that the ground-state wall $\gswall_{\la}$ is mapped onto the
ground-state path $\path_{\la}$.

Conversely, to each $\la$-path $\path = (\path(k))_{k=0}^{\infty}
\in \pathspace(\la)$, we can associate a proper Young wall
$\mathbf{Y}=(\mathbf{Y}(k))_{k=0}^{\infty}$ built on
$\mathbf{Y}_{\la}$ such that $\psi(\path(k)) = \mathbf{Y}(k)$ for
all $k \ge 0$. By removing an appropriate number of
$\delta$'s (proceeding from left to right), one can easily see that there
exists a unique reduced proper Young wall
$\mathbf{Y}=(\mathbf{Y}(k))_{k=0}^{\infty}$ with this property,
which shows that $\Phi: \rpwspace(\la) \longrightarrow
\pathspace(\la)$ is a bijection.

Here is our main realization theorem,
which shows that the above bijection is a crystal isomorphism.

\begin{thm}\label{thm:main thm}
The bijection $\Phi: \rpwspace(\la) \longrightarrow
\pathspace(\la)$ defined by \eqref{eq:main} is a $\uq(\g)$-crystal
isomorphism.
Therefore, we have a $\uq(\g)$-crystal isomorphism
\begin{equation}
\rpwspace(\la) \overset{\sim} \longrightarrow \pathspace(\la)
\overset {\sim} \longrightarrow \hwc(\la).
\end{equation}
\end{thm}
\begin{proof}
We shall focus our efforts on showing that the set $\rpwspace(\la)$
is a $\uq(\mathfrak{g})$-subcrystal of $\pwspace(\la)$ and that
the map $\Phi$ commutes with the Kashiwara operators $\fit$ and $\eit$.
Other parts of the proof are similar or easy.

First, let us assume that the set $\rpwspace(\la)$ is a
$\uq(\mathfrak{g})$-subcrystal of $\pwspace(\la)$.
Then the action of Kashiwara operators on $\pwspace$
defined in Section~\ref{sec:5} is well-defined on $\rpwspace(\la)$.
The fact that the map $\Phi$ is
determined by the $\uq'(\g)$-crystal isomorphism $\psi:
\oldcrystal \overset{\sim} \longrightarrow \newcrystal$,
and that the action of Kashiwara operators on $\rpwspace(\la)$
follows the process for defining it on
$\pathspace(\la)$,
implies naturally that the map $\Phi$ commutes with the Kashiwara
operators $\fit$ and $\eit$.

Now, it remains to show that the set $\rpwspace(\la)$
is a $\uq(\mathfrak{g})$-subcrystal of $\pwspace(\la)$.
It suffices to prove that
the action of Kashiwara operators on $\rpwspace(\la)$
satisfies the following properties :
\begin{equation}
\fit \rpwspace(\la) \subset \rpwspace(\la) \cup \{0\}, \qquad \eit
\rpwspace(\la) \subset \rpwspace(\la) \cup \{0\} \qquad \text{for
all} \ \ i\in I.
\end{equation}

Suppose that there exists some $\mathbf{Y}\in\rpwspace(\la)$ for which
$\fit \mathbf{Y} \not\in \rpwspace(\la) \cup \{0\}$.
We assume that $\fit$ has acted on the $j$th column of $\mathbf{Y}$
and set $\path = \spaceosi(\mathbf{Y})$.
Then, the action of $\fit$ on $\path$ would also have been
on the $j$th tensor component of $\path$.

Since $\fit \mathbf{Y} \in \pwspace(\la)$,
we may remove finitely many $\delta$'s from $\fit \mathbf{Y}$
to obtain a reduce proper Young wall $\mathbf{Y}'$.
The number of $\delta$'s removed is nonzero since $\fit \mathbf{Y}
\not\in \rpwspace(\la)$.
Since every column of $\fit \mathbf{Y}$ is related to
the corresponding column of $\mathbf{Y}'$ under the previously defined
equivalence relation~(\ref{eq:rel}), we have
\begin{equation*}
\fit(\path) = \spaceosi(\mathbf{Y}').
\end{equation*}

Let us apply $\eit$ to both $\mathbf{Y}'$ and $\fit(\path)$.
We have $\path = \eit(\fit(\path))$ and the action of
$\eit$ on $\fit(\path)$ would have been on the
$j$th tensor component.
The action of $\eit$ of $\mathbf{Y}'$ will also be on
the $j$th column.
Hence
the proper Young wall $\eit \mathbf{Y}'$ may be obtained from the
reduced proper Young wall $\mathbf{Y}$ by removing finitely many
$\delta$'s.
We may now remove finitely many $\delta$'s from $\eit
\mathbf{Y}'$ to obtain a reduced proper Young wall $\mathbf{Y}''$
which also
corresponds to $\path$ under the map
$\spaceosi$.

Recall that we started out with a reduced proper Young wall
$\mathbf{Y}$,
added an $i$-block to the $j$th column of $\mathbf{Y}$ to obtain
$\fit \mathbf{Y}$,
removed finitely many $\delta$'s from $\fit \mathbf{Y}$
to obtain a reduced proper Young wall $\mathbf{Y}'$, and removed an
$i$-block from the $j$th column of $\mathbf{Y}'$ to obtain
$\eit \mathbf{Y}'$, and finally, removed finitely
many $\delta$'s from $\eit \mathbf{Y}'$ to obtain a reduced
proper Young wall $\mathbf{Y}''$ which corresponds to $\path$ under
the map $\spaceosi$.
Therefore, we have
$\Phi(\mathbf{Y})=\path=\Phi(\mathbf{Y''})$, but $\mathbf{Y} \neq
\mathbf{Y''}$, which is a contradiction. Hence $\fit \mathbf{Y}$
must be reduced.

Similarly, one can show that $\eit \rpwspace(\la) \subset
\rpwspace(\la) \cup \{0\}$ for all $i\in I$, which completes the
proof of our claim.
\end{proof}

We close this paper with an example of Young wall realization of
irreducible highest weight crystals.

\begin{example}
The top part of the crystal graph
$\rpwspace(\La_0 +\La_4)$ for type $D_4^{(1)}$ is given below.
In each Young wall, we have shaded the part that has received
some change through the action of $\fit$.\\
\begin{center}
\begin{texdraw}
\drawdim in
\arrowheadsize l:0.065 w:0.03
\arrowheadtype t:F
\fontsize{5}{5}\selectfont
\textref h:C v:C
\drawdim em
\setunitscale 1.7
\move(0.4 0)
\bsegment
\move(0 0)%%%%%%%%%%%%
\bsegment
\move(0 0)
 \bsegment
 \move(0 0)\rlvec(1 0)
 \move(0 1)\rlvec(1 0)
 \move(0 2)\rlvec(1 0)
 \move(0 3)\rlvec(1 0)
 \move(1 0)\rlvec(0 3)
 \move(0 -0.2)\rlvec(0 3.2)
 \move(0 0)\rlvec(1 1)
 \move(0 2)\rlvec(1 1)
 \htext(0.25 0.75){$1$}
 \htext(0.75 0.25){$0$}
 \htext(0.75 2.25){$3$}
 \htext(0.25 2.75){$4$}
 \htext(0.5 1.5){$2$}
\esegment
\move(-0.6 0)
\bsegment
 \move(0 -0.2)\rlvec(0 3.2)
 \move(0.3 0)\rlvec(0 2)
 \move(0.15 0)\rlvec(0 1)
 \move(0.45 0)\rlvec(0 1)
 \move(0 0)\rlvec(0.6 0)
 \move(0 1)\rlvec(0.3 0)
 \move(0.45 1)\rlvec(0.15 0)
 \move(0 2)\rlvec(0.3 0)
 \move(0 3)\rlvec(0.15 0)\rlvec(0 -1)
\esegment
\move(-1.2 0)
\bsegment
 \move(0 -0.2)\rlvec(0 3.2)
 \move(0.3 0)\rlvec(0 2)
 \move(0.15 0)\rlvec(0 1)
 \move(0.45 0)\rlvec(0 1)
 \move(0 0)\rlvec(0.6 0)
 \move(0 1)\rlvec(0.3 0)
 \move(0.45 1)\rlvec(0.15 0)
 \move(0 2)\rlvec(0.3 0)
 \move(0 3)\rlvec(0.15 0)\rlvec(0 -1)
\esegment
\move(-1.8 0)
\bsegment
 \move(0 -0.2)\rlvec(0 3.2)
 \move(0.3 0)\rlvec(0 2)
 \move(0.15 0)\rlvec(0 1)
 \move(0.45 0)\rlvec(0 1)
 \move(-0.4 0)\rlvec(1 0)
 \move(0 1)\rlvec(0.3 0)
 \move(0.45 1)\rlvec(0.15 0)
 \move(0 2)\rlvec(0.3 0)
 \move(0 3)\rlvec(0.15 0)\rlvec(0 -1)
\esegment
\esegment
\move(-4 -6)%%%%%%%%%%%%%%%%
\bsegment
\move(0 0)
 \bsegment
 \move(0 0)\rlvec(1 0)
 \move(0 1)\rlvec(1 0)
 \move(0 2)\rlvec(1 0)
 \move(0 3)\rlvec(1 0)
 \move(1 0)\rlvec(0 3)
 \move(0 -0.2)\rlvec(0 3.2)
 \move(0 0)\rlvec(1 1)
 \move(0 2)\rlvec(1 1)
 \htext(0.25 0.75){$1$}
 \htext(0.75 0.25){$0$}
 \htext(0.75 2.25){$3$}
 \htext(0.25 2.75){$4$}
 \htext(0.5 1.5){$2$}
\esegment
\move(-0.6 0)
\bsegment
\move(0.3 0)\rlvec(0.15 0)\rlvec(0 1)\rlvec(-0.15 0)\ifill f:0.5
 \move(0 -0.2)\rlvec(0 3.2)
 \move(0.3 0)\rlvec(0 2)
 \move(0.15 0)\rlvec(0 1)
 \move(0.45 0)\rlvec(0 1)
 \move(0 0)\rlvec(0.6 0)
 \move(0 1)\rlvec(0.3 0)
 \move(0.3 1)\rlvec(0.3 0)
 \move(0 2)\rlvec(0.3 0)
 \move(0 3)\rlvec(0.15 0)\rlvec(0 -1)
\esegment
\move(-1.2 0)
\bsegment
 \move(0 -0.2)\rlvec(0 3.2)
 \move(0.3 0)\rlvec(0 2)
 \move(0.15 0)\rlvec(0 1)
 \move(0.45 0)\rlvec(0 1)
 \move(0 0)\rlvec(0.6 0)
 \move(0 1)\rlvec(0.3 0)
 \move(0.45 1)\rlvec(0.15 0)
 \move(0 2)\rlvec(0.3 0)
 \move(0 3)\rlvec(0.15 0)\rlvec(0 -1)
\esegment
\move(-1.8 0)
\bsegment
 \move(0 -0.2)\rlvec(0 3.2)
 \move(0.3 0)\rlvec(0 2)
 \move(0.15 0)\rlvec(0 1)
 \move(0.45 0)\rlvec(0 1)
 \move(-0.4 0)\rlvec(1 0)
 \move(0 1)\rlvec(0.3 0)
 \move(0.45 1)\rlvec(0.15 0)
 \move(0 2)\rlvec(0.3 0)
 \move(0 3)\rlvec(0.15 0)\rlvec(0 -1)
\esegment
\esegment
\move(4 -6)%%%%%%%%%%%%%%%
\bsegment
\move(0 0)
 \bsegment
 \move(0 0)\rlvec(1 0)
 \move(0 1)\rlvec(1 0)
 \move(0 2)\rlvec(1 0)
 \move(0 3)\rlvec(1 0)
 \move(1 0)\rlvec(0 3)
 \move(0 -0.2)\rlvec(0 3.2)
 \move(0 0)\rlvec(1 1)
 \move(0 2)\rlvec(1 1)
 \htext(0.25 0.75){$1$}
 \htext(0.75 0.25){$0$}
 \htext(0.75 2.25){$3$}
 \htext(0.25 2.75){$4$}
 \htext(0.5 1.5){$2$}
\esegment
\move(-0.6 0)
\bsegment
\move(0.15 2)\rlvec(0.15 0)\rlvec(0 1)\rlvec(-0.15 0)\ifill f:0.5
 \move(0 -0.2)\rlvec(0 3.2)
 \move(0.3 0)\rlvec(0 2)
 \move(0.15 0)\rlvec(0 1)
 \move(0.45 0)\rlvec(0 1)
 \move(0 0)\rlvec(0.6 0)
 \move(0 1)\rlvec(0.3 0)
 \move(0.45 1)\rlvec(0.15 0)
 \move(0 2)\rlvec(0.3 0)
 \move(0 3)\rlvec(0.15 0)\rlvec(0 -1)
 \move(0.15 3)\rlvec(0.15 0)\rlvec(0 -1)
\esegment
\move(-1.2 0)
\bsegment
 \move(0 -0.2)\rlvec(0 3.2)
 \move(0.3 0)\rlvec(0 2)
 \move(0.15 0)\rlvec(0 1)
 \move(0.45 0)\rlvec(0 1)
 \move(0 0)\rlvec(0.6 0)
 \move(0 1)\rlvec(0.3 0)
 \move(0.45 1)\rlvec(0.15 0)
 \move(0 2)\rlvec(0.3 0)
 \move(0 3)\rlvec(0.15 0)\rlvec(0 -1)
\esegment
\move(-1.8 0)
\bsegment
 \move(0 -0.2)\rlvec(0 3.2)
 \move(0.3 0)\rlvec(0 2)
 \move(0.15 0)\rlvec(0 1)
 \move(0.45 0)\rlvec(0 1)
 \move(-0.4 0)\rlvec(1 0)
 \move(0 1)\rlvec(0.3 0)
 \move(0.45 1)\rlvec(0.15 0)
 \move(0 2)\rlvec(0.3 0)
 \move(0 3)\rlvec(0.15 0)\rlvec(0 -1)
\esegment
\esegment
\move(0 -12)%%%%%%%%%%%
\bsegment
\move(0 0)
 \bsegment
 \move(0 0)\rlvec(1 0)
 \move(0 1)\rlvec(1 0)
 \move(0 2)\rlvec(1 0)
 \move(0 3)\rlvec(1 0)
 \move(1 0)\rlvec(0 3)
 \move(0 -0.2)\rlvec(0 3.2)
 \move(0 0)\rlvec(1 1)
 \move(0 2)\rlvec(1 1)
 \htext(0.25 0.75){$1$}
 \htext(0.75 0.25){$0$}
 \htext(0.75 2.25){$3$}
 \htext(0.25 2.75){$4$}
 \htext(0.5 1.5){$2$}
\esegment
\move(-0.6 0)
\bsegment
\move(0.3 0)\rlvec(0.15 0)\rlvec(0 1)\rlvec(-0.15 0)\ifill f:0.5
\move(0.15 2)\rlvec(0.15 0)\rlvec(0 1)\rlvec(-0.15 0)\ifill f:0.5
 \move(0 -0.2)\rlvec(0 3.2)
 \move(0.3 0)\rlvec(0 2)
 \move(0.15 0)\rlvec(0 1)
 \move(0.45 0)\rlvec(0 1)
 \move(0 0)\rlvec(0.6 0)
 \move(0 1)\rlvec(0.3 0)
 \move(0.3 1)\rlvec(0.3 0)
 \move(0 2)\rlvec(0.3 0)
 \move(0 3)\rlvec(0.15 0)\rlvec(0 -1)
 \move(0.15 3)\rlvec(0.15 0)\rlvec(0 -1)
\esegment
\move(-1.2 0)
\bsegment
 \move(0 -0.2)\rlvec(0 3.2)
 \move(0.3 0)\rlvec(0 2)
 \move(0.15 0)\rlvec(0 1)
 \move(0.45 0)\rlvec(0 1)
 \move(0 0)\rlvec(0.6 0)
 \move(0 1)\rlvec(0.3 0)
 \move(0.45 1)\rlvec(0.15 0)
 \move(0 2)\rlvec(0.3 0)
 \move(0 3)\rlvec(0.15 0)\rlvec(0 -1)
\esegment
\move(-1.8 0)
\bsegment
 \move(0 -0.2)\rlvec(0 3.2)
 \move(0.3 0)\rlvec(0 2)
 \move(0.15 0)\rlvec(0 1)
 \move(0.45 0)\rlvec(0 1)
 \move(-0.4 0)\rlvec(1 0)
 \move(0 1)\rlvec(0.3 0)
 \move(0.45 1)\rlvec(0.15 0)
 \move(0 2)\rlvec(0.3 0)
 \move(0 3)\rlvec(0.15 0)\rlvec(0 -1)
\esegment
\esegment
\move(8 -12)%%%%%%%%%%%%%
\bsegment
\move(0 0)
 \bsegment
 \move(0 0)\rlvec(1 0)
 \move(0 1)\rlvec(1 0)
 \move(0 2)\rlvec(1 0)
 \move(0 3)\rlvec(1 0)
 \move(0 4)\rlvec(1 0)
 \move(1 0)\rlvec(0 4)
 \move(0 -0.2)\rlvec(0 4.2)
 \move(0 0)\rlvec(1 1)
 \move(0 2)\rlvec(1 1)
 \htext(0.25 0.75){$1$}
 \htext(0.75 0.25){$0$}
 \htext(0.75 2.25){$3$}
 \htext(0.25 2.75){$4$}
 \htext(0.5 1.5){$2$}
 \htext(0.5 3.5){$2$}
\esegment
\move(-0.6 0)
\bsegment
\move(0 3)\rlvec(0.3 0)\rlvec(0 1)\rlvec(-0.3 0)\ifill f:0.5
 \move(0 -0.2)\rlvec(0 3.2)
 \move(0.3 0)\rlvec(0 2)
 \move(0.15 0)\rlvec(0 1)
 \move(0.45 0)\rlvec(0 1)
 \move(0 0)\rlvec(0.6 0)
 \move(0 1)\rlvec(0.3 0)
 \move(0.45 1)\rlvec(0.15 0)
 \move(0 2)\rlvec(0.3 0)
 \move(0 3)\rlvec(0.15 0)\rlvec(0 -1)
 \move(0.15 3)\rlvec(0.15 0)\rlvec(0 -1)
 \move(0 3)\rlvec(0 1)\rlvec(0.3 0)\rlvec(0 -1)
\esegment
\move(-1.2 0)
\bsegment
 \move(0 -0.2)\rlvec(0 3.2)
 \move(0.3 0)\rlvec(0 2)
 \move(0.15 0)\rlvec(0 1)
 \move(0.45 0)\rlvec(0 1)
 \move(0 0)\rlvec(0.6 0)
 \move(0 1)\rlvec(0.3 0)
 \move(0.45 1)\rlvec(0.15 0)
 \move(0 2)\rlvec(0.3 0)
 \move(0 3)\rlvec(0.15 0)\rlvec(0 -1)
\esegment
\move(-1.8 0)
\bsegment
 \move(0 -0.2)\rlvec(0 3.2)
 \move(0.3 0)\rlvec(0 2)
 \move(0.15 0)\rlvec(0 1)
 \move(0.45 0)\rlvec(0 1)
 \move(-0.4 0)\rlvec(1 0)
 \move(0 1)\rlvec(0.3 0)
 \move(0.45 1)\rlvec(0.15 0)
 \move(0 2)\rlvec(0.3 0)
 \move(0 3)\rlvec(0.15 0)\rlvec(0 -1)
\esegment
\esegment
\move(-8 -12)%%%%%%%%%%%%
\bsegment
\move(-0.6 0)
\bsegment
\move(0.3 1)\rlvec(0.3 0)\rlvec(0 1)\rlvec(-0.3 0)\ifill f:0.5
 \move(0 -0.2)\rlvec(0 3.2)
 \move(0.3 0)\rlvec(0 2)
 \move(0.15 0)\rlvec(0 1)
 \move(0.45 0)\rlvec(0 1)
 \move(0 0)\rlvec(0.6 0)
 \move(0 1)\rlvec(0.3 0)
 \move(0.3 1)\rlvec(0.3 0)
 \move(0 2)\rlvec(0.6 0)
 \move(0 3)\rlvec(0.15 0)\rlvec(0 -1)
\esegment
\move(0 0)
 \bsegment
 \move(0 0)\rlvec(1 0)
 \move(0 1)\rlvec(1 0)
 \move(0 2)\rlvec(1 0)
 \move(0 3)\rlvec(1 0)
 \move(1 0)\rlvec(0 3)
 \move(0 -0.2)\rlvec(0 3.2)
 \move(0 0)\rlvec(1 1)
 \move(0 2)\rlvec(1 1)
 \htext(0.25 0.75){$1$}
 \htext(0.75 0.25){$0$}
 \htext(0.75 2.25){$3$}
 \htext(0.25 2.75){$4$}
 \htext(0.5 1.5){$2$}
\esegment
\move(-1.2 0)
\bsegment
 \move(0 -0.2)\rlvec(0 3.2)
 \move(0.3 0)\rlvec(0 2)
 \move(0.15 0)\rlvec(0 1)
 \move(0.45 0)\rlvec(0 1)
 \move(0 0)\rlvec(0.6 0)
 \move(0 1)\rlvec(0.3 0)
 \move(0.45 1)\rlvec(0.15 0)
 \move(0 2)\rlvec(0.3 0)
 \move(0 3)\rlvec(0.15 0)\rlvec(0 -1)
\esegment
\move(-1.8 0)
\bsegment
 \move(0 -0.2)\rlvec(0 3.2)
 \move(0.3 0)\rlvec(0 2)
 \move(0.15 0)\rlvec(0 1)
 \move(0.45 0)\rlvec(0 1)
 \move(-0.4 0)\rlvec(1 0)
 \move(0 1)\rlvec(0.3 0)
 \move(0.45 1)\rlvec(0.15 0)
 \move(0 2)\rlvec(0.3 0)
 \move(0 3)\rlvec(0.15 0)\rlvec(0 -1)
\esegment
\esegment
\move(-12 -19.5)%%%%%%%%%%%%
\bsegment
\move(-0.6 0)
\bsegment
 \move(0 -0.2)\rlvec(0 3.2)
 \move(0.3 0)\rlvec(0 2)
 \move(0.15 0)\rlvec(0 1)
 \move(0.45 0)\rlvec(0 1)
 \move(0 0)\rlvec(0.6 0)
 \move(0 1)\rlvec(0.3 0)
 \move(0.3 1)\rlvec(0.3 0)
 \move(0 2)\rlvec(0.6 0)
 \move(0 3)\rlvec(0.15 0)\rlvec(0 -1)
\esegment
\move(0 0)
 \bsegment
 \move(0 0)\rlvec(1 0)
 \move(0 1)\rlvec(1 0)
 \move(0 2)\rlvec(1 0)
 \move(0 3)\rlvec(1 0)
 \move(1 0)\rlvec(0 3)
 \move(0 -0.2)\rlvec(0 3.2)
 \move(0 0)\rlvec(1 1)
 \move(0 2)\rlvec(1 1)
 \htext(0.25 0.75){$1$}
 \htext(0.75 0.25){$0$}
 \htext(0.75 2.25){$3$}
 \htext(0.25 2.75){$4$}
 \htext(0.5 1.5){$2$}
\esegment
\move(-1.2 0)
\bsegment
\move(0.3 1)\rlvec(0.15 0)\rlvec(0 -1)\rlvec(-0.15 0)\ifill f:0.5
 \move(0 -0.2)\rlvec(0 3.2)
 \move(0.3 0)\rlvec(0 2)
 \move(0.15 0)\rlvec(0 1)
 \move(0.45 0)\rlvec(0 1)
 \move(0 0)\rlvec(0.6 0)
 \move(0 1)\rlvec(0.3 0)
 \move(0.3 1)\rlvec(0.3 0)
 \move(0 2)\rlvec(0.3 0)
 \move(0 3)\rlvec(0.15 0)\rlvec(0 -1)
\esegment
\move(-1.8 0)
\bsegment
 \move(0 -0.2)\rlvec(0 3.2)
 \move(0.3 0)\rlvec(0 2)
 \move(0.15 0)\rlvec(0 1)
 \move(0.45 0)\rlvec(0 1)
 \move(-0.4 0)\rlvec(1 0)
 \move(0 1)\rlvec(0.3 0)
 \move(0.45 1)\rlvec(0.15 0)
 \move(0 2)\rlvec(0.3 0)
 \move(0 3)\rlvec(0.15 0)\rlvec(0 -1)
\esegment
\esegment
\move(-8 -19.5)%%%%%%%%%%%%
\bsegment
\move(0 0)
 \bsegment
 \move(0 0)\rlvec(1 0)
 \move(0 1)\rlvec(1 0)
 \move(0 2)\rlvec(1 0)
 \move(0 3)\rlvec(1 0)
 \move(1 0)\rlvec(0 3)
 \move(0 -0.2)\rlvec(0 3.2)
 \move(0 0)\rlvec(1 1)
 \move(0 2)\rlvec(1 1)
 \htext(0.25 0.75){$1$}
 \htext(0.75 0.25){$0$}
 \htext(0.75 2.25){$3$}
 \htext(0.25 2.75){$4$}
 \htext(0.5 1.5){$2$}
\esegment
\move(-0.6 0)
\bsegment
\move(0.3 2)\rlvec(0 1)\rlvec(0.15 0)\rlvec(0 -1)\ifill f:0.5
 \move(0 -0.2)\rlvec(0 3.2)
 \move(0.3 0)\rlvec(0 2)
 \move(0.15 0)\rlvec(0 1)
 \move(0.45 0)\rlvec(0 1)
 \move(0 0)\rlvec(0.6 0)
 \move(0 1)\rlvec(0.3 0)
 \move(0.3 1)\rlvec(0.3 0)
 \move(0 2)\rlvec(0.6 0)
 \move(0 3)\rlvec(0.15 0)\rlvec(0 -1)
 \move(0.3 2)\rlvec(0 1)\rlvec(0.15 0)\rlvec(0 -1)
\esegment
\move(-1.2 0)
\bsegment
 \move(0 -0.2)\rlvec(0 3.2)
 \move(0.3 0)\rlvec(0 2)
 \move(0.15 0)\rlvec(0 1)
 \move(0.45 0)\rlvec(0 1)
 \move(0 0)\rlvec(0.6 0)
 \move(0 1)\rlvec(0.3 0)
 \move(0.45 1)\rlvec(0.15 0)
 \move(0 2)\rlvec(0.3 0)
 \move(0 3)\rlvec(0.15 0)\rlvec(0 -1)
\esegment
\move(-1.8 0)
\bsegment
 \move(0 -0.2)\rlvec(0 3.2)
 \move(0.3 0)\rlvec(0 2)
 \move(0.15 0)\rlvec(0 1)
 \move(0.45 0)\rlvec(0 1)
 \move(-0.4 0)\rlvec(1 0)
 \move(0 1)\rlvec(0.3 0)
 \move(0.45 1)\rlvec(0.15 0)
 \move(0 2)\rlvec(0.3 0)
 \move(0 3)\rlvec(0.15 0)\rlvec(0 -1)
\esegment
\esegment
\move(-4 -19.5)%%%%%%%%%%%%
\bsegment
\move(-0.6 0)
\bsegment
\move(0 2)\rlvec(0.6 0)\rlvec(0 0.5)\rlvec(-0.6 0)\ifill f:0.5
 \move(0 -0.2)\rlvec(0 2.7)
 \move(0.3 0)\rlvec(0 2)
 \move(0.15 0)\rlvec(0 1)
 \move(0.45 0)\rlvec(0 1)
 \move(0 0)\rlvec(0.6 0)
 \move(0 1)\rlvec(0.3 0)
 \move(0.3 1)\rlvec(0.3 0)
 \move(0 2)\rlvec(0.6 0)
 \move(0.15 2.5)\rlvec(0 -0.5)
 \move(0.3 2.5)\rlvec(0 -0.5)
 \move(0.45 2.5)\rlvec(0 -0.5)
 \lpatt(0.02 0.1)
 \move(0 2.5)\rlvec(0.6 0)
\esegment
\move(0 0)
 \bsegment
 \move(0 0)\rlvec(1 0)
 \move(0 1)\rlvec(1 0)
 \move(0 2)\rlvec(1 0)
 \move(0 3)\rlvec(1 0)
 \move(1 0)\rlvec(0 3)
 \move(0 -0.2)\rlvec(0 3.2)
 \move(0 0)\rlvec(1 1)
 \move(0 2)\rlvec(1 1)
 \htext(0.25 0.75){$1$}
 \htext(0.75 0.25){$0$}
 \htext(0.75 2.25){$3$}
 \htext(0.25 2.75){$4$}
 \htext(0.5 1.5){$2$}
\esegment
\move(-1.2 0)
\bsegment
 \move(0 -0.2)\rlvec(0 3.2)
 \move(0.3 0)\rlvec(0 2)
 \move(0.15 0)\rlvec(0 1)
 \move(0.45 0)\rlvec(0 1)
 \move(0 0)\rlvec(0.6 0)
 \move(0 1)\rlvec(0.3 0)
 \move(0.45 1)\rlvec(0.15 0)
 \move(0 2)\rlvec(0.3 0)
 \move(0 3)\rlvec(0.15 0)\rlvec(0 -1)
\esegment
\move(-1.8 0)
\bsegment
 \move(0 -0.2)\rlvec(0 3.2)
 \move(0.3 0)\rlvec(0 2)
 \move(0.15 0)\rlvec(0 1)
 \move(0.45 0)\rlvec(0 1)
 \move(-0.4 0)\rlvec(1 0)
 \move(0 1)\rlvec(0.3 0)
 \move(0.45 1)\rlvec(0.15 0)
 \move(0 2)\rlvec(0.3 0)
 \move(0 3)\rlvec(0.15 0)\rlvec(0 -1)
\esegment
\esegment
\move(0 -19.5)%%%%%%%%%%%
\bsegment
\move(-0.6 0)
\bsegment
\move(0.3 1.5)\rlvec(0 -0.5)\rlvec(0.3 0)\rlvec(0 0.5)\ifill f:0.5
\move(0 3.5)\rlvec(0 -0.5)\rlvec(0.3 0)\rlvec(0 0.5)\ifill f:0.5
 \move(0 -0.2)\rlvec(0 3.2)
 \move(0.3 0)\rlvec(0 2)
 \move(0.15 0)\rlvec(0 1)
 \move(0.45 0)\rlvec(0 1)
 \move(0 0)\rlvec(0.6 0)
 \move(0 1)\rlvec(0.3 0)
 \move(0.3 1)\rlvec(0.3 0)
 \move(0 2)\rlvec(0.3 0)
 \move(0 3)\rlvec(0.15 0)\rlvec(0 -1)
 \move(0.15 3)\rlvec(0.15 0)\rlvec(0 -1)
 \move(0 3)\rlvec(0 0.5)
 \move(0.3 3)\rlvec(0 0.5)
 \lpatt(0.02 0.1)
 \move(0 3.5)\rlvec(0.3 0)
 \move(0.3 1.5)\rlvec(0.3 0)
\esegment
\move(0 0)
 \bsegment
 \move(0 0)\rlvec(1 0)
 \move(0 1)\rlvec(1 0)
 \move(0 2)\rlvec(1 0)
 \move(0 3)\rlvec(1 0)
 \move(0 4)\rlvec(1 0)
 \move(1 0)\rlvec(0 4)
 \move(0 -0.2)\rlvec(0 4.2)
 \move(0 0)\rlvec(1 1)
 \move(0 2)\rlvec(1 1)
 \htext(0.25 0.75){$1$}
 \htext(0.75 0.25){$0$}
 \htext(0.75 2.25){$3$}
 \htext(0.25 2.75){$4$}
 \htext(0.5 1.5){$2$}
 \htext(0.5 3.5){$2$}
\esegment
\move(-1.2 0)
\bsegment
 \move(0 -0.2)\rlvec(0 3.2)
 \move(0.3 0)\rlvec(0 2)
 \move(0.15 0)\rlvec(0 1)
 \move(0.45 0)\rlvec(0 1)
 \move(0 0)\rlvec(0.6 0)
 \move(0 1)\rlvec(0.3 0)
 \move(0.45 1)\rlvec(0.15 0)
 \move(0 2)\rlvec(0.3 0)
 \move(0 3)\rlvec(0.15 0)\rlvec(0 -1)
\esegment
\move(-1.8 0)
\bsegment
 \move(0 -0.2)\rlvec(0 3.2)
 \move(0.3 0)\rlvec(0 2)
 \move(0.15 0)\rlvec(0 1)
 \move(0.45 0)\rlvec(0 1)
 \move(-0.4 0)\rlvec(1 0)
 \move(0 1)\rlvec(0.3 0)
 \move(0.45 1)\rlvec(0.15 0)
 \move(0 2)\rlvec(0.3 0)
 \move(0 3)\rlvec(0.15 0)\rlvec(0 -1)
\esegment
\esegment
\move(4 -19.5)%%%%%%%%%%%%%
\bsegment
\move(0 0)
 \bsegment
 \move(0 0)\rlvec(1 0)
 \move(0 1)\rlvec(1 0)
 \move(0 2)\rlvec(1 0)
 \move(0 3)\rlvec(1 0)
 \move(0 4)\rlvec(1 0)
 \move(0 5)\rlvec(1 0)
 \move(1 0)\rlvec(0 5)
 \move(0 -0.2)\rlvec(0 5.2)
 \move(0 0)\rlvec(1 1)
 \move(0 2)\rlvec(1 1)
 \move(0 4)\rlvec(1 1)
 \htext(0.25 0.75){$1$}
 \htext(0.75 0.25){$0$}
 \htext(0.75 2.25){$3$}
 \htext(0.25 2.75){$4$}
 \htext(0.5 1.5){$2$}
 \htext(0.5 3.5){$2$}
 \htext(0.25 4.75){$1$}
 \htext(0.75 4.25){$0$}
\esegment
\move(-0.6 0)
\bsegment
 \move(0 4)\rlvec(0 0.5)\rlvec(0.3 0)\rlvec(0 -0.5)\ifill f:0.5
 \move(0.3 0)\rlvec(0 0.5)\rlvec(0.3 0)\rlvec(0 -0.5)\ifill f:0.5
 \move(0 -0.2)\rlvec(0 3.2)
 \move(0.3 0)\rlvec(0 2)
 \move(0.15 0)\rlvec(0 1)
 \move(0.45 0)\rlvec(0 0.5)
 \move(0 0)\rlvec(0.6 0)
 \move(0 1)\rlvec(0.3 0)
 \move(0 2)\rlvec(0.3 0)
 \move(0 3)\rlvec(0.15 0)\rlvec(0 -1)
 \move(0.15 3)\rlvec(0.15 0)\rlvec(0 -1)
 \move(0 3)\rlvec(0 1)\rlvec(0.3 0)\rlvec(0 -1)
 \move(0 4)\rlvec(0 0.5)
 \move(0.3 4)\rlvec(0 0.5)
 \move(0.15 4)\rlvec(0 0.5)
 \lpatt(0.02 0.1)
 \move(0 4.5)\rlvec(0.3 0)
 \move(0.3 0.5)\rlvec(0.3 0)
\esegment
\move(-1.2 0)
\bsegment
 \move(0 -0.2)\rlvec(0 3.2)
 \move(0.3 0)\rlvec(0 2)
 \move(0.15 0)\rlvec(0 1)
 \move(0.45 0)\rlvec(0 1)
 \move(0 0)\rlvec(0.6 0)
 \move(0 1)\rlvec(0.3 0)
 \move(0.45 1)\rlvec(0.15 0)
 \move(0 2)\rlvec(0.3 0)
 \move(0 3)\rlvec(0.15 0)\rlvec(0 -1)
\esegment
\move(-1.8 0)
\bsegment
 \move(0 -0.2)\rlvec(0 3.2)
 \move(0.3 0)\rlvec(0 2)
 \move(0.15 0)\rlvec(0 1)
 \move(0.45 0)\rlvec(0 1)
 \move(-0.4 0)\rlvec(1 0)
 \move(0 1)\rlvec(0.3 0)
 \move(0.45 1)\rlvec(0.15 0)
 \move(0 2)\rlvec(0.3 0)
 \move(0 3)\rlvec(0.15 0)\rlvec(0 -1)
\esegment
\esegment
\move(8 -19.5)%%%%%%%%%%%%%
\bsegment
\move(0 0)
 \bsegment
 \move(0 0)\rlvec(1 0)
 \move(0 1)\rlvec(1 0)
 \move(0 2)\rlvec(1 0)
 \move(0 3)\rlvec(1 0)
 \move(0 4)\rlvec(1 0)
 \move(0 5)\rlvec(1 0)
 \move(1 0)\rlvec(0 5)
 \move(0 -0.2)\rlvec(0 5.2)
 \move(0 0)\rlvec(1 1)
 \move(0 4)\rlvec(1 1)
 \move(0 2)\rlvec(1 1)
 \htext(0.25 0.75){$1$}
 \htext(0.75 0.25){$0$}
 \htext(0.75 2.25){$3$}
 \htext(0.25 2.75){$4$}
 \htext(0.5 1.5){$2$}
 \htext(0.5 3.5){$2$}
 \htext(0.25 4.75){$1$}
 \htext(0.75 4.25){$0$}
\esegment
\move(-0.6 0)
\bsegment
 \move(0.15 4)\rlvec(0 1)\rlvec(0.15 0)\rlvec(0 -1)\ifill f:0.5
 \move(0 -0.2)\rlvec(0 3.2)
 \move(0.3 0)\rlvec(0 2)
 \move(0.15 0)\rlvec(0 1)
 \move(0.45 0)\rlvec(0 1)
 \move(0 0)\rlvec(0.6 0)
 \move(0 1)\rlvec(0.3 0)
 \move(0.45 1)\rlvec(0.15 0)
 \move(0 2)\rlvec(0.3 0)
 \move(0 3)\rlvec(0.15 0)\rlvec(0 -1)
 \move(0.15 3)\rlvec(0.15 0)\rlvec(0 -1)
 \move(0 3)\rlvec(0 1)\rlvec(0.3 0)\rlvec(0 -1)
 \move(0.15 4)\rlvec(0 1)\rlvec(0.15 0)\rlvec(0 -1)
\esegment
\move(-1.2 0)
\bsegment
 \move(0 -0.2)\rlvec(0 3.2)
 \move(0.3 0)\rlvec(0 2)
 \move(0.15 0)\rlvec(0 1)
 \move(0.45 0)\rlvec(0 1)
 \move(0 0)\rlvec(0.6 0)
 \move(0 1)\rlvec(0.3 0)
 \move(0.45 1)\rlvec(0.15 0)
 \move(0 2)\rlvec(0.3 0)
 \move(0 3)\rlvec(0.15 0)\rlvec(0 -1)
\esegment
\move(-1.8 0)
\bsegment
 \move(0 -0.2)\rlvec(0 3.2)
 \move(0.3 0)\rlvec(0 2)
 \move(0.15 0)\rlvec(0 1)
 \move(0.45 0)\rlvec(0 1)
 \move(-0.4 0)\rlvec(1 0)
 \move(0 1)\rlvec(0.3 0)
 \move(0.45 1)\rlvec(0.15 0)
 \move(0 2)\rlvec(0.3 0)
 \move(0 3)\rlvec(0.15 0)\rlvec(0 -1)
\esegment
\esegment
\move(12 -19.5)%%%%%%%%%%%%%
\bsegment
\move(0 0)
 \bsegment
 \move(0 0)\rlvec(1 0)
 \move(0 1)\rlvec(1 0)
 \move(0 2)\rlvec(1 0)
 \move(0 3)\rlvec(1 0)
 \move(0 4)\rlvec(1 0)
 \move(1 0)\rlvec(0 4)
 \move(0 -0.2)\rlvec(0 4.2)
 \move(0 0)\rlvec(1 1)
 \move(0 2)\rlvec(1 1)
 \htext(0.25 0.75){$1$}
 \htext(0.75 0.25){$0$}
 \htext(0.75 2.25){$3$}
 \htext(0.25 2.75){$4$}
 \htext(0.5 1.5){$2$}
 \htext(0.5 3.5){$2$}
\esegment
\move(-0.6 0)
\bsegment
 \move(0 -0.2)\rlvec(0 3.2)
 \move(0.3 0)\rlvec(0 2)
 \move(0.15 0)\rlvec(0 1)
 \move(0.45 0)\rlvec(0 1)
 \move(0 0)\rlvec(0.6 0)
 \move(0 1)\rlvec(0.3 0)
 \move(0.45 1)\rlvec(0.15 0)
 \move(0 2)\rlvec(0.3 0)
 \move(0 3)\rlvec(0.15 0)\rlvec(0 -1)
 \move(0.15 3)\rlvec(0.15 0)\rlvec(0 -1)
 \move(0 3)\rlvec(0 1)\rlvec(0.3 0)\rlvec(0 -1)
\esegment
\move(-1.2 0)
\bsegment
\move(0.15 2)\rlvec(0.15 0)\rlvec(0 1)\rlvec(-0.15 0)\ifill f:0.5
 \move(0 -0.2)\rlvec(0 3.2)
 \move(0.3 0)\rlvec(0 2)
 \move(0.15 0)\rlvec(0 1)
 \move(0.45 0)\rlvec(0 1)
 \move(0 0)\rlvec(0.6 0)
 \move(0 1)\rlvec(0.3 0)
 \move(0.45 1)\rlvec(0.15 0)
 \move(0 2)\rlvec(0.3 0)
 \move(0 3)\rlvec(0.15 0)\rlvec(0 -1)
 \move(0.15 3)\rlvec(0.15 0)\rlvec(0 -1)
\esegment
\move(-1.8 0)
\bsegment
 \move(0 -0.2)\rlvec(0 3.2)
 \move(0.3 0)\rlvec(0 2)
 \move(0.15 0)\rlvec(0 1)
 \move(0.45 0)\rlvec(0 1)
 \move(-0.4 0)\rlvec(1 0)
 \move(0 1)\rlvec(0.3 0)
 \move(0.45 1)\rlvec(0.15 0)
 \move(0 2)\rlvec(0.3 0)
 \move(0 3)\rlvec(0.15 0)\rlvec(0 -1)
\esegment
\esegment
\esegment

\move(0 0)
\bsegment
\setsegscale 0.8
\move(-1.5 -1.3)\ravec(-1.3 -1.8)
\htext(-2.5 -1.9){$0$}
\move(1.5 -1.3)\ravec(1.3 -1.8)
\htext(2.5 -1.9){$4$}
\move(-6.3 -8.8)\ravec(-1.3 -1.8)
\htext(-7.2 -9.4){$2$}
\move(-3.3 -8.8)\ravec(1.3 -1.8)
\htext(-2.3 -9.4){$4$}
\move(6.3 -8.8)\ravec(1.3 -1.8)
\htext(7.2 -9.4){$2$}
\move(3.3 -8.8)\ravec(-1.3 -1.8)
\htext(2.3 -9.4){$0$}
\move(-11.7 -16.8)\ravec(-1.3 -1.8)
\htext(-12.7 -17.5){$1$}
\move(-10.1 -16.8)\ravec(0 -2)
\htext(-9.6 -17.5){$3$}
\move(-8.5 -16.8)\ravec(1.3 -1.8)
\htext(-7.5 -17.5){$4$}
\move(0 -16.8)\ravec(0 -2)
\htext(0.4 -17.5){$2$}
\move(11.7 -16.8)\ravec(1.3 -1.8)
\htext(12.8 -17.5){$3$}
\move(10.1 -16.8)\ravec(0 -1)
\htext(10.6 -17.5){$1$}
\move(8.5 -16.8)\ravec(-1.3 -1.8)
\htext(7.5 -17.5){$0$}
\vtext(-15.2 -27){$\cdots$}
\vtext(-10.2 -27){$\cdots$}
\vtext(-5.2 -27){$\cdots$}
\vtext(-0.2 -27){$\cdots$}
\vtext(4.8 -27){$\cdots$}
\vtext(9.8 -27){$\cdots$}
\vtext(14.8 -27){$\cdots$}
\esegment
\end{texdraw}%
\end{center}

\vspace{5mm}

\end{example}

%%%%%%%%%%%%%%%%%%%%%%%%%%%%%%%%%%%%%%%%%%%%%%%%%%%%%%%%%%%%%%%%%%%%%%%%%%
\bibliographystyle{amsplain}
%\bibliography{dnone}

\def\cprime{$'$}
\providecommand{\bysame}{\leavevmode\hbox
to3em{\hrulefill}\thinspace}
\providecommand{\MR}{\relax\ifhmode\unskip\space\fi MR }
% \MRhref is called by the amsart/book/proc definition of \MR.
\providecommand{\MRhref}[2]{%
  \href{http://www.ams.org/mathscinet-getitem?mr=#1}{#2}
} \providecommand{\href}[2]{#2}

%%%%%%%%%%%%%%%%%%%%%%%%%%%%%%%%%%%%%%%%%%%%%%%%%%%%%%%%%%%%%%%%%%%%%%%%%%
\end{document}